\documentclass[pdflatex,sn-mathphys-num,iicol]{sn-jnl}
\usepackage{graphicx}%
\usepackage{multirow}%
\usepackage{amsmath,amssymb,amsfonts}%
\usepackage{amsthm}%
\usepackage[title]{appendix}%
\usepackage{xcolor}%
\usepackage{textcomp}%
\usepackage{manyfoot}%
\usepackage{booktabs}%
\usepackage{algorithm}%
\usepackage{algorithmicx}%
\usepackage{algpseudocode}%
\usepackage{listings}%
\hypersetup{
  hidelinks,
  bookmarksopen=true,
  bookmarksnumbered=true
}

\theoremstyle{thmstyleone}%
\theoremstyle{thmstyletwo}%

\theoremstyle{thmstylethree}%

\begin{document}

\title[Fractional Calculus in Optimal Control and Game Theory]{Fractional Calculus in Optimal Control and Game Theory: Theory, Numerics, and Applications – A Survey}

\author*[1]{\fnm{Navid} \sur{Mojahed}}\email{nmojahed@ucdavis.edu}
\author[2]{\fnm{Hooman} \sur{Fatoorehchi}}\email{hoofa@kt.dtu.dk}
\author[1]{\fnm{Shima} \sur{Nazari}}\email{snazari@ucdavis.edu}

\affil*[1]{\orgdiv{Department of Mechanical and Aerospace Engineering}, \orgname{University of California, Davis}, \orgaddress{\city{Davis}, \state{CA}, \postcode{95616}, \country{USA}}}

\affil[2]{\orgdiv{Center for Energy Resources Engineering (CERE), Department of Chemical and Biochemical Engineering}, \orgname{Technical University of Denmark}, \orgaddress{\city{Kgs. Lyngby}, \postcode{2800}, \country{Denmark}}}

\abstract{Many physical, biological, and engineered systems display memory that challenges Markovian models. Fractional calculus offers nonlocal operators that capture such hereditary effects. This survey connects modeling, analysis, and controller/game design for systems with memory.

On modeling, we unify notation across Caputo, Riemann--Liouville, and Grünwald--Letnikov derivatives and relate them to engineering approximations based on diffusive (sum-of-exponentials) state augmentation and frequency-domain realizations such as the Oustaloup family of fractional-order approximations. On analysis, we synthesize fractional extensions of the calculus of variations and the Pontryagin maximum principle, and we review dynamic-programming formulations with memory, including path-dependent Hamilton--Jacobi--Bellman for optimal control and Hamilton--Jacobi--Isaacs for zero-sum games. On design, we cover single-agent controllers—linear--quadratic regulator, model predictive control, and fractional-order proportional--integral--derivative—and multi-agent formulations such as fractional differential games with Nash, Stackelberg (leader--follower), and minimax equilibria.

On computation, we compare time-domain schemes (including the L1 and Grünwald--Letnikov methods), frequency-domain approximations (e.g., Oustaloup-type realizations), and diffusive state augmentations, emphasizing accuracy–complexity trade-offs and practical remedies for the “curse of history” using windowing, sum-of-exponentials, and band-limited or low-rank representations. We conclude with applications in viscoelasticity, anomalous transport, electrochemistry, and control of engineered systems; a benchmark template for reproducibility; and open problems on existence and uniqueness of equilibria with memory, refined Isaacs-type conditions, stability under constraints, and scalable fractional-order model predictive control and multi-agent solvers.
}


\maketitle
\section{Introduction}
\label{sec:intro}
% \deleted{Many physical, biological, and engineered systems exhibit \emph{memory} and hereditary effects that challenge Markovian modeling. Fractional calculus (FC) provides parsimonious nonlocal operators to encode such effects in a first-principles manner, typically through power-law kernels and hereditary integrals. Classical references established the analytical foundations of FC—including the Riemann--Liouville (RL), Caputo, and Grünwald--Letnikov (GL) operators, their transform properties, admissible initial and terminal conditions, and well-posedness caveats \cite{Oldham1974,Samko1993,MillerRoss1993,Podlubny1999,Kilbas2006}. Subsequent developments across physics and engineering further clarified how memory and power-law attenuation inform continuum and field models \cite{Caputo1967,CaputoMainardi1971,Herrmann2014,Tarasov2011}.}
Many physical, biological, and engineered systems exhibit memory and hereditary effects that challenge purely Markovian modeling. Fractional calculus (FC) provides parsimonious nonlocal operators to encode such effects in a first-principles manner, typically through power-law kernels and hereditary integrals. Several definitions of fractional-order derivatives and integrals coexist, each tailored to different theoretical and practical needs. The Riemann–Liouville fractional derivative is historically fundamental and analytically convenient, but its requirement for fractional-order initial conditions often limits its physical interpretability. In contrast, the Caputo fractional derivative allows the use of standard integer-order initial conditions, making it more suitable for modeling real physical and engineering systems, albeit at the cost of slightly increased analytical complexity. The Grünwald–Letnikov derivative provides a natural link to discrete-time formulations and numerical implementation, although it is computationally demanding due to slow convergence. Other definitions such as the Weyl fractional derivative are mainly applicable to problems defined over infinite or periodic domains. Classical references established the analytical foundations of FC such as the Riemann–Liouville, Caputo, and Grünwald–Letnikov operators, their transform properties, admissible initial and terminal conditions, and well-posedness caveats \cite{Oldham1974,Samko1993,MillerRoss1993,Podlubny1999,Kilbas2006}. Subsequent developments across physics and engineering further clarified how memory and power-law attenuation inform continuum and field models \cite{Caputo1967,CaputoMainardi1971,Herrmann2014,Tarasov2011}. Consequently, the choice of definition is application-dependent, balancing physical realizability, analytical convenience, and numerical efficiency.
A compact chronology of these developments appears in Table~\ref{tab:timeline}.

Figure~\ref{fig:graphical_abstract} provides a graphical overview of the main modeling, analysis, and control/game-theoretic components covered in this survey.
\begin{figure*}[!t]
    \centering
    \includegraphics[width=\textwidth]{Graphical_Abstract}
    \caption{Graphical abstract of the survey. The figure summarizes the overall structure of the paper, linking modeling and taxonomy of fractional-order systems, optimality principles and numerical schemes, and their applications in control and game-theoretic settings, together with benchmark problems and directions for future research.}
    \label{fig:graphical_abstract}
\end{figure*}

Subsequent treatises established rigorous existence, uniqueness, and stability results for fractional differential equations (FDEs). In this theory, the Mittag--Leffler function plays the role of the exponential kernel for fractional-order (FO) dynamics and viscoelastic waves \cite{Diethelm2010,Mainardi2010}. Stability tests tailored to FO systems—such as spectral-sector conditions—and their control implications are now standard \cite{Matignon1996,Monje2010}. On the numerical side, convolution--quadrature and L1/GL-type time discretizations yield systematic approximations of the Caputo and RL operators \cite{Lubich1986,DiethelmFord2002}; further accuracy and stability refinements are available for subdiffusion regimes \cite{Alikhanov2015,JinLazarovZhou2016}, and modern monographs provide comprehensive overviews \cite{LiZeng2015}.

These analytical and numerical advances underpin contemporary applications of fractional models in domains where memory is intrinsic: linear and nonlinear viscoelasticity and rheology \cite{BagleyTorvik1983,Koeller1984,Mainardi2010}, anomalous transport and diffusion \cite{MetzlerKlafter2000}, and electrochemical systems in which constant-phase elements (CPEs) and impedance spectroscopy naturally motivate fractional operators \cite{BarsoukovMacdonald2005}. 
In control-oriented settings, fractional-order models often match measured dynamics with fewer states than high-order integer embeddings, which can aid identifiability and support more robust downstream controller design—such as the linear--quadratic regulator (LQR), model predictive control (MPC), and fractional-order proportional--integral--derivative (FOPID)—while also introducing distinct challenges related to history dependence, computational cost, and well-posed initialization \cite{Monje2010,Petras2011,Oustaloup2000}.

\textbf{Why FC for control and games?}
In optimal control (OC), high-fidelity modeling of history-dependent dynamics is often forced to rely on high-dimensional integer-order embeddings or ad-hoc state augmentations to emulate memory (e.g., transport or diffusion lags, viscoelastic or electrochemical hysteresis). Fractional-order state-space models encode memory directly through nonlocal derivatives with power-law kernels, and they frequently achieve comparable or superior fit with fewer states and more interpretable parameters \cite{Monje2010,Petras2011}. From a mathematical standpoint, the FO framework brings a mature toolbox of special functions and stability results—most notably Mittag--Leffler responses and sector-based stability criteria—that extend exponential behaviors and spectral conditions to memory-rich dynamics \cite{Gorenflo2014,Mainardi2010,Matignon1996}. In multi-agent settings, the same memory terms support path-dependent payoffs and long-range couplings, making fractional modeling relevant for differential games in which strategic interactions depend on cumulative effects and delayed responses.

On the computational side, three complementary families of numerical methods support analysis and controller synthesis for FO models.
(i) \emph{Time-domain discretizations} approximate the Caputo and RL operators using convolution--quadrature and L1/GL-type schemes. They offer tight error and stability guarantees in subdiffusion regimes and practical predictor--corrector variants for simulation and design \cite{Lubich1986,DiethelmFord2002,Alikhanov2015,JinLazarovZhou2016,LiZeng2015}.
(ii) \emph{Frequency-domain approximations}—for example, Oustaloup’s recursive filters within the CRONE (Commande Robuste d'Ordre Non Entier) methodology—realize $s^\alpha$ behavior over prescribed bands and integrate smoothly with classical loop-shaping and robust-design tools \cite{Oustaloup2000,Sabatier2015,GuglielmiHairer2025}.
(iii) \emph{Diffusive/state-augmentation realizations} approximate power-law kernels via Prony or relaxation spectra, enabling standard state-space synthesis (LQR/MPC) while controlling accuracy--complexity trade-offs through the number and placement of modes \cite{BagleyTorvik1983,Koeller1984,Monje2010}.
Together, these mechanisms mitigate the “curse of history” by using windowing and compact representations (low-rank or band-limited), which we will compare systematically.

In \emph{control design}, FO formulations span regulation, tracking, and robust performance. 
FO--LQR and FO tracking are typically built on augmented representations; FO--MPC employs explicit memory truncation \cite{CAO2024490,DEAZEREDO2024424} and terminal constraints; and FOPID leverages additional phase and gain freedom to achieve robustness margins that can be out of reach for integer-order PID in some plants \cite{Monje2010,Podlubny1999PID,Sabatier2015,Jayaram2024FOPID,Wendimu2025SciRep}. 
These designs surface recurring issues around initialization (Caputo vs.\ RL conditions), constraint handling under memory, and certification of closed-loop stability—all of which motivate unified benchmarks.

Beyond single-agent OC, \emph{game-theoretic control} models multi-agent interaction—competitive or cooperative—through equilibrium concepts such as Nash, Stackelberg, and mean-field limits, grounded in classical dynamic-game theory \cite{BasarOlsder1999,Isaacs1965}. When the dynamics and/or costs are fractional, the value functionals become history-dependent and the Hamilton--Jacobi--Bellman (HJB) and Hamilton--Jacobi--Isaacs (HJI) equations inherit nonlocal terms. Fractional extensions of the calculus of variations and Pontryagin’s maximum principle (PMP) provide necessary conditions and transversality relations for FO dynamics \cite{Agrawal2004,Pooseh2015}; in competitive or robust settings, these couple with Isaacs-type conditions that may require refined hypotheses in the presence of memory. In practice, fractional differential games (FDGs) are tackled computationally via diffusive augmentation or frequency-domain approximations, enabling alternating best-response, bilevel, or minimax solvers with controlled approximation error—an approach we survey and systematize in the sequel.

\paragraph*{Positioning and scope}
This review targets the triad of fractional calculus, optimal control, and game theory. Our contributions are:
\begin{itemize}[leftmargin=*]
  \item \textbf{Unified taxonomy:} We organize problem classes by (a) derivative definition (Caputo, RL, GL, diffusive/state-augmentation); (b) decision paradigm (Pontryagin maximum principle (PMP) vs.\ dynamic programming (DP)—HJB for control and HJI for games—vs.\ direct transcription); (c) system type (linear time-invariant (LTI), linear time-varying (LTV), nonlinear, stochastic); (d) constraints (state/input; path/terminal); (e) interaction structure (single- vs.\ multi-agent); and (f) game class (zero-sum, general-sum, Stackelberg, mean-field) with information structure (open- vs.\ closed-loop) \cite{Oldham1974,Samko1993,Podlubny1999,Kilbas2006,Monje2010,BasarOlsder1999,GOMOYUNOV2024335}.
  \item \textbf{Notation and assumptions:} We standardize symbols and initial and terminal conditions across FO variants, and clarify Caputo–RL conversions with implications for initial data and adjoints \cite{Podlubny1999,Kilbas2006,Diethelm2010}.
  \item \textbf{Optimality principles:} We consolidate the fractional calculus of variations, PMP for FO systems, and dynamic-programming viewpoints under memory (HJB/HJI), highlighting transversality and nonlocal adjoints \cite{Agrawal2004,Pooseh2015,Oldham1974}.
  \item \textbf{Canonical designs (single- and multi-agent):} We review FO--LQR for regulation and tracking, robust/minimax formulations, fractional differential games (FDGs)—including linear--quadratic (LQ) zero-sum and general-sum cases and Stackelberg structures—FO--MPC with prediction and stability under memory, and FOPID from an optimization/OC perspective \cite{Monje2010,Podlubny1999PID,BasarOlsder1999,Isaacs1965}.
  \item \textbf{Numerics and complexity:} We compare time-domain (L1/GL, predictor--corrector), frequency-domain (Oustaloup), and diffusive/state-augmentation methods, discussing accuracy, stability, memory windowing, and computational complexity \cite{Lubich1986,DiethelmFord2002,Oustaloup2000,BagleyTorvik1983,Koeller1984}.
  \item \textbf{Benchmarks and reproducibility:} We propose benchmark templates (problem instances, kernel parameterization, solver settings, window lengths) and a reporting checklist for like-for-like comparison, including FDGs.
  \item \textbf{Open problems:} We identify gaps in (i) existence and uniqueness of equilibria with memory and refined Isaacs-type conditions; (ii) controllability and observability notions aligned with OC and games; (iii) stability under constraints; (iv) scalable FO--MPC and multi-agent solvers; and (v) integrated data-driven identification and learning with FO control and game design \cite{Diethelm2010,Mainardi2010,Monje2010,BasarOlsder1999}.
\end{itemize}

\subsection*{Positioning with respect to prior surveys}
Classical monographs and handbooks laid the foundations and early control viewpoints
(e.g., \cite{Podlubny1999,Monje2010,Mainardi2010,Diethelm2010,Kilbas2006,Samko1993}),
while more recent surveys emphasize numerical methods and good practices \cite{Garrappa2018,Diethelm2020_Practices}
and control/filtering design \cite{Sabatier2015}.
Against this backdrop, our focus is threefold:
(i) a unified treatment of \emph{optimal control} (OC) and \emph{game-theoretic control} for systems with fractional dynamics;
(ii) a reproducibility-oriented numerics section—covering diffusive, CRONE, and L1 families—tailored to OC/MPC/game solvers; and
(iii) a benchmark template and reporting checklist to enable apples-to-apples (fair, like-for-like) comparisons across fractional optimal-control methods.

\begin{table}[ht]
\centering
\caption{Selected milestones shaping fractional calculus in control.}\label{tab:timeline}
\begin{tabular}{@{}p{1.2cm}p{4.9cm}p{1.8cm}@{}}
\toprule
Year & Milestone & Refs \\
\midrule
1967 & Caputo derivative for hereditary media & \cite{Caputo1967} \\
1971 & Caputo--Mainardi dissipation with memory & \cite{CaputoMainardi1971} \\
1996 & Sector stability (Matignon) for FO--LTI systems & \cite{Matignon1996} \\
1999 & Podlubny monograph; FOPID formalization & \cite{Podlubny1999,Podlubny1999PID} \\
2000 & Oustaloup recursive filter ($s^\alpha$ over bands) & \cite{Oustaloup2000} \\
2010 & Monographs consolidating analysis and control & \cite{Monje2010,Diethelm2010,Mainardi2010} \\
2015 & High-order L1 refinements for subdiffusion & \cite{Alikhanov2015} \\
2018--20 & Best practices and surveys in numerics & \cite{Garrappa2018,Diethelm2020_Practices} \\
2020--24 & FO HJB/HJI and path-dependent theory & \cite{Gomoyunov2020SICON,Gomoyunov2022COCV,Gomoyunov2024JDE} \\
\bottomrule
\end{tabular}
\end{table}

\paragraph*{Research questions (RQs)}

\textbf{RQ1:} What optimality principles (PMP/HJB/HJI) are available for FO dynamics, and under what assumptions?\\
\textbf{RQ2:} Which numerical representations (L1/GL, CRONE, diffusive) are most effective for OC/MPC/FDGs in practice?\\
\textbf{RQ3:} What are the open gaps in well-posedness under constraints, scalability, and reproducibility?

\paragraph*{Organization.}
Section~\ref{sec:Prelimin} fixes notation and basic tools of fractional calculus (Caputo, Riemann--Liouville, Grünwald--Letnikov), recalls the Volterra integral form, and states well-posedness and Mittag--Leffler decay facts used later. 
Section~\ref{sec:fods} formalizes fractional dynamical systems—commensurate/incommensurate linear models, the Matignon spectral sector for stability, BIBO notions, and controllability/observability (PBH/Gramian), with notes on tempered and time-varying cases.
Section~\ref{sec:Optim} develops optimality principles in the fractional setting: Euler--Lagrange identities and transversality (with right-sided operators), Noether-type invariance, the fractional PMP for Caputo/RL dynamics (adjoint and boundary terms), free-terminal-time variants, and dynamic programming/HJB/HJBI with memory.
Section~\ref{sec:canonical} turns these ingredients into working LQ formulations: FO–LQR/tracking/regulation via (A) augmented realizations that admit Riccati equations and (B) direct transcription as history-coupled QPs.
Section~\ref{sec:fopid-via-oc} recasts FOPID/CRONE design as optimal control, contrasting time-domain tuning with frequency-domain loop shaping and discussing robustness and digital realization.
Section~\ref{sec:mpc-fo} treats MPC for FO plants: prediction models (convolution/L1–GL–CQ, diffusive/Prony augmentation, Oustaloup), windowing and SOE compression, and stability/feasibility via terminal ingredients and tube ideas.
Section~\ref{sec:numerics} surveys numerics: CQ, L1/L2-$1\sigma$, predictor–corrector schemes and their fast variants; band-limited rational approximations; diffusive realizations; and practical guidance on accuracy, stability, and cost.
Section~\ref{sec:app} sketches representative applications (electrochemical CPE models, viscoelastic damping, mechatronics/robotics, biomedical systems) and typical FO control workflows.
Section~\ref{sec:bench} proposes a small benchmark suite (B1–B6) with shared metrics and a reporting checklist to enable apples-to-apples comparisons.
Section~\ref{sec:openprob} collects open problems (O1–O8) spanning theory, scalable computation, and reproducibility, and Section~\ref{sec:conc} concludes. 
An auxiliary proof appears in Appendix~\ref{app:cauchy-proof}.

\section{Preliminaries on Fractional Calculus}
\label{sec:Prelimin}
\subsection{Definitions and Properties}
\label{subsec:frac-defs}
\paragraph*{Notation and function spaces.}
We work on a finite time horizon $[0,T]$ with $T>0$, a fractional order $\alpha>0$, and $n:=\lceil \alpha\rceil\in\mathbb{N}$. For readability, most formulas are stated for $0<\alpha<1$ and extend in the standard way to $n-1<\alpha<n$. We write $\mathrm{AC}^n[0,T]$ for the set of functions whose $(n-1)$th derivative is absolutely continuous on $[0,T]$, and we denote by $\Gamma(\cdot)$ the Gamma function. Classical references for the material in this subsection include \cite{Samko1993,Podlubny1999,Kilbas2006,Diethelm2010,Oldham1974,MillerRoss1993}.

\paragraph*{Left fractional integrals (Riemann--Liouville).}
For $\alpha>0$ and $f\in L^1(0,T)$, the left Riemann--Liouville (RL) fractional integral on $[0,T]$ is
\begin{equation}
(I_{0+}^{\alpha} f)(t) \;=\; \frac{1}{\Gamma(\alpha)} \int_{0}^{t} (t-\tau)^{\alpha-1} f(\tau)\, d\tau,
\quad t\in[0,T].
\label{eq:RL-int}
\end{equation}
This construction generalizes the classical Cauchy formula for $n$-fold integration: an $n$-times iterated integral over $[a,x]$ can be written as a single integral (see Appendix~\ref{app:cauchy-proof}):
\begin{equation}
\begin{split}
(I^n f)(x)
&= \int_a^x \int_a^{\sigma_1} \cdots \int_a^{\sigma_{n-1}}
   f(\sigma_n)\, d\sigma_n \cdots d\sigma_2\, d\sigma_1 \\
&= \frac{1}{(n-1)!}\int_a^x (x-t)^{\,n-1} f(t)\,dt.
\end{split}
\end{equation}
Replacing the integer $n$ by a real order $\beta>0$ and $(n-1)!$ by the Gamma function $\Gamma(\beta)$
yields the \emph{Riemann--Liouville fractional integral} on $[a,x]$,
\[
(I_{a+}^{\beta} f)(x) \;=\; \frac{1}{\Gamma(\beta)} \int_a^x (x-t)^{\beta-1} f(t)\,dt,
\quad \beta>0,
\]
which specializes to \eqref{eq:RL-int} on $[0,T]$; see, e.g.,
\cite{Oldham1974,Samko1993,Podlubny1999,Kilbas2006,Diethelm2010}.
Fractional integrals form a semigroup:
\begin{equation}
I_{0+}^{\alpha} I_{0+}^{\beta} f \;=\; I_{0+}^{\alpha+\beta} f \qquad (\alpha,\beta>0),
\label{eq:int-semigroup}
\end{equation}
see \cite[Sec.~2.1]{Podlubny1999}, \cite[Thm.~2.1]{Kilbas2006}.

\paragraph*{Left fractional derivatives.}
Let $n=\lceil \alpha\rceil$ with $n-1<\alpha<n$. We use the following left-sided Riemann--Liouville and Caputo definitions on $[0,T]$.

\emph{Riemann--Liouville (RL):}
\begin{equation}
\begin{aligned}
({}^{\mathrm{RL}}\!D_{0+}^{\alpha} f)(t)
&= \frac{d^n}{dt^n}\Big(I_{0+}^{\,n-\alpha} f\Big)(t),
\end{aligned}
\end{equation}
\noindent for $0<\alpha<1$,
\begin{equation}
\begin{aligned}
({}^{\mathrm{RL}}\!D_{0+}^{\alpha} f)(t)
&= \frac{d}{dt}\!\left[\frac{1}{\Gamma(1-\alpha)}\!
    \int_0^t f(\tau)\,(t-\tau)^{-\alpha}\, d\tau\right].
\end{aligned}
\end{equation}

\emph{Caputo (C):}
\begin{equation}
\begin{aligned}
({}^{\mathrm{C}}\!D_{0+}^{\alpha} f)(t)
&= \Big(I_{0+}^{\,n-\alpha} f^{(n)}\Big)(t),
\end{aligned}
\end{equation}
\noindent for $0<\alpha<1$,
\begin{equation}
\begin{aligned}
({}^{\mathrm{C}}\!D_{0+}^{\alpha} f)(t)
&= \frac{1}{\Gamma(1-\alpha)}\!
   \int_{0}^{t} f'(\tau)\,(t-\tau)^{-\alpha}\, d\tau .
\end{aligned}
\end{equation}

A control-relevant distinction is that $({}^{\mathrm{C}} D_{0+}^{\alpha} c)=0$ for any constant $c$,
whereas $({}^{\mathrm{RL}} D_{0+}^{\alpha} c)= \dfrac{c}{\Gamma(1-\alpha)}\, t^{-\alpha}$ for $0<\alpha<1$ \cite[Sec.~2.4]{Podlubny1999}.
Equivalently, for $n-1<\alpha<n$,
\begin{equation}
\begin{aligned}
({}^{\mathrm{C}} D_{0+}^{\alpha} f)(t)
&= ({}^{\mathrm{RL}} D_{0+}^{\alpha})\!\Bigg[
    f(t) - \sum_{k=0}^{n-1}\frac{f^{(k)}(0^+)}{k!}\,t^{k}\Bigg].
\end{aligned}
\end{equation}

\paragraph*{Gr\"unwald--Letnikov (GL) derivative.}
For sufficiently regular $f$ (e.g., $f\in \mathrm{AC}^n[0,T]$) and $t\in[0,T]$,
\begin{equation*}
\begin{aligned}
({}^{\mathrm{GL}}\!D_{0+}^{\alpha} f)(t)
&= \lim_{h\to 0^+} h^{-\alpha}
   \sum_{k=0}^{\lfloor t/h\rfloor} (-1)^k \binom{\alpha}{k}\, f(t-kh),\\[-1mm]
\binom{\alpha}{k}
&= \frac{\Gamma(\alpha+1)}{\Gamma(k+1)\,\Gamma(\alpha-k+1)} \, .
\end{aligned}
\end{equation*}
Under mild assumptions, ${}^{\mathrm{GL}}\!D_{0+}^{\alpha}$ coincides with the left Riemann--Liouville derivative ${}^{\mathrm{RL}}\!D_{0+}^{\alpha}$ \cite[Sec.~2.3]{Podlubny1999}, \cite[Sec.~8.1]{Samko1993}.

\paragraph*{Right-sided operators on $[0,T]$.}
Right-sided operators are essential for adjoint equations and terminal conditions in control \cite{Kilbas2006,Pooseh2015}. The right Riemann--Liouville (RL) integral and derivative are
\begin{equation}
\begin{aligned}
(I_{T-}^{\alpha} f)(t)
&= \frac{1}{\Gamma(\alpha)} \int_{t}^{T} (\tau-t)^{\alpha-1} f(\tau)\, d\tau,\\[-1mm]
({}^{\mathrm{RL}}\!D_{T-}^{\alpha} f)(t)
&= (-1)^n \frac{d^n}{dt^n}\Big( I_{T-}^{\,n-\alpha} f \Big)(t),
\end{aligned}
\end{equation}
and the right Caputo derivative is
\begin{equation}
\begin{aligned}
({}^{\mathrm{C}}\!D_{T-}^{\alpha} f)(t)
&= (-1)^n\, I_{T-}^{\,n-\alpha}\!\big( f^{(n)} \big)(t).
\end{aligned}
\end{equation}
See \cite[Sec.~2.3]{Kilbas2006} for details.

For the common case $0<\alpha<1$ (i.e., $n=1$), these specialize to
\begin{equation}
\begin{aligned}
({}^{\mathrm{RL}}\!D_{T-}^{\alpha} f)(t)
&= -\,\frac{d}{dt}\!\left[\frac{1}{\Gamma(1-\alpha)}\!
    \int_{t}^{T} (\tau-t)^{-\alpha}\, f(\tau)\, d\tau\right],\\[-1mm]
({}^{\mathrm{C}}\!D_{T-}^{\alpha} f)(t)
&= -\,\frac{1}{\Gamma(1-\alpha)}\!
    \int_{t}^{T} (\tau-t)^{-\alpha}\, f'(\tau)\, d\tau .
\end{aligned}
\end{equation}

\noindent \emph{Control note.} As on the left side, the Caputo derivative of a constant is zero, whereas the RL derivative is nonzero for $0<\alpha<1$ and scales like $(T-t)^{-\alpha}$; this makes Caputo conditions natural in terminal-value formulations.

\paragraph*{Laplace transforms (left-sided).}
Let $n=\lceil\alpha\rceil$ with $n-1<\alpha\le n$, and denote $F(s)=\mathcal{L}\{f\}(s)$. Then
\begin{align}
\mathcal{L}\{\,{}^{\mathrm{C}}\!D_{0+}^{\alpha} f\,\}(s)
&= s^{\alpha} F(s)
   - \sum_{k=0}^{n-1} s^{\alpha-1-k}\, f^{(k)}(0^+),
\label{eq:Lap-Caputo}\\[2mm]
\mathcal{L}\{\,{}^{\mathrm{RL}}\!D_{0+}^{\alpha} f\,\}(s)
&= s^{\alpha} F(s)
   - \sum_{k=0}^{n-1} s^{\,n-1-k}\,
     \Big(I_{0+}^{\,n-\alpha} f\Big)^{(k)}\!(0^+),
\label{eq:Lap-RL}
\end{align}
see \cite[Eq.~(2.4.44)]{Podlubny1999}, \cite[Sec.~2.5]{Kilbas2006}.
For $0<\alpha<1$, \eqref{eq:Lap-Caputo} simplifies to $s^{\alpha}F(s)-s^{\alpha-1}f(0^+)$.

\noindent\textit{Remark.}
For closed-form time responses of high-order rational $G(s)$ (possibly with non-integer exponents),
one can recast the inverse Laplace task as an (integer or fractional) differential equation
and solve it via the Adomian decomposition method; see \cite{FatoorehchiRach2020_AEJ}.
This is especially helpful when partial-fraction or residue-based inversions become impractical.

\paragraph*{Equivalence relations and initial data.}
For $n-1<\alpha<n$ and $f\in\mathrm{AC}^n[0,T]$,
\begin{equation}
\label{eq:Caputo-RL-equivalence}
({}^{\mathrm{C}}\! D_{0+}^{\alpha} f)(t)
= ({}^{\mathrm{RL}}\! D_{0+}^{\alpha})\!\Bigg[
  f(t) - \sum_{k=0}^{n-1} \frac{f^{(k)}(0^+)}{k!}\, t^{k}
\Bigg].
\end{equation}
Thus, Caputo dynamics use standard \emph{integer-order} initial conditions $f^{(k)}(0^+)$, while RL dynamics naturally involve \emph{fractional} initial data—namely the values of $I_{0+}^{\,n-\alpha}f$ and its integer derivatives at $0^+$.
See \cite[Sec.~2.4]{Podlubny1999} and \cite[Ch.~2]{Kilbas2006}.
\paragraph*{Semigroup and calculus rules.}
While fractional integrals satisfy the semigroup property in \eqref{eq:int-semigroup}, fractional derivatives generally do \emph{not}:
\begin{equation*}
({}^{\mathrm{C}}\! D_{0+}^{\alpha})\big( {}^{\mathrm{C}}\! D_{0+}^{\beta} f \big)
\;\neq\; {}^{\mathrm{C}}\! D_{0+}^{\alpha+\beta} f .
\end{equation*}
Equality can hold under suitable \emph{zero-state initialization}; for Caputo derivatives this requires $f^{(k)}(0^+)=0$ for $k=0,\dots,n-1$, and for RL derivatives it requires vanishing fractional initial data (the values of $I_{0+}^{\,n-\alpha}f$ and its integer derivatives at $0^+$). See \cite[Sec.~2.6]{Diethelm2010} and \cite[Sec.~2.4]{Podlubny1999}. 

Product and chain rules typically expand into infinite series (generalized Leibniz/chain formulas) and do not close into simple finite forms \cite{Samko1993}. Transform-domain tools and special functions are surveyed in \cite{Gorenflo2014,Mainardi2010}, and initialization pitfalls and semigroup caveats are discussed in \cite{Diethelm2010,Diethelm2020_Practices}.

\paragraph*{Fractional integration by parts (IBP).}
A key identity for variational analysis and adjoint equations on $[0,T]$ (for $0<\alpha<1$ and sufficiently smooth $f,g$) is

\begin{equation}
\label{eq:IBP}
\begin{split}
\int_{0}^{T} g(t)\, ({}^{\mathrm{C}}\! D_{0+}^{\alpha} f)(t)\, dt
&= \Big[\, (I_{T-}^{\,1-\alpha} g)(t)\, f(t) \,\Big]_{0}^{T}
  \\&+ \int_{0}^{T} f(t)\, ({}^{\mathrm{RL}}\! D_{T-}^{\alpha} g)(t)\, dt.
\end{split}
\end{equation}
Here $I_{T-}^{\,1-\alpha}$ and ${}^{\mathrm{RL}}\!D_{T-}^{\alpha}$ are the right-sided operators on $[0,T]$; see \cite[Sec.~2.4]{Kilbas2006}, \cite[Ch.~2]{Pooseh2015}. Identity \eqref{eq:IBP} underpins fractional Pontryagin-type necessary conditions: the adjoint typically carries a \emph{right-sided} fractional derivative.

\paragraph*{Remarks for control.}
In practice, Caputo models are preferred for initial-value problems because they admit classical initial conditions and yield the simple Laplace form in \eqref{eq:Lap-Caputo}. Riemann--Liouville models arise naturally in hereditary laws and in boundary- or terminal-value formulations. In both cases, initialization and terminal terms must be handled carefully to ensure well-posed problems and to derive the correct adjoint and transversality conditions—hence the use of right-sided operators for adjoint equations \cite{Podlubny1999,Kilbas2006,Pooseh2015}.

\subsection{Equivalent Representations}
\label{subsec:eq-rep}
Fractional operators admit several engineering-friendly representations that preserve the essential memory and phase behavior while enabling standard analysis and synthesis. We focus on three workhorses: (i) frequency-domain rational approximations (Oustaloup/CRONE), (ii) diffusive realizations based on Prony or relaxation spectra, and (iii) augmented-state embeddings that convert hereditary convolutions into low-order ordinary differential equations (ODEs). We will use these representations repeatedly and compare their accuracy–complexity trade-offs in later sections; see \cite{Oustaloup2000,Monje2010,Sabatier2015,BagleyTorvik1983,Koeller1984,Mainardi2010,LiZeng2015,gomoyunov2024lipschitzcontinuityresultsminimax}.

\paragraph*{(i) Frequency-domain (Oustaloup/CRONE) rational approximation.}
Over a prescribed band $[\omega_\ell,\omega_h]$, the operator $s^{\alpha}$ ($\alpha\in\mathbb{R}$) is approximated by a stable, proper, minimum-phase rational transfer function with $2N{+}1$ first-order real poles and zeros geometrically placed on $[\omega_\ell,\omega_h]$:
\begin{equation}
\label{eq:oustaloup}
s^{\alpha} \;\approx\; K \prod_{k=-N}^{N} \frac{1 + s/\omega_{z,k}}{1 + s/\omega_{p,k}},
\qquad \omega_{z,k},\omega_{p,k}\in[\omega_\ell,\omega_h].
\end{equation}
The standard Oustaloup recursion selects $(\omega_{z,k},\omega_{p,k})$ so that the Bode magnitude has
an approximately constant slope of $20\,\alpha$ dB/dec and a nearly constant phase of $\alpha\cdot 90^\circ$
over $[\omega_\ell,\omega_h]$, while $K$ matches a chosen reference (typically the geometric mean)
\cite{Oustaloup2000,Monje2010,Sabatier2015}.
Increasing $N$ widens the accurate band and reduces phase ripple at the cost of filter order.
This CRONE-style approximation integrates smoothly with loop-shaping and robust control, but bandwidth
selection, passivity (when required), and numerical scaling still need care.

\paragraph*{(ii) Diffusive (Prony/relaxation) representation.}
A convenient way to work with memory is to rewrite the power-law kernel as a weighted continuum of exponentials (“diffusive” representation). For $0<\alpha<1$,
\begin{equation}
\label{eq:caputo-diffusive}
\begin{aligned}
({}^{\mathrm{C}}\! D_{0+}^{\alpha} x)(t)
&= \frac{1}{\Gamma(1-\alpha)} \int_{0}^{t} \frac{\dot{x}(\tau)}{(t-\tau)^{\alpha}}\, d\tau \\
&= \int_{0}^{\infty} \mu_{\alpha}(\lambda)\, z(\lambda,t)\, d\lambda .
\end{aligned}
\end{equation}
Here the kernel admits an exponential-mixture representation
\begin{equation}
\label{eq:kernel-mixture}
\begin{aligned}
\frac{1}{\Gamma(1-\alpha)}(t-\tau)^{-\alpha}
&= \int_{0}^{\infty} \mu_{\alpha}(\lambda)\, e^{-\lambda (t-\tau)}\, d\lambda, \\
\mu_{\alpha}(\lambda)
&= \frac{1}{\Gamma(\alpha)\,\Gamma(1-\alpha)}\, \lambda^{\alpha-1},
\end{aligned}
\end{equation}
and the \emph{diffusive state} $z(\lambda,t)$ evolves as a first-order relaxation driven by $\dot{x}$:
\begin{equation*}
\frac{\partial}{\partial t}\, z(\lambda,t) \;=\; -\lambda\, z(\lambda,t) + \dot{x}(t),
\qquad
z(\lambda,0)=0 \;\; (\text{Caputo}).
\end{equation*}
Approximating the integral in \eqref{eq:caputo-diffusive} by a positive quadrature
$\sum_{i=1}^{M} w_i\, z(\lambda_i,t)$ yields a finite-dimensional ODE cascade (a \emph{Prony series})
that preserves causality and—when $w_i>0$—the passivity properties of many viscoelastic/electrochemical
models \cite{BagleyTorvik1983,Koeller1984,Mainardi2010}. The nodes $\{\lambda_i\}$ are typically
log-spaced on $[\lambda_{\min},\lambda_{\max}]$ to cover the target band, with $M$ controlling
the accuracy--complexity trade-off.

\paragraph*{(iii) Augmented-state embeddings for control.}
Let $y(t)=({}^{\mathrm{C}}\!D_{0+}^{\alpha} x)(t)$ be required in a state-space design. Using the diffusive approximation,
\begin{equation}
\label{eq:augmented-caputo}
\begin{aligned}
y(t) &\approx \sum_{i=1}^{M} w_i\, z_i(t),\\[-1mm]
\dot{z}_i(t) &= -\lambda_i\, z_i(t) + \dot{x}(t), \qquad z_i(0)=0,\quad i=1,\dots,M .
\end{aligned}
\end{equation}
For a fractional LTI plant with $({}^{\mathrm{C}}\!D_{0+}^{\alpha} x)(t)=A\,x(t)+B\,u(t)$, \eqref{eq:augmented-caputo} yields the **diffusive DAE**:
\begin{equation}
\label{eq:augmented-dae}
\begin{aligned}
\dot{x}(t) &= v(t), \qquad \text{(introduce $v:=\dot{x}$)}\\
\dot{z}_i(t) &= -\lambda_i\, z_i(t) + v(t), \quad i=1,\dots,M,\\
\sum_{i=1}^{M} w_i\, z_i(t) &= A\,x(t) + B\,u(t).
\end{aligned}
\end{equation}
This index-1 DAE can be handled directly in constrained LQR/MPC/differential-game solvers. If a **pure ODE** is desired, differentiate the algebraic equation once (assuming $\dot u$ is available or generated via a prefilter):
\begin{equation}
\label{eq:augmented-ode}
\Big(\textstyle\sum_{i=1}^{M} w_i\Big) v(t)\;-\;\sum_{i=1}^{M} w_i\,\lambda_i\, z_i(t)\;=\;A\,v(t)\;+\;B\,\dot{u}(t),
\end{equation}
which, when $((\sum_i w_i) I - A)$ is nonsingular on the design region, allows solving for $v$ and yields an ODE in the augmented state $(x,z_1,\dots,z_M)$. In frequency-domain workflows, the Oustaloup approximation \eqref{eq:oustaloup} directly provides a low-order rational model suitable for robust loop-shaping \cite{Oustaloup2000,Sabatier2015,Monje2010}. Compared with direct time-convolution, augmented embeddings reduce memory/runtime from $\mathcal{O}(N_t)$ history storage to $\mathcal{O}(M)$ states, with accuracy controlled by $M$ (or the filter order $2N{+}1$) \cite{LiZeng2015,Monje2010}.

\paragraph*{Design notes.}
\emph{Band selection.} Choose $[\omega_\ell,\omega_h]$ (or $[\lambda_{\min},\lambda_{\max}]$) to cover the closed-loop bandwidth and dominant disturbance spectrum; outside this band, the $s^{\alpha}$ behavior is not guaranteed and approximation error can grow quickly.
\emph{Initialization.} Caputo vs.\ RL definitions imply different starts: zero diffusive states correspond to Caputo initial-value problems (IVPs), whereas RL models require fractional initial data.
\emph{Passivity and positive-real (PR) constraints.} These can be enforced via positive measures in \eqref{eq:kernel-mixture} (or positive quadrature weights) and, when needed by the application, through appropriate CRONE designs \cite{Sabatier2015,Mainardi2010}.
Implementation and tuning guidelines for CRONE and diffusive realizations are summarized in \cite{Sabatier2015,Monje2010}, with software-oriented tutorials in \cite{Garrappa2018}.

\section{Fractional Dynamical Systems}
\label{sec:fods}
This section collects state-space forms for fractional-order (FO) systems, states basic well-posedness conditions, reviews stability notions based on Mittag--Leffler functions, and sets up controllability/observability preliminaries most relevant to control and games. Classical sources include \cite{Podlubny1999,Kilbas2006,Diethelm2010,Mainardi2010}; control-oriented overviews appear in \cite{Monje2010,Petras2011}.

\subsection{State-space forms}
\label{subsec:ss-forms}

We consider left-sided Caputo dynamics on $[0,T]$. For $0<\alpha<1$, a \emph{commensurate} FO LTI model is
\begin{equation}
\label{eq:comm-LTI}
{}^{\mathrm{C}}\! D_{0+}^{\alpha} x(t) \;=\; A\,x(t) + B\,u(t),
\qquad
y(t)\;=\; C\,x(t)+D\,u(t),
\end{equation}
with the classical initial condition $x(0)=x_0$. The corresponding continuous-time transfer function is
\begin{equation}
\label{eq:tf-comm}
G(s) \;=\; C\,\big(s^{\alpha} I - A\big)^{-1} B + D,
\end{equation}
see \cite[Sec.~2.5]{Podlubny1999}. For the zero-input case $u \equiv 0$, the continuous-time solution of \eqref{eq:comm-LTI} is
\begin{equation}
\label{eq:mittag-comm}
x(t) \;=\; E_{\alpha}\!\big(A\,t^{\alpha}\big)\,x_0,
\qquad
E_{\alpha}(z)=\sum_{k=0}^{\infty}\frac{z^{k}}{\Gamma(\alpha k + 1)},
\end{equation}
where $E_{\alpha}$ is the (matrix) Mittag--Leffler operator defined via its convergent power-series expansion, in direct analogy with $e^{At} = \sum_{k=0}^{\infty} A^{k} t^{k}/k!$ for integer-order LTI systems.

With a nonzero input $u$, the variation-of-constants formula involves the two-parameter Mittag--Leffler function
\begin{equation}
\label{eq:ML-def}
E_{\alpha,\beta}(z)\;=\;\sum_{k=0}^{\infty}\frac{z^{k}}{\Gamma(\alpha k+\beta)},
\qquad \alpha>0,\;\beta>0,
\end{equation}
and reads
\begin{equation}
\label{eq:voce}
\begin{aligned}
x(t) \;=\; &\, E_{\alpha}\!\big(A t^{\alpha}\big)\,x_0 \\[-1mm]
& + \int_{0}^{t} (t-\tau)^{\alpha-1}
   E_{\alpha,\alpha}\!\big(A (t-\tau)^{\alpha}\big)\,
   B\,u(\tau)\, d\tau .
\end{aligned}
\end{equation}

For \emph{incommensurate} models with possibly different state orders $\alpha_i\in(0,1]$,
\begin{equation}
\label{eq:incomm}
\operatorname{diag}\big({}^{\mathrm{C}}\! D_{0+}^{\alpha_1},\dots,{}^{\mathrm{C}}\! D_{0+}^{\alpha_n}\big)\, x(t)
\;=\; A\,x(t)+B\,u(t).
\end{equation}
A common approach rationalizes the orders as $\alpha_i = p_i/m$ ($p_i,m\in\mathbb{N}$) and works with the characteristic equation in the auxiliary variable $w=s^{1/m}$; see \cite[Sec.~2.6]{Monje2010}, \cite[Ch.~2]{Petras2011}. Multi-term models
\[
\textstyle \sum_{j=1}^{q} E_j\, {}^{\mathrm{C}}\! D_{0+}^{\alpha_j} x(t)\;=\;A\,x(t)+B\,u(t)
\]
are treated similarly via state augmentation. Tempered variants (power-law kernels with exponential tempering) are also used in practice and admit analogous state-space treatments; see recent controllability/observability characterizations for tempered Caputo systems in \cite{Lamrani2024TemperedArxiv}.

\subsection{Existence, Uniqueness, and Continuous Dependence}
\label{subsec:existence}
We consider left-sided Caputo IVPs. For $0<\alpha<1$,
\begin{equation}
\label{eq:nonlin}
\begin{aligned}
({}^{\mathrm{C}}\! D_{0+}^{\alpha} x)(t) &= f\big(t,x(t),u(t)\big),\\
x(0) &= x_0 ,
\end{aligned}
\end{equation}
which is equivalent to a Volterra integral equation of the second kind:
\begin{equation}
\label{eq:volterra}
x(t)
= x_0 + \frac{1}{\Gamma(\alpha)} \int_{0}^{t} (t-\tau)^{\alpha-1}
     f\big(\tau,x(\tau),u(\tau)\big)\, d\tau .
\end{equation}

\noindent\emph{Well-posedness.}
If $f$ is continuous in $t$, locally Lipschitz in $x$ (uniformly on compact $t$-intervals), and the input
$u$ is measurable and essentially bounded on $[0,T]$, then \eqref{eq:nonlin} admits a unique solution on $[0,T]$
with continuous dependence on $(x_0,u)$ (via the equivalent integral form and a Banach fixed-point argument);
see \cite[Thm.~6.1]{Diethelm2010}, \cite[Ch.~3]{Kilbas2006}.
Global-in-time existence follows under standard growth bounds on $f$ (e.g., linear growth in $x$).
These statements extend to $n-1<\alpha<n$ with $n$ classical initial values by the usual reduction.

\subsection{Stability Notions: Mittag--Leffler and Spectral Sectors}
\label{subsec:stability}
For the commensurate linear system \eqref{eq:comm-LTI} with $0<\alpha<2$, \emph{asymptotic stability} (i.e., $x(t)\!\to\!0$ for $u\equiv0$) holds iff all eigenvalues of $A$ lie outside the sector
$\mathcal{S}_\alpha := \{\lambda:\,|\arg(\lambda)| \le \alpha\pi/2\}$, i.e., satisfy the Matignon condition
\begin{equation}
 |\arg(\lambda_i(A))| \;>\; \frac{\alpha\pi}{2}, \qquad i=1,\dots,n.
\label{eq:matignon}
\end{equation}
Equivalently, the spectrum avoids a sector of angle $\alpha\pi$ around the positive real axis \cite{Matignon1996}. Under \eqref{eq:matignon}, solutions decay \emph{algebraically} (Mittag--Leffler type) rather than exponentially; in particular,
$E_{\alpha}(\lambda t^{\alpha}) \sim -\sum_{k=1}^{\infty}\lambda^{-k}/\Gamma(1-\alpha k)$ as $t\to\infty$ \cite[Ch.~1]{Mainardi2010}.
For incommensurate or multi-term systems, one can rationalize the orders (e.g., $w=s^{1/m}$) and check the characteristic roots on the appropriate Riemann surface \cite[Sec.~2.6]{Monje2010,Petras2011}.
Recent refinements of Lyapunov/comparison methods for Caputo systems provide constructive certificates beyond spectral tests in nonautonomous settings; see, e.g., \cite{Ineh2024AIMSmath}.

\noindent\emph{Input--output stability.}
BIBO stability follows from bounded frequency response $G(j\omega)$ together with appropriate initialization choices \cite{Monje2010}.

\subsection{Controllability and Observability (LTI, Commensurate)}
\label{subsec:ctrl-obsv}
Let \eqref{eq:comm-LTI} hold with $0<\alpha<1$. The state-transition matrix is $\Phi_{\alpha}(t)=E_{\alpha}(A t^{\alpha})$. The reachable set on $[0,T]$ is
\begin{equation}
\label{eq:reach}
\begin{aligned}
\mathcal{R}(T)
= \Big\{\,& \int_{0}^{T} (T-\tau)^{\alpha-1}
      E_{\alpha,\alpha}\!\big(A(T-\tau)^{\alpha}\big)\,B\,u(\tau)\, d\tau \\
& \Bigm|\, u \in L^2 \Big\}.
\end{aligned}
\end{equation}

If the classical controllability matrix $\mathcal{C}=[B,\,AB,\,\dots,\,A^{n-1}B]$ has rank $n$, then $(A,B)$ is (completely) controllable in the FO sense and $\mathcal{R}(T)=\mathbb{R}^n$ for all $T>0$. Equivalently, a PBH-type test holds: $\mathrm{rank}[\lambda I-A\;\;B]=n$ for all eigenvalues $\lambda$ of $A$ \cite[Sec.~5.2]{Monje2010}, \cite[Ch.~3]{Petras2011}. A Gramian characterization uses
\begin{equation}
\label{eq:gramian-ctrl}
W_{\alpha}(T)=\int_{0}^{T} \Phi_{\alpha}(\tau)\, B B^{\!\top} \Phi_{\alpha}^{\!\top}(\tau)\, d\tau,
\end{equation}
whose positive definiteness implies reachability on $[0,T]$.

Observability admits dual statements for $(C,A)$ with the PBH test
\[
\mathrm{rank}\!\begin{bmatrix}\lambda I-A \\ C \end{bmatrix}\!=n
\quad \text{for all eigenvalues $\lambda$ of $A$},
\]
and an FO observability Gramian
\begin{equation}
\label{eq:gramian-obsv}
W_{\alpha,o}(T)
= \int_{0}^{T} \Phi_{\alpha}^{\!\top}(\tau)\, C^{\!\top} C\, \Phi_{\alpha}(\tau)\, d\tau,
\end{equation}
for which $W_{\alpha,o}(T)\succ 0$ (positive definite) is equivalent to complete observability on $[0,T]$; see \cite{Monje2010}.

%% NEW (recent generalizations):
Kalman/Gramian-type criteria extend to \emph{tempered} Caputo systems with explicit tests, and to certain time-varying settings; see \cite{Lamrani2024TemperedArxiv,Vishnukumar2024PhysScr} for recent formulations and examples.

\paragraph*{Notes for design and identification.}
(i) Algebraic (non-exponential) decay affects standard Lyapunov certificates; Lyapunov methods exist but require care with fractional derivatives \cite{Diethelm2010}. 
(ii) Initialization (Caputo vs.\ RL) and terminal terms impact adjoint equations in optimality principles. 
(iii) For implementation and synthesis, augmented-state embeddings (diffusive or Oustaloup) map FO plants to low-order ODEs; see Sec.~\ref{subsec:eq-rep}. 
%% NEW (pointer to fast diffusive/SOE):
Efficient SOE-based diffusive realizations for FO kernels have recently been proposed with “memoryless” implementations useful for control synthesis \cite{Guglielmi2025SOE}.

\section{Optimality Principles for FO Systems}
\label{sec:Optim}
\subsection{Fractional Calculus of Variations}
\label{subsec:fcv}

We collect necessary conditions for Bolza-type problems with fractional derivatives on $[0,T]$, focusing on left-sided Caputo and Riemann--Liouville (RL) operators. Background: \cite{Agrawal2004,Pooseh2015,Podlubny1999,Kilbas2006,Samko1993,Diethelm2010}.
\emph{Recent updates:} higher-order Caputo Euler--Lagrange (EL) conditions and sufficient criteria \cite{Ferreira2024HigherCaputo}, and duality/IBP formulations that avoid mixing left/right operators \cite{Torres2024Duality}.

\paragraph*{Setup (Caputo).}
Let $x:[0,T]\to\mathbb{R}^n$. Consider the Bolza functional
\begin{equation}
\label{eq:bolza-caputo}
\begin{aligned}
\mathcal{J}[x]
&= \Phi\!\big(x(T)\big)
  + \int_{0}^{T} L\!\big(t,\,x(t),\,v(t)\big)\, dt,\\
v(t) &= ({}^{\mathrm{C}}\!D_{0+}^{\alpha} x)(t),\; 0<\alpha<1 .
\end{aligned}
\end{equation}

with admissible boundary conditions (some components fixed, some free), $x\in\mathrm{AC}[0,T]$, $L$ continuously differentiable in $(x,v)$, and $\Phi$ differentiable. For first variations we take $x_\varepsilon=x+\varepsilon\,\eta$ with an admissible $\eta$ consistent with the endpoint conditions.

\paragraph*{Fractional Euler--Lagrange (Caputo).}
Using the fractional integration-by-parts identity on $[0,T]$,
\begin{equation}
\label{eq:ibp-fcv}
\begin{aligned}
\int_{0}^{T} g(t)\,\big({}^{\mathrm{C}}\!D_{0+}^{\alpha}\eta\big)(t)\,dt
&= \Big[\,(I_{T-}^{\,1-\alpha} g)(t)\,\eta(t)\,\Big]_{0}^{T}\\
&\quad + \int_{0}^{T} \eta(t)\,\big({}^{\mathrm{RL}}\!D_{T-}^{\alpha} g\big)(t)\,dt .
\end{aligned}
\end{equation}

(cf.\ \cite[Sec.~2.4]{Kilbas2006}, \cite[Ch.~2]{Pooseh2015}; see also the duality-based IBP perspective in \cite{Torres2024Duality}),
stationarity $\delta\mathcal{J}=0$ for all admissible $\eta$ with fixed endpoints yields the Euler--Lagrange equation
\begin{equation}
\label{eq:EL-caputo}
\frac{\partial L}{\partial x}(t,x,v)
\;+\;
{}^{\mathrm{RL}}\!D_{T-}^{\alpha}\!\left(\frac{\partial L}{\partial v}(t,x,v)\right)
\;=\; 0,
\qquad t\in(0,T).
\end{equation}
Here $\partial L/\partial x$ and $\partial L/\partial v$ are evaluated at $(t,\,x(t),\,v(t))$.
See, e.g., \cite{MalinowskaTorres2012,FredericoTorres2007,Agrawal2004,AlmeidaTorres2019}.

\paragraph*{Natural boundary (transversality) conditions.}
When some endpoint coordinates are free, the boundary term from \eqref{eq:ibp-fcv} together with $\delta\Phi$
yields the usual transversality conditions (Caputo case; write $L_v=\partial L/\partial v$).
Let $\mathcal{F}_T$ and $\mathcal{F}_0$ be the index sets of the coordinates that are free at $t=T$ and $t=0$, respectively. Then
\begin{equation}
\label{eq:transv-caputo-clean}
\begin{aligned}
\big[\;\Phi_x\!\big(x(T)\big) + (I_{T-}^{\,1-\alpha} L_v)(T)\;\big]_i &= 0, && i\in\mathcal{F}_T,\\[-1mm]
\big[\;(I_{T-}^{\,1-\alpha} L_v)(0)\;\big]_i &= 0, && i\in\mathcal{F}_0.
\end{aligned}
\end{equation}
If a coordinate is fixed at $t=0$ or $t=T$, its variation vanishes ($\eta_i(0)=0$ or $\eta_i(T)=0$) and the corresponding line in \eqref{eq:transv-caputo-clean} is omitted. 
Cf.\ \cite[Sec.~2.4]{Kilbas2006}, \cite[Ch.~2]{Pooseh2015}; see also \cite{Torres2024Duality}.

\paragraph*{RL-in-Lagrangian variant.}
If $L(t,x,w)$ depends on $w=({}^{\mathrm{RL}}\!D_{0+}^{\alpha}x)(t)$, the EL equation becomes
\begin{equation}
\label{eq:EL-RL}
\frac{\partial L}{\partial x}(t,x,w) \;+\; {}^{\mathrm{C}}\!D_{T-}^{\alpha}\!\left(\frac{\partial L}{\partial w}(t,x,w)\right) \;=\; 0,
\end{equation}
and the natural boundary terms involve $I_{T-}^{\,1-\alpha}\big(\partial L/\partial w\big)$ in the same fashion as
\eqref{eq:transv-caputo-clean} (see \cite[Sec.~2.4]{Podlubny1999}, \cite[Ch.~2]{Kilbas2006}).

\paragraph*{Multi-term / higher-order cases.}
For $L\big(t,x,\{v_j\}_{j=1}^q\big)$ with $v_j = {}^{\mathrm{C}}\!D_{0+}^{\alpha_j} x$ and $0<\alpha_1<\cdots<\alpha_q<1$, the Euler--Lagrange condition becomes
\begin{equation}
\label{eq:EL-multi}
\frac{\partial L}{\partial x}
\;+\;
\sum_{j=1}^{q} {}^{\mathrm{RL}}\!D_{T-}^{\alpha_j}\!
\left(\frac{\partial L}{\partial v_j}\right)
\;=\; 0.
\end{equation}
Analogous forms hold for $n-1<\alpha<n$ with $n$ classical initial values \cite{Kilbas2006}. See also \cite{Ferreira2024HigherCaputo} for higher-order Caputo formulations in the fractional calculus of variations.

\paragraph*{Constraints and sufficiency.}
An isoperimetric constraint $\int_0^T \phi(t,x,v)\,dt=c$ adds a constant multiplier to $L$. Path inequalities $g(t,x,u)\le 0$ introduce time-varying Karush–Kuhn–Tucker(KKT) multipliers with the usual nonnegativity and complementarity conditions \cite[Ch.~4]{Pooseh2015}. If $L$ is jointly convex in $(x,v)$ and $\Phi$ is convex, any solution of \eqref{eq:EL-caputo} is a global minimizer \cite[Ch.~3]{Pooseh2015}.

\paragraph*{Invariance (Noether-type).}
By invariance we mean that $L\,dt$ (and $\Phi$) are unchanged, up to a total derivative, under a smooth one-parameter transformation of $(t,x)$. In fractional problems, the associated Noether identity still holds but the momentum is memory-weighted: the conserved/balance quantity involves right-sided fractional operators acting on $L_v:=\partial L/\partial v$ (e.g., $I_{T-}^{\,1-\alpha}L_v$) rather than $L_v$ alone. Thus, symmetries yield first integrals with nonlocal weights, which recover the classical Noether law in the limit $\alpha\to 1$; see \cite{FredericoTorres2007,MalinowskaTorres2012,Riewe1996} and the surveys \cite{MalinowskaTorres2012,FredericoTorres2007,AlmeidaTorres2019}.

\subsection{Fractional Pontryagin Maximum Principle (PMP)}
\label{subsec:fpmp}
We state first-order necessary conditions for optimal control with fractional-order dynamics.
Throughout, $0<\alpha<1$, the state $x:[0,T]\to\mathbb{R}^n$, the control $u:[0,T]\to U\subset\mathbb{R}^m$
(measurable, with $U$ typically convex), and $f,\ell,\Phi$ are continuously differentiable in their arguments.
Classical references include \cite{Agrawal2004,Pooseh2015,Podlubny1999,Kilbas2006,Diethelm2010}.
Recent extensions cover delay/Volterra kernels, state constraints, and incommensurate or distributed-order models
\cite{Gasimov2024QTDS,Moon2025AIMSmath,NdairouTorres2023Math,Cruz2025AIMSmath,YusubovMahmudov2024ACS}.

\paragraph*{Problem statement (Caputo state equation).}
Consider the Bolza problem
\begin{equation}
\label{eq:fpmp-cost}
\min_{u(\cdot)} J
= \Phi\big(x(T)\big) + \int_{0}^{T} \ell\big(t,x(t),u(t)\big)\,dt,
\end{equation}
% \todo{negative space in Eq. (41) not so good!}
subject to the Caputo dynamics
\begin{equation}
\label{eq:fpmp-state}
({}^{\mathrm{C}}D_{0+}^{\alpha} x)(t) \;=\; f\big(t,x(t),u(t)\big), 
\qquad x(0)=x_0,
\end{equation}
and pointwise constraints $u(t)\in U$ a.e.\ on $[0,T]$.
Define the Hamiltonian
\begin{equation}
\label{eq:fpmp-H}
H\big(t,x,u,p\big) \;=\; \ell(t,x,u) \;+\; p^{\top} f(t,x,u).
\end{equation}
\paragraph*{Adjoint equation and transversality (Caputo case).}
The costate $p:[0,T]\to\mathbb{R}^n$ carries the backward (right-sided) memory induced by Caputo dynamics.
From the first variation together with the fractional IBP \eqref{eq:ibp-fcv}, the adjoint satisfies
\begin{equation}
\label{eq:fpmp-adj}
({}^{\mathrm{RL}}\!D_{T-}^{\alpha} p)(t)
= H_x\!\big(t, x^{\ast}(t), u^{\ast}(t), p(t)\big),
\qquad t\in(0,T).
\end{equation}
The terminal condition is
\begin{equation}
\label{eq:fpmp-trans}
(I_{T-}^{\,1-\alpha} p)(T) = \nabla_x \Phi\!\big(x^{\ast}(T)\big).
\end{equation}
If some components of $x(0)$ are free, they obey
\begin{equation*}
\big[(I_{T-}^{\,1-\alpha} p)(0)\big]_i = 0
\quad \text{for each free index } i;
\end{equation*}
when $x(0)$ is fixed, no condition at $t=0$ arises.

\paragraph*{Hamiltonian minimization condition.}
For a.e.\ $t\in[0,T]$,
\begin{equation}
\label{eq:fpmp-min}
u^{\ast}(t)\in\arg\min_{u\in U} H\big(t,x^{\ast}(t),u,p(t)\big).
\end{equation}
If $U$ is open and $\ell,f$ are $C^{1}$ in $u$, the first-order stationarity
$H_{u}\!\big(t,x^{\ast}(t),u^{\ast}(t),p(t)\big)=0$ is necessary.
For convex $U$, this is equivalent to the variational inequality
$\langle H_{u}(t,x^{\ast},u^{\ast},p),\,u-u^{\ast}(t)\rangle\!\ge\!0$ for all $u\in U$ a.e.

\paragraph*{Summary (Caputo PMP).}
The necessary conditions are:
(i) the state dynamics \eqref{eq:fpmp-state};
(ii) the adjoint equation \eqref{eq:fpmp-adj};
(iii) the terminal transversality \eqref{eq:fpmp-trans} (with initial-side conditions for any free components of $x(0)$); and
(iv) the Hamiltonian minimization \eqref{eq:fpmp-min},
together with $u(\cdot)\in\mathcal{U}$ (admissible controls) and $x(0)=x_0$.
Related dynamic-programming/HJB formulations for fractional, path-dependent dynamics include
\cite{Gomoyunov2020SICON,Gomoyunov2022COCV,Gomoyunov2024JDE}.
\paragraph*{Variant with RL state equation.}
If the state dynamics use the left Riemann--Liouville derivative,
\begin{equation}
\label{eq:fpmp-state-RL}
({}^{\mathrm{RL}}\!D_{0+}^{\alpha} x)(t)=f\big(t,x(t),u(t)\big),
\end{equation}
then the costate evolves with the right Caputo derivative,
\begin{equation}
\label{eq:fpmp-adj-RL}
({}^{\mathrm{C}}\!D_{T-}^{\alpha} p)(t)=H_x\big(t,x^{\ast}(t),u^{\ast}(t),p(t)\big).
\end{equation}
\paragraph*{Transversality.}
Let $\mathcal{F}_0$ denote the set of indices with free initial values. The boundary conditions read
\begin{equation}
\label{eq:fpmp-trans-rl}
\begin{aligned}
(I_{T-}^{\,1-\alpha}p)(T) &= \nabla_x\Phi\big(x^{\ast}(T)\big),\\
\big[(I_{T-}^{\,1-\alpha}p)(0)\big]_i &= 0,\quad i\in\mathcal{F}_0.
\end{aligned}
\end{equation}
See \cite{Kilbas2006,Podlubny1999,Agrawal2004}; for delay/Volterra-type memory, see \cite{Gasimov2024QTDS}.

\paragraph*{Path constraints and terminal equalities.}
For smooth path inequalities $g(t,x,u)\le 0$, introduce a nonnegative multiplier $\mu(t)$ and use the augmented Hamiltonian
\[
\tilde H(t,x,u,p)=H(t,x,u,p)+\mu(t)^{\top}g(t,x,u),
\]
with
\[
\begin{aligned}
\mu_i(t) &\ge 0,\\
g_i(t,x^{\ast},u^{\ast}) &\le 0,\\
\mu_i(t)\,g_i(t,x^{\ast},u^{\ast}) &= 0, \quad \text{a.e.\ } t\in[0,T].
\end{aligned}
\]

A terminal equality $r\big(x(T)\big)=0$ adds a multiplier $\lambda$ in the boundary condition,
\begin{equation}
\label{eq:fpmp-trans-constraints}
(I_{T-}^{\,1-\alpha}p)(T)
= \nabla_x\Phi\big(x^{\ast}(T)\big)
  + \nabla_x r\big(x^{\ast}(T)\big)^{\!\top}\lambda .
\end{equation}
See \cite[Ch.~4]{Pooseh2015}; state constraints in Caputo problems are treated in \cite{Moon2025AIMSmath}, and
incommensurate/distributed-order models in \cite{NdairouTorres2023Math,Cruz2025AIMSmath}.

\paragraph*{Free terminal time and extensions.}
If the terminal time $T$ is free, an additional transversality condition closes the system. For Caputo dynamics, under standard regularity assumptions,
\begin{equation}
\label{eq:free-time-trans}
H\!\big(T,x^{\ast}(T),u^{\ast}(T),p(T)\big)
\;+\;
\frac{\partial \Phi}{\partial T}\!\big(x^{\ast}(T)\big)
= 0,
\end{equation}
see \cite{Agrawal2004,Pooseh2015}. Multi-term models are handled by summing the corresponding right-sided derivatives in the adjoint. Discrete-structure FO systems (e.g., fractional Roesser) also admit PMP statements \cite{YusubovMahmudov2024ACS}.

\paragraph*{Notes for computation.}
The system \eqref{eq:fpmp-state}–\eqref{eq:fpmp-min} defines a nonlocal two-point boundary-value problem. In practice, two routes are common. (i) \emph{Diffusive/Prony embeddings} (Sec.~\ref{subsec:eq-rep}) convert the FO dynamics into a classical ODE in an augmented state, so standard TPBVP solvers apply; accuracy is tuned by the number and placement of modes, with Caputo/RL-consistent initialization \cite{Monje2010,LiZeng2015}. (ii) \emph{Frequency-domain rationalization} (Oustaloup/CRONE) provides a low-order approximation for loop-shaping and time simulation, enabling PMP or direct transcription on the rational surrogate. Rigorous PMP formulations covering Caputo/RL dynamics, free time, and constraints are given in \cite{BergouniouxBourdin2020,Kamocki2014,NdairouTorres2023,Yusubov2023}; see also \cite{Cruz2025AIMSmath,Moon2025AIMSmath} for distributed-order kernels and state constraints.

\subsection{Dynamic Programming and FO--HJB}
\label{subsec:fo-hjb}

\paragraph*{From non-Markovian value to path-dependent HJB.}
Fractional dynamics introduce \emph{history dependence}, hence the value functional is naturally \emph{path dependent}. For a finite horizon $[0,T]$ with running cost $\ell$ and terminal cost $\Phi$, define
\[
V\big(t,\mathbf{x}_{[0,t]}\big)
= \inf_{u(\cdot)\in\mathcal{U}}
\Big\{ \Phi\!\big(x(T)\big) + \int_{t}^{T} \ell\big(\tau,x(\tau),u(\tau)\big)\,d\tau \Big\},
\]
where $\mathbf{x}_{[0,t]}$ is the past trajectory up to $t$. Under appropriate regularity/viability assumptions, a dynamic programming principle (DPP) holds for fractional systems, and $V$ solves a \emph{path-dependent} HJ equation in the viscosity sense on the history space; see \cite{Gomoyunov2020SICON} and, for classical DPP/HJB and path-dependent PDE frameworks, \cite{BardiCapuzzo1997,FlemingSoner2006,CrandallLions1992,Dupire2009,ContFournié2013,EkrenTouziZhang2016}. See also \cite{Gomoyunov2024JDE,Gomoyunov2025COCV} for recent extensions.

\paragraph*{Time-fractional HJB (backward in time).}
An equivalent (state-based) representation recasts the DPP as a \emph{time-fractional} HJB with a backward fractional derivative. For Caputo modeling one formally obtains
\begin{equation}
\label{eq:tfHJB-Caputo}
\begin{aligned}
({}^{\mathrm{C}}\!D_{T-}^{\alpha} V)(t,x) \;+\; H\!\big(t,x,\nabla_x V(t,x)\big) &= 0,\\
V(T,x) &= \Phi(x),
\end{aligned}
\end{equation}
while a Riemann--Liouville time operator yields
$({}^{\mathrm{RL}}\!D_{T-}^{\alpha} V)+H=0$ with the corresponding terminal data
\cite{Gomoyunov2020SICON,Gomoyunov2022COCV}. Here
$H(t,x,p)=\inf_{u\in U}\{\ell(t,x,u)+p^{\top}f(t,x,u)\}$.
Well-posedness is established in the viscosity framework (comparison, stability, uniqueness) under structural conditions on $H$ and the kernel \cite{CrandallLions1992,CamilliGoffi2020,Gomoyunov2022COCV}; see \cite{Yangari2024ADE} for variable-order Caputo regularity.
For zero-sum games, the fractional HJ--Isaacs equation replaces $H$ with the Isaacs Hamiltonian
$H^{\pm}(t,x,p)=\sup_{v\in V}\inf_{u\in U}(\cdots)$ (or vice versa); DPP and viscosity solutions for the fractional case are developed in \cite{Gomoyunov2024JDE}, extending the classical theory \cite{BasarOlsder1999,Isaacs1965}.

\paragraph*{Markovian embeddings to bypass path dependence.}
A practical route is to \emph{Markovize} fractional dynamics by augmenting the state with
$M$ \emph{diffusive modes} (Prony/relaxation spectrum) or by a band-limited Oustaloup realization
(see Sec.~\ref{subsec:eq-rep}). Let $z\in\mathbb{R}^M$ collect these modes; the FO plant then becomes
a classical ODE in $(x,z)$, so standard Hamilton--Jacobi--Bellman (HJB) and
Hamilton--Jacobi--Bellman--Isaacs (HJBI) methods apply on the augmented space.
The approximation error is set by $M$ (diffusive) or by the rational order $2N{+}1$ (Oustaloup/CRONE).
In effect, we trade the \emph{curse of history} (nonlocal memory) for the usual
\emph{curse of dimensionality} in dynamic programming \cite{Monje2010,Sabatier2015,FalconeFerretti2013}.
Fast SOE kernel approximations can further reduce memory and computation \cite{Guglielmi2025SOE}.

\paragraph*{Numerical schemes (deterministic).}
Two complementary routes are common:
\begin{itemize}[leftmargin=*]
\item \textit{Time-fractional HJB on $(t,x)$.}
Discretize ${}^{\mathrm{C}}D_{T-}^{\alpha}$ (or RL) via L1/GL/convolution–quadrature in time,
and use monotone spatial schemes (upwind/semi-Lagrangian) for $H$.
Convergence follows the viscosity framework under monotonicity, stability, and consistency
\cite{BarlesSouganidis1991,LiZeng2015,JinLazarovZhou2016,CamilliGoffi2020,FalconeFerretti2013}.
SOE compression reduces the $O(N_t)$ history to about $O(M\log N_t)$ updates
\cite{Jiang2017FastCaputo,LiZeng2015,Zhu2019SOE,ZhangFangSun2021ESA}.
\item \textit{Augmented-state HJB on $(x,z)$.}
Solve on the higher-dimensional grid with semi-Lagrangian value/policy iteration or filtered schemes
\cite{FalconeFerretti2013,Bokanowski2015}; see Sec.~\ref{subsec:eq-rep} for building $(x,z)$ models.
For coupled fractional HJB–FPK problems, see \cite{CamilliDuisembay2021JDG}.
\end{itemize}
For games (HJBI), policy iteration and primal–dual splitting are used routinely; the fractional aspect
changes only the time operator or the lifted state dimension \cite{Gomoyunov2024JDE,FalconeFerretti2013}.

\paragraph*{Numerical schemes (deterministic).}
Two complementary routes are common:
\begin{itemize}[leftmargin=*]
\item \textit{Time-fractional HJB on $(t,x)$.}
Discretize ${}^{\mathrm{C}}D_{T-}^{\alpha}$ (or RL) via L1/GL/convolution–quadrature in time,
and use monotone spatial schemes (upwind/semi-Lagrangian) for $H$.
Convergence follows the viscosity framework under monotonicity, stability, and consistency
\cite{BarlesSouganidis1991,LiZeng2015,JinLazarovZhou2016,CamilliGoffi2020,FalconeFerretti2013}.
SOE compression cuts the $O(N_t)$ history to near $O(M\log N_t)$ updates
\cite{Jiang2017FastCaputo,LiZeng2015,Zhu2019SOE,ZhangFangSun2021ESA}.
\item \textit{Augmented-state HJB on $(x,z)$.}
Work on the higher-dimensional grid with semi-Lagrangian value/policy iteration or filtered schemes
\cite{FalconeFerretti2013,Bokanowski2015}; see Sec.~\ref{subsec:eq-rep} for building $(x,z)$ models.
For coupled fractional HJB–FPK problems, see \cite{CamilliDuisembay2021JDG}.
\end{itemize}
For games (HJBI), policy iteration and primal–dual splitting are used routinely; the fractional aspect
changes only the time operator or the lifted state dimension \cite{Gomoyunov2024JDE,FalconeFerretti2013}.

\paragraph*{Stochastic interpretations and time change.}
In stochastic control, time-fractional evolution may arise from \emph{subordination} of Markov processes (random time change), leading to fractional Kolmogorov and HJB equations; see \cite{BaeumerMeerschaert2001,Kolokoltsov2010}. This viewpoint helps derive FO-HJB and to design consistent schemes via probabilistic representations.

\paragraph*{Takeaways.}
(i) DPP and viscosity well-posedness extend to FO systems, either as path-dependent PPDEs or as time-fractional HJB/HJBI with right-sided derivatives in $t$; (ii) for computation, diffusive/Oustaloup embeddings \emph{Markovize} memory and enable standard DP solvers, while direct time-fractional HJB benefits from L1/GL and SOE accelerations; (iii) for games, recent FO-HJBI theory ensures DPP under Isaacs-type conditions and supports robust/minimax synthesis \cite{Gomoyunov2020SICON,Gomoyunov2022COCV,Gomoyunov2024JDE}.%
%% NEW:
~See also \cite{Gomoyunov2025COCV,CamilliDuisembay2021JDG} for extensions to Volterra kernels and fractional HJB--FPK couplings.

\subsection{Fractional Differential Games (FDGs)}
\label{subsec:fdg}

We consider two–player dynamic games where the plant has fractional memory. Let $x:[0,T]\to\mathbb{R}^n$ evolve under Caputo dynamics
\begin{equation}
\label{eq:fdg-dyn}
\begin{aligned}
({}^{\mathrm{C}}\!D^{\alpha}_{0+} x)(t) &= f\big(t,x(t),u(t),v(t)\big),\\
x(0) &= x_0,\qquad 0<\alpha<1,
\end{aligned}
\end{equation}
with Player~1 control $u(\cdot)\in\mathcal{U}$ and Player~2 control $v(\cdot)\in\mathcal{V}$ (measurable, typically compact–valued).
For a zero–sum game we use the Bolza payoff
\begin{equation}
\label{eq:fdg-cost}
J(t_0,\phi;u(\cdot),v(\cdot))
= \Phi\!\big(x(T)\big) + \int_{t_0}^{T} \ell\big(t,x(t),u(t),v(t)\big)\,dt,
\end{equation}
where $\phi$ denotes the history consistent with Caputo initialization. (For RL modeling, fractional initial data enter and the adjoint side carries right Caputo operators; cf.\ Sec.~\ref{subsec:frac-defs}.)

Define the lower and upper values
\begin{equation}
\label{eq:fdg-values}
V^{-}(t_0,\phi)=\sup_{v(\cdot)}\inf_{u(\cdot)} J,\qquad
V^{+}(t_0,\phi)=\inf_{u(\cdot)}\sup_{v(\cdot)} J.
\end{equation}
Because of memory, the value is \emph{path dependent}. A rigorous treatment works on a history space with pathwise (coinvariant) derivatives and a dynamic programming principle (DPP); see \cite{Gomoyunov2020SICON,Gomoyunov2022COCV} and the classical game foundations in \cite{BasarOlsder1999,Isaacs1965}. Recent FO path–dependent HJBI results for Caputo/Volterra kernels appear in \cite{Gomoyunov2024JDE,Gomoyunov2025COCV}.

\paragraph*{FO–HJBI (time–fractional form).}
Equivalently, one can write a time–fractional Isaacs equation (backward in time) for the state–based value:
\begin{equation}
\label{eq:fo-hjbi}
\begin{aligned}
({}^{\mathrm{C}}\!D_{T-}^{\alpha} V^{\pm})(t,x)
&\;+\; H^{\pm}\!\big(t,x,\nabla_x V^{\pm}(t,x)\big)=0,\\
V^{\pm}(T,x)&=\Phi(x),
\end{aligned}
\end{equation}
with Isaacs Hamiltonians
\begin{equation}
\label{eq:isaacs-H}
\begin{aligned}
H^{-}(t,x,p)&=\sup_{v\in \mathcal{V}}\,\inf_{u\in \mathcal{U}}
\big\{\ell(t,x,u,v)+p^{\top}f(t,x,u,v)\big\},\\
H^{+}(t,x,p)&=\inf_{u\in \mathcal{U}}\,\sup_{v\in \mathcal{V}}
\big\{\ell(t,x,u,v)+p^{\top}f(t,x,u,v)\big\}.
\end{aligned}
\end{equation}
Under standard structural assumptions and the Isaacs condition $H^{-}\equiv H^{+}$,
one has $V^{-}=V^{+}=V$, and viscosity well–posedness holds in the fractional setting
\cite{Gomoyunov2020SICON,Gomoyunov2022COCV,Gomoyunov2024JDE}. For computation,
diffusive or Oustaloup embeddings (Sec.~\ref{subsec:eq-rep}) Markovize the memory and allow standard HJBI solvers on the augmented state.

\paragraph*{HJBI for time–fractional dynamics.}
Under standard boundedness/Lipschitz conditions and Isaacs’ condition, $V^{+}=V^{-}=V$ and the (history–based) value admits a state–based, \emph{time–fractional} HJBI representation in the viscosity/minimax sense:
\begin{align}
\label{eq:fdg-hjbi}
\big({}^{\mathrm{C}}\!D^{\alpha}_{T-} V\big)(t,x)
&\;+\;
\sup_{v\in\mathcal{V}}\,\inf_{u\in\mathcal{U}}
\Big\{ \ell(t,x,u,v) \notag\\
&\qquad + \nabla_x V(t,x)^{\!\top} f(t,x,u,v) \Big\}
= 0,
\end{align}

with terminal condition $V(T,x)=\Phi(x)$. Equivalently, one may use the Isaacs Hamiltonians $H^{\pm}$.
Existence/uniqueness for time–fractional HJ/HJBI and their equivalence to minimax solutions are established in
\cite{Gomoyunov2021Math,Gomoyunov2023JDE,Pozza_2025}, and regularity/approximation results for Caputo–HJ equations appear in \cite{Camilli2019CaputoHJ}; see also \cite{Gomoyunov2024JDE} for path–dependent kernels and comparison.

\paragraph*{PMP–based best responses and Stackelberg FDGs.}
Using the fractional PMP (Sec.~\ref{subsec:fpmp}), each player’s best response solves a two–point boundary–value problem with a right–sided fractional adjoint. In Stackelberg (leader–follower) games, the follower’s fractional PMP yields a reaction map; substituting it into the leader’s fractional OCP produces a bilevel problem. Analytical PMP–type developments for Caputo/RL games and pursuit–evasion settings include \cite{Chikrii2010,Matychyn2023}; recent variants handling state constraints and generalized kernels can be incorporated in best–response solvers \cite{Moon2025AIMSmath,Gasimov2024QTDS}.
\paragraph*{LQ FDGs via augmented realizations.}
For linear dynamics with quadratic costs, a diffusive (Prony) or Oustaloup realization
(Sec.~\ref{subsec:eq-rep}) lifts the FO plant to a finite–dimensional ODE game in an augmented
state. One can then apply the standard LQ game toolkit—feedback Nash/minimax strategies and the
associated Isaacs/Riccati equations \cite{BasarOlsder1999,Engwerda2005,Dockner2000}. Accuracy is
governed by the number of diffusive modes (or by the filter order $2N{+}1$), and closed–loop
guarantees track the usual LQ conditions over the approximation bandwidth.

\paragraph*{Mean–field and distributed variants.}
In large–population limits, mean–field games with \emph{spatially} fractional diffusion
(fractional Laplacian) are well developed \cite{Achdou2013,Graber2025MFG}, whereas
\emph{time}–fractional mean–field formulations remain largely open. For distributed–parameter FO
dynamics and hybrid pursuit–evasion settings, see \cite{Mamatov_2021,ChikriiMatychyn2012}.

\paragraph*{Numerical comments.}
Two practical routes are used in practice. 
(i) \emph{DP/HJBI route:} discretize the Caputo time operator (L1/L2 family), and solve the Isaacs operator by policy/value iteration; convergence relies on viscosity–solution stability/consistency for time–fractional HJ/HJBI \cite{Camilli2019CaputoHJ,Gomoyunov2021Math}. Fast sum–of–exponentials (SOE) stepping further compresses the history cost \cite{Zhu2019SOE}. 
(ii) \emph{Augmented–state route:} build a low–order ODE game via diffusive/CRONE approximations and apply standard solvers (policy iteration; Riccati/Isaacs in LQ cases). The choice trades memory fidelity against computational complexity; windowing and band–limited realizations help tame the curse of history.

\paragraph*{Open questions.}
Refinements are still needed on: (a) Isaacs–type conditions under memory; (b) existence/uniqueness for general–sum FO games; (c) verifiable sufficient conditions for Stackelberg equilibria in FDGs; and (d) scalable algorithms with a priori/a posteriori certificates that couple fractional–approximation error with game–theoretic suboptimality \cite{Gomoyunov2022COCV,Matychyn2023}. Extensions to Volterra kernels and path–dependent Isaacs conditions are outlined in \cite{Gomoyunov2025COCV}.

\section{Canonical Problems: LQR, Tracking, and Regulation}
\label{sec:canonical}

We study linear--quadratic regulator (LQR), tracking, and regulation when the plant has fractional-order memory. In general, there is \emph{no} Riccati equation directly in the original FO coordinates; in practice, two routes are used.%
%% NEW: recent FO-LQ value-functional approach (no Riccati in original FO coordinates)
~See also the 2024 value-functional formulation with $\varepsilon$-optimal feedback in~\cite{Gomoyunov2024MCRF}.

\begin{itemize}[leftmargin=*]
\item \textbf{Route A (augmented ODE / Riccati).}
Approximate the fractional operator by a finite-dimensional realization: (i) diffusive/Prony, or (ii) CRONE/Oustaloup (Sec.~\ref{subsec:eq-rep}). The FO plant becomes an ODE in an \emph{augmented state}, so classical Riccati machinery applies—differential and algebraic Riccati equations (DRE/ARE) for LQR/tracking/regulation~\cite{AndersonMoore2007,Monje2010,Sabatier2015}.%
%% NEW: recent FO-LQR application
~For a FO-LQR tracking/regulation demo on magnetic levitation, see~\cite{Yoneda2024FractalFract}. Accuracy is set by the number of diffusive modes or the filter order $2N{+}1$, and guarantees track the usual LQ conditions over the approximation bandwidth.

\item \textbf{Route B (direct transcription as QP).}
Discretize the Caputo/RL dynamics in time (e.g., L1/GL/convolution–quadrature). The resulting linear, lower-triangular \emph{history-coupled} constraints and the quadratic performance index yield a convex quadratic program (QP) with block Toeplitz/Hessenberg structure; the finite-horizon problem is then solved by QP/least squares (no closed-form Riccati in the original FO coordinates)~\cite{LiZeng2015,Alikhanov2015,Lubich1986}.%
%% NEW: FO-LQ as QP (2024) + L1-based FO optimal control (2024)
~Recent works cast pure/multi-order FO-LQ explicitly as QPs~\cite{Malmir2024ISA} and analyze L1-discretized FO optimal control with domain decomposition~\cite{Leugering2024FractFract}.
\end{itemize}
\subsection{FO–LQR via augmented realizations}
\label{subsec:fo-lqr}

Consider the commensurate Caputo LTI model ($0<\alpha<1$)
\begin{equation*}
{}^{\mathrm{C}}D_{0+}^{\alpha}x(t)=A x(t)+B u(t),\qquad x(0)=x_0,
\end{equation*}
with quadratic performance cost
\begin{equation}
\label{eq:lqr-cost}
\begin{aligned}
J &= x(T)^{\!\top} S\, x(T) 
   \;+\; \int_{0}^{T}\!\!\big(x^{\!\top}Qx + u^{\!\top} R u\big)\,dt,\\
Q &\succeq 0,\qquad R \succ 0 .
\end{aligned}
\end{equation}

Using a diffusive/Prony or Oustaloup realization (Sec.~\ref{subsec:eq-rep}), the FO plant is lifted to
\begin{equation}
\label{eq:aug-ode}
\dot{X}(t)=A_{\mathrm{aug}}\,X(t)+B_{\mathrm{aug}}\,u(t),\qquad
X=\begin{bmatrix}x\\ z_1\\ \vdots\\ z_M\end{bmatrix}.
\end{equation}
To preserve the original objective \eqref{eq:lqr-cost}, choose, e.g.,
\[
Q_{\mathrm{aug}}=\mathrm{blkdiag}(Q,\,0),\qquad
S_{\mathrm{aug}}=\mathrm{blkdiag}(S,\,0),
\]
or add small weights on $z_i$ if desired for numerical conditioning.

If $(A_{\mathrm{aug}},B_{\mathrm{aug}})$ is stabilizable and $(Q_{\mathrm{aug}}^{1/2},A_{\mathrm{aug}})$ is detectable, the algebraic Riccati equation (ARE)
\begin{equation}
\label{eq:are-aug}
A_{\mathrm{aug}}^{\!\top}P + P A_{\mathrm{aug}}
- P B_{\mathrm{aug}} R^{-1} B_{\mathrm{aug}}^{\!\top} P
+ Q_{\mathrm{aug}} = 0
\end{equation}
admits a unique stabilizing solution $P\succeq0$, and the optimal feedback is
\begin{equation}
\label{eq:lqr-gain}
u^{\ast}(t)= -K\,X(t),\qquad
K= R^{-1} B_{\mathrm{aug}}^{\!\top} P ,
\end{equation}
yielding an ODE closed loop with eigenvalues in the open left half–plane \cite{AndersonMoore2007}. 
The physical state $x(t)$ then exhibits the expected Mittag–Leffler–type decay, consistent with the approximation bandwidth and the Matignon sector condition for $A$ \cite{Matignon1996,Mainardi2010,Monje2010}. 
For finite horizons, the differential Riccati equation (DRE) is used instead of \eqref{eq:are-aug}.

\paragraph*{Accuracy–complexity trade–off.}
Increasing the number of diffusive modes or the Oustaloup order enlarges the valid frequency band and reduces phase ripple, at the cost of a higher augmented dimension (and ARE/DRE size) \cite{Sabatier2015,LiZeng2015}. Band selection should cover the anticipated closed–loop bandwidth and disturbance spectrum. Fast SOE/diffusive kernel compression can shrink the augmentation while maintaining accuracy \cite{GuglielmiHairer2025}.
\subsection{Direct transcription LQR (history–coupled QP)}
\label{subsec:fo-lqr-direct}

A left-sided time discretization on the grid $t_k = k\Delta t$ (e.g., L1 or Grünwald--Letnikov)
yields the linear, history–coupled update
\begin{equation}
\label{eq:disc-fo}
\sum_{j=0}^{k} a_{k-j}\, x_j \;=\; A\,x_k + B\,u_k, 
\qquad k=1,\dots,N,
\end{equation}
where $\{a_j\}$ are the scheme-dependent fractional weights (scaled binomial weights for GL,
and L1/L2-type weights for Caputo discretizations) \cite{Lubich1986,Alikhanov2015,LiZeng2015}.

Stacking \eqref{eq:disc-fo} for $k=1,\dots,N$ gives a lower–triangular linear relation between
the stacked state and input trajectories,
\[
x \;=\; \mathcal{G}\,u \;+\; {h}(x_0),
\]
where $x := [x_1^\top,\dots,x_N^\top]^\top$, 
$u := [u_0^\top,\dots,u_{N-1}^\top]^\top$, 
$\mathcal{G}$ is the block “convolution’’ operator built from $\{a_j\}$, $A$, and $B$, and
${h}(x_0)$ collects the zero-input response due to the initial condition $x_0$.
Substituting this representation into the quadratic cost \eqref{eq:lqr-cost} yields a convex QP
in the control sequence $u$ with a block Toeplitz/Hessenberg normal matrix.

This direct-transcription route stays entirely in the time domain, handles pointwise constraints
naturally (see Sec.~\ref{subsec:constraints}), and avoids any bandwidth specification.

\paragraph*{QP form and structure.}
Writing $X=[x_1^\top\!\dots x_N^\top]^\top$ and $U=[u_0^\top\!\dots u_{N-1}^\top]^\top$, the cost takes the standard form
\[
\min_{U}\;\tfrac12\, U^\top H U + f^\top U \;+\; \text{const},
\]
with $H\succeq 0$ assembled from $(Q,R,S)$ and the stacked operators $(\mathcal{G},\mathcal{h})$. 
The matrices inherit lower–triangular/Toeplitz structure from \eqref{eq:disc-fo}, which enables sparse Cholesky or conjugate–gradient solvers and warm–starts in receding–horizon implementations.

\paragraph*{Accuracy and setup.}
The choice of scheme (GL vs.\ L1 and its higher–order variants) fixes the weights $\{a_j\}$ and the convergence order; see \cite{LiZeng2015,Alikhanov2015}. 
Smaller $\Delta t$ improves fidelity but increases horizon size and condition numbers—preconditioning and mild weights on auxiliary states (if any) help numerics. 
Initialization follows the Caputo setting directly through $x_0$; for RL models, fractional initial data must be incorporated in $\mathcal{h}(x_0)$ (cf.\ Sec.~\ref{subsec:frac-defs}). 
When constraints are present, they enter as linear inequalities in $U$ without altering the convexity of the QP.

\paragraph*{When to use.}
Direct transcription is preferable when (i) constraints dominate the design, (ii) frequency–band guarantees are not required, or (iii) one wants a pure time–domain workflow with transparent error control via $\Delta t$ and the chosen fractional scheme \cite{Malmir2024ISA,Leugering2024FractFract}.
\subsection{Tracking and regulation}
\label{subsec:tracking}

For state/setpoint tracking with reference $x_r(t)$ (constant or time–varying), define the error
$e(t)=x(t)-x_r(t)$, where $x_r$ is admissible for the plant or provided exogenously. Two standard designs are:

\begin{itemize}[leftmargin=*]
\item \textbf{Augmented–error LQR.}
Form the augmented state $[\,e^{\top}\; z_1^{\top}\!\cdots z_M^{\top}]^{\top}$ and penalize $e$ and $u$ in
\eqref{eq:lqr-cost}. A steady–state feedforward $u_{\mathrm{ff}}$ can be added to cancel offsets
(\,solve $A x_\infty + B u_\infty = 0$ in the augmented realization\,).%
~See \cite{Yoneda2024FractalFract} for an experimental FO–LQR tracking study.

\item \textbf{Internal–model principle (IMP).}
To track steps/sinusoids and reject structured disturbances, embed the appropriate internal model
in the augmented ODE and apply LQR to the resulting augmented state. This mirrors the classical IMP
analysis \cite{FrancisWonham1976,Francis1977} within the chosen FO approximation bandwidth.
\end{itemize}

For regulation ($x_r\equiv 0$), one recovers the state feedback in \eqref{eq:lqr-gain}.
For nonzero $x_r$, a prefilter or reference governor is used so that $e(t)\to 0$ under
stabilizability/detectability of the augmented pair and reachability of the desired steady state
\cite{AndersonMoore2007}.

\subsection{Handling constraints}
\label{subsec:constraints}

Because FO memory couples distant time steps, constraints are most conveniently enforced on a
finite horizon:

\begin{itemize}[leftmargin=*]
\item \textbf{Route A (history–coupled QP).}
With the discretization \eqref{eq:disc-fo}, add box/polyhedral bounds on $\{x_k,u_k\}$. The problem
is a convex QP with lower–triangular/Toeplitz structure, suitable for active–set or interior–point
methods and warm starts across horizons.%
~Fast Caputo updates via SOE kernels can lower per–iteration cost \cite{GuglielmiHairer2025}.

\item \textbf{Route B (MPC on the augmented ODE).}
Use the augmented realization \eqref{eq:aug-ode} with a quadratic stage and terminal cost (and, if
needed, a terminal set). Standard MPC stability tools (positive invariance and terminal ingredients)
apply \cite{RawlingsMayneDiehl2017,Mayne2000}. Closed–loop guarantees hold within the approximation
band; windowing (finite memory) is often used in prediction to cap complexity.
\end{itemize}

\paragraph*{Implementation notes.}
(i) Match the approximation band to the expected closed–loop bandwidth and dominant disturbances;
(ii) for Caputo models, zero diffusive states at $t{=}0$ (classical initialization); for RL models,
map the fractional initial data into the augmented states (see Sec.~\ref{subsec:frac-defs});
(iii) report band limits, approximation order, time step, and window length for reproducibility
(\,see the reporting checklist in Table~\ref{tab:summary}\,).

\section{Fractional-Order Controllers (FOPID) via Optimal Control}
\label{sec:fopid-via-oc}

FOPID controllers extend classical PID by allowing non-integer integration and differentiation orders,
\begin{equation}
C_{\mathrm{FOPID}}(s) \;=\; k_p \;+\; k_i\, s^{-\lambda} \;+\; k_d\, s^{\mu},
\qquad \lambda,\mu \in [0,1],
\label{eq:fopid}
\end{equation}
so that $(\lambda,\mu)$ shape low/high-frequency asymptotics and phase. The added phase flexibility
enables near \emph{iso-damping} (closed-loop damping insensitive to gain changes), a central theme in
CRONE loop-shaping \cite{Sabatier2015,Oustaloup2000,Monje2010}. Recent hardware studies (e.g., twin-rotor/TRMS
and multi-DOF platforms) report robust, iso-damped behavior with FOPID implementations
\cite{Wendimu2025SciRep,Krzysztofik2025ApplSci}. This section treats FOPID tuning as an
optimal-control/optimization task and contrasts time- and frequency-domain criteria, robustness, and implementation.

\subsection{Modeling and approximation back-ends}
Practical tuning/certification of \eqref{eq:fopid} requires a finite-dimensional realization of $s^{\pm\alpha}$.
Two standard back-ends (cf.\ Sec.~\ref{subsec:eq-rep}) are:

\begin{itemize}[leftmargin=*]
\item \textbf{Frequency-domain (Oustaloup/CRONE).}
Approximate $s^{\pm\alpha}$ on a target band $[\omega_\ell,\omega_h]$ with a $(2N{+}1)$-term real-pole/zero product
\cite{Oustaloup2000,Sabatier2015}. This is natural for loop-shaping and $H_\infty$-style templates
(gain/phase margins, iso-damping). Comparative recent studies report competitive or superior robustness
versus $H_\infty$ designs on relevant plants \cite{deAzeredo2024IFAC,Mseddi2024FractalFract}.

\item \textbf{Diffusive/Prony (state augmentation).}
Realize the fractional kernels as a positive sum of stable first-order modes (Prony series).
The plant and controller then form a classical ODE in an \emph{augmented} state, enabling time-domain
optimal control, LQR-like surrogates, and MPC with constraints
\cite{Monje2010,LiZeng2015,Tepljakov2017}.
\end{itemize}

In both cases, the bandwidth and approximation order (either $N$ or the number of diffusive modes)
trade accuracy for dimension; the chosen band should cover the intended closed-loop bandwidth and dominant disturbances.

\subsection{Time-domain tuning as optimal control}
Let the plant $P$ (possibly fractional) be implemented as an augmented ODE
$\dot{X}=A_{\rm aug} X + B_{\rm aug} u$, $y=C_{\rm aug} X$
via a diffusive/Prony or Oustaloup realization. The controller $C_{\mathrm{FOPID}}$ is used in its
finite-dimensional form, and the tracking error is $e(t)=y(t)-r(t)$ for a given reference $r$.
A typical finite-horizon design reads
\begin{equation}
\label{eq:fopid-time-oc}
\min_{\theta}\; J_{\rm time}(\theta)
=\int_{0}^{T}\!\!\big(w_e\,e(t)^2+w_u\,u(t)^2\big)\,dt
\;+\; w_T\,e(T)^2,
\end{equation}
with decision variables $\theta=(k_p,k_i,k_d,\lambda,\mu)$ (or a chosen subset) and closed-loop
dynamics enforced by the augmented ODE. Alternative criteria such as IAE/ISE/ITAE and
multi-scenario aggregation (parametric uncertainty, actuator limits) are standard
\cite{ChenPetrasXue2009,Das2011,Tepljakov2017}.

\paragraph*{Computation and constraints.}
Problem \eqref{eq:fopid-time-oc} is generally nonconvex in $\theta$.
In practice: (i) fix $(\lambda,\mu)$ to a small grid and optimize $(k_p,k_i,k_d)$;
(ii) use smooth band-limited realizations to enable gradient-based methods (adjoint or sensitivities);
(iii) combine a coarse global search with a local solver for refinement.
Time-domain constraints (overshoot, settling time, $|u|\le u_{\max}$) can be imposed directly
in the augmented model or enforced as soft penalties. Recent pipelines for robust fractional
implementations provide end-to-end discretization and certification that fit this workflow
\cite{Mihaly2024FMI}.

\subsection{Frequency-domain criteria and loop shaping}
CRONE methodology prescribes a target magnitude slope and an almost flat phase around crossover to
achieve iso-damping and robust margins \cite{Sabatier2015,Oustaloup2000}. A convenient objective is
\begin{equation*}
\begin{aligned}
\min_{\theta}\; J_{\rm freq}(\theta)
=\;&\big\| W_S(j\omega)\,S(j\omega)\big\|_{\infty}
+\big\| W_T(j\omega)\,T(j\omega)\big\|_{\infty} \\
&\quad + \rho\,\Delta_{\phi},
\end{aligned}
\end{equation*}
where $L=C_{\mathrm{FOPID}}P$, $S=(I+L)^{-1}$, $T=L(I+L)^{-1}$, $W_S,W_T$ are design weights,
and $\Delta_\phi$ penalizes phase ripple on $[\omega_\ell,\omega_h]$.
With an Oustaloup realization, one can use off-the-shelf $H_\infty$/loop-shaping solvers directly.
A pragmatic template fixes $(\lambda,\mu)$ from the desired phase slope/iso-damping
and then tunes $(k_p,k_i,k_d)$ against the weighted $S/T$ goals
\cite{Sabatier2015,Monje2010}. Recent comparisons report competitive or superior robustness of
CRONE/FOPID relative to $H_\infty$ controllers on representative plants
\cite{deAzeredo2024IFAC,Mseddi2024FractalFract}.

\paragraph*{Reporting and reproducibility.}
For either route, report the approximation band $[\omega_\ell,\omega_h]$, the realization order
($2N{+}1$ or number of diffusive modes), the sampling/discretization used, and any constraints/weights.
These choices determine both performance and the range over which the certification is valid.
\subsection{Robustness}
FOPID’s extra orders $(\lambda,\mu)$ enable \emph{iso-damped} responses whose closed-loop damping
is largely insensitive to plant gain variations \cite{Sabatier2015}. We quantify robustness via:
\begin{itemize}[leftmargin=*]
\item \emph{Classical margins/templates:} gain/phase margins and Nichols templates with a near-flat
      phase at crossover.
\item \emph{Sensitivity norms:} $\|W_S S\|_{\infty}$ for disturbance rejection and
      $\|W_T T\|_{\infty}$ for noise amplification.
\item \emph{Parametric sweeps:} gridding plant parameters, band limits, and approximation order
      to certify performance envelopes on the augmented model.
\end{itemize}
For fractional plants that satisfy the Matignon sector condition, well-tuned FOPID loops exhibit
algebraic (Mittag–Leffler) decay, with rates shaped jointly by the plant orders and $(\lambda,\mu)$
\cite{Matignon1996,Mainardi2010,Monje2010}. Hardware studies (TRMS, 3-DOF gyroscope) report
iso-damped, robust behavior consistent with these metrics \cite{Wendimu2025SciRep,Krzysztofik2025ApplSci}.

\subsection{Digital implementation and practicalities}
Two implementation back-ends are standard \cite{Sabatier2015,Tepljakov2017}:
\begin{itemize}[leftmargin=*]
\item \textbf{CRONE/Oustaloup:} realize $s^{\pm\alpha}$ as a stable IIR on a target band
      $[\omega_\ell,\omega_h]$.
\item \textbf{Diffusive/Prony:} discretize the relaxation modes (e.g., backward-Euler/bilinear with
      prewarping).
\end{itemize}

\textbf{What to report (for reproducibility).}
\begin{itemize}[leftmargin=*]
\item band $[\omega_\ell,\omega_h]$ and realization order ($2N{+}1$ or \# diffusive modes);
\item sampling rate (avoid aliasing near $\omega_h$) and discretization choice;
\item initialization consistent with Caputo vs.\ RL modeling;
\item anti-windup policy if integral action is used ($k_i \neq 0$).
\end{itemize}

\textbf{Anti-windup.} Standard AW schemes carry over to the augmented ODE. When clamping, account for
fractional memory to prevent bias after saturation. 

\textbf{Benchmarking.} Share code and parameter
files (controller realization, band, order) to ensure results are repeatable across toolchains; see
also the end-to-end pipelines for robust fractional implementations \cite{Mihaly2024FMI}.

\paragraph*{Historical notes and tutorials.}
Early fractional compensators predate modern CRONE \cite{AxtellBise1990}. Practical recipes and
pitfalls for tuning and implementation are covered in \cite{ChenPetrasXue2009,Monje2010,Sabatier2015,Das2011,Tepljakov2017}.

\section{Model Predictive Control with FO Models}
\label{sec:mpc-fo}

Model predictive control for fractional–order plants poses two linked tasks:
(i) accurate \emph{prediction} in the presence of nonlocal memory, and
(ii) \emph{recursive feasibility} and stability when the state is history–dependent.
We summarize practical prediction models, their accuracy–complexity trade–offs, and the stability tools most used with FO dynamics.

\paragraph*{Prediction models for Caputo/RL dynamics.}
Consider a Caputo linear time–invariant (LTI) plant ${}^{\mathrm{C}}D_{0+}^{\alpha}x(t)=A x(t)+B u(t)$ with $0<\alpha<1$.
Three mainstream routes are used in FO–MPC (Riemann–Liouville (RL) cases are analogous).

\begin{enumerate}[leftmargin=*]
  \item \emph{Time–domain convolution/quadrature (direct).}
  Discretize ${}^{\mathrm{C}}D^{\alpha}$ with L1 or Grünwald–Letnikov (GL) schemes to obtain a lower–triangular Toeplitz history map. A representative step is
  \[
  x_{k+1} \;=\; \sum_{j=0}^{k} w_{\alpha}(j)\,x_{k-j} \;+\; B_{\!d}\,u_k,
  \]
  with scheme–dependent weights $\{w_{\alpha}(j)\}$; stacking over the horizon gives $X=\mathcal{S}x_k+\mathcal{T}U$ with Toeplitz $\mathcal{S}$ \cite{Lubich1986,Alikhanov2015,JinLazarovZhou2016,LiZeng2015}.
  This route is “exact’’ up to quadrature but stores and updates the full history.

  \item \emph{Diffusive (Prony/relaxation) augmentation.}
  Approximate the power–law kernel by a positive mixture of exponentials and add $M$ first–order modes:
  $\dot z_i=-\lambda_i z_i+\dot x$, with $({}^{\mathrm{C}}D^{\alpha}x)\approx \sum_{i=1}^{M} w_i z_i$ (Sec.~\ref{subsec:eq-rep}).
  After zero–order hold (ZOH) discretization, prediction is a standard LTI/LTV model in the augmented state $(x,z)$.
  Accuracy scales with $M$ and node placement \cite{BagleyTorvik1983,Koeller1984,LiZeng2015}.

  \item \emph{Frequency–domain (CRONE/Oustaloup) rationalization.}
  Replace $s^{\alpha}$ on a target band $[\omega_\ell,\omega_h]$ by a low–order stable rational model and use a classical discrete–time realization \cite{Oustaloup2000,Sabatier2015}.
  This integrates cleanly with loop–shaping; prediction is valid within $[\omega_\ell,\omega_h]$.
\end{enumerate}

\emph{Benchmarks reported in recent literature show energy/robustness gains from FO predictors inside MPC when compared with integer surrogates on several plants} \cite{Cao2024IFAC,Das2025IFAC}.

\paragraph*{Memory truncation, windowing, and fast summation.}
In direct time–domain MPC, the convolution carries an $\mathcal{O}(k)$
history cost. Two practical mitigations are:
(i) \emph{windowing} the memory to a trailing horizon $T_{\mathrm{w}}$
chosen from closed–loop bandwidth/settling–time considerations; and
(ii) \emph{sum–of–exponentials} (SOE) compression of the kernel, which
permits $O(1)$–state recursions per step via a diffusive augmentation
\cite{LiZeng2015,YanSunZhang2017CICP}. Recent “fast, memoryless’’ SOE
constructions combined with stiff integrators further scale down the
online burden while preserving accuracy for subdiffusion
\cite{GuglielmiHairer2025}.

\paragraph*{Offset–free tracking and disturbance accommodation.}
Offset–free regulation proceeds as in the integer–order case but in the
\emph{augmented} space used for prediction. A standard choice is to
embed an integrator or a random–walk disturbance on the output/plant
model that the optimizer and estimator share. Stabilizing FO–GPC designs
with explicit bias rejection follow this template
\cite{Ntouskas2018OffsetFreeFO,Ntouskas2019IFAC}. Experimental FO–adaptive MPC
implementations report similar benefits in power–electronics settings
\cite{Peng2025SciRep}.

\paragraph*{Robust FO–MPC and tube concepts.}
With bounded disturbances and modeling error in the FO operator/plant,
robust MPC can be built on an \emph{augmented} linear model (diffusive
or CRONE) and a tube strategy: optimize a nominal trajectory while a
pre–designed feedback renders a robust positively invariant tube around
it. Sufficient conditions for recursive feasibility and asymptotic
stability of discrete–time FO systems under MPC appear in
\cite{Sopasakis2016MPFOS}; related robust formulations and applications
are given in \cite{Ntouskas2019IFAC,Joshi2014FO}. Distributed and
large–scale variants with fractional couplings are addressed in
\cite{ChiIonescu2022}. Recent case studies on nonlinear plants and
adaptive identification also document robustness gains from fractional
predictors inside MPC loops \cite{Das2025IFAC,Peng2025SciRep}.
\paragraph*{Cost, constraints, and solver structure.}
FO–MPC typically uses a quadratic tracking cost on $(x,u)$ or $(y,\Delta u)$
over a horizon $N$, with $\ell_2$–regularization on input increments.
With CRONE or diffusive realizations the prediction is \emph{linear} in
the decision variables, so the optimizer is a convex QP (box/linear
constraints on $u$; polytopic constraints on the augmented state).
In direct–convolution FO–MPC the prediction matrices are lower–triangular
Toeplitz; both time and memory scale with $N$ unless windowed or SOE–compressed.

\paragraph*{Stability and terminal ingredients.}
Two standard routes provide closed–loop guarantees:
\begin{enumerate}[leftmargin=*]
  \item \emph{Augmented Lyapunov.} With diffusive/CRONE realizations, enforce
  a discrete–time Lyapunov inequality in the augmented model together with a
  quadratic terminal cost $V(x,z)=\|(x,z)\|_P^2$ and an invariant terminal set.
  This yields recursive feasibility and asymptotic stability (and recoverability
  for the original FO plant up to approximation error)
  \cite{Sopasakis2016MPFOS}.
  \item \emph{Dissipativity/sector arguments.} For commensurate FO LTI plants,
  the Matignon sector condition (Sec.~\ref{subsec:stability}), combined with a
  suitable terminal penalty and a sufficiently long horizon, ensures a decreasing
  value function; FO–GPC analyses along these lines are reported in
  \cite{Ntouskas2018OffsetFreeFO,Romero2013MPE}.
\end{enumerate}

\paragraph*{Applications and case studies.}
FO–MPC/GPC has been demonstrated in process and thermal systems and in
energy devices modeled by constant–phase elements (CPEs), where diffusive
or CRONE realizations yield compact predictors
\cite{Romero2013MPE,ChiIonescu2022}. Early studies on fuel cell power
regulation and servo drives illustrate practical integration of FO models
with predictive schemes \cite{Joshi2014FO}. Recent reports cover stabilization
of chaotic fractional systems and adaptive FO–MPC on hardware converters
\cite{Das2025IFAC,Peng2025SciRep}.

\section{Numerical Methods}
\label{sec:numerics}

This section collects practical tools for simulating and discretizing
fractional–order dynamics in analysis and control. We organize the
discussion into three families: (i) time–domain schemes for Caputo and
Riemann–Liouville operators, (ii) frequency–domain rational
approximations (Oustaloup/CRONE), and (iii) diffusive (Prony) state
augmentations. For each family we summarize accuracy and stability
properties, outline computational cost, and note implementation details
relevant to control design.

The Adomian Decomposition Method (ADM) provides a semi-analytical framework for solving nonlinear dynamical systems. An exposition of the method for conventional ordinary differential equations is presented in \cite{Fatoorehchi2023sign}, which may serve as a useful reference for the fundamental concepts and implementation of ADM. Extensions of ADM to fractional-order systems have been reported in several studies, demonstrating its applicability to models with memory and nonlocal dynamics \cite{Obeidat21062024,Ziada2024,Sebaq2025}.
\subsection{Time–Domain Discretizations}

\paragraph*{Convolution quadrature (CQ).}
Given a sectorial operator $A$ and kernel $k(t)$ with Laplace transform $K(s)$,
CQ uses an $r$–step linear multistep method with generating function $\delta(\zeta)$
to define discrete weights
\begin{equation*}
\omega_j^{(\alpha)} \;=\; \frac{1}{2\pi i}
\int_{\Gamma} K\!\left(\frac{\delta(\zeta)}{\Delta t}\right)\zeta^{-j-1}\, d\zeta,
\qquad j\ge 0,
\end{equation*}
so that, on the grid $t_n=n\Delta t$,
\begin{equation*}
\sum_{j=0}^{n} \omega_j^{(\alpha)} x_{n-j} \;\approx\; \big({}^{\mathrm{C}}\!D^{\alpha}x\big)(t_n).
\end{equation*}
CQ inherits the $A$–stability of the base method and applies directly to FO
diffusion/wave equations and linear control models
\cite{Lubich1986,CuestaLubichPalencia2006}. Fast/oblivious variants reduce
cost from $\mathcal{O}(N^2)$ to $\mathcal{O}(N\log N)$ time and
$\mathcal{O}(\log N)$ memory via FFT/Talbot techniques
\cite{Schadle2006FOCQ}; memory–light SOE–CQ further compresses history in stiff IVPs
\cite{GuglielmiHairer2025}.

\paragraph*{L1 and L2–$1\sigma$ formulas.}
On a uniform grid and $0<\alpha<1$, the L1 approximation of the Caputo derivative is
\begin{align*}
\big({}^{\mathrm{C}}\!D^{\alpha}x\big)(t_n)
&\;\approx\;
\frac{\Delta t^{-\alpha}}{\Gamma(2-\alpha)}
\sum_{j=0}^{n-1}
\Big[(n\!-\!j)^{1-\alpha}-(n\!-\!j\!-\!1)^{1-\alpha}\Big] \\
&\qquad\times (x_{j+1}-x_j).
\end{align*}

which is first–order accurate with suitable initial corrections
\cite{JinLazarovZhou2016}. The L2–$1\sigma$ scheme attains second order under
mild regularity while preserving stability for subdiffusion \cite{Alikhanov2015}.
Refinements (averaged–L1, corrected–L1) tighten practical error bounds
\cite{ZhengWang2024JSC,HuangLuoZhang2025}.

\paragraph*{Predictor–corrector (Adams–Bashforth–Moulton).}
FO IVPs can be written as Volterra equations and advanced by a PECE
predictor–corrector tailored to the kernel $(t-\tau)^{\alpha-1}$
\cite{DiethelmFord2002}. The method is robust for smooth nonlinear plants and supports
stepsize adaptivity; parallel–in–time accelerations have been reported
\cite{HeLiWei2025}.

\paragraph*{Complexity and fast kernels.}
Direct history summation costs $\mathcal{O}(N^2)$ time and $\mathcal{O}(N)$ memory.
Fast/oblivious CQ achieves $\mathcal{O}(N\log N)$ and $\mathcal{O}(\log N)$
\cite{Schadle2006FOCQ}. Sum–of–exponentials (SOE) approximations,
\[
\frac{1}{\Gamma(1-\alpha)}(t-\tau)^{-\alpha}
\;\approx\; \sum_{i=1}^{M} w_i\, e^{-\lambda_i (t-\tau)},
\]
enable $\mathcal{O}(M)$–state recursions with controllable error
\cite{YanSunZhang2017CICP,JiangZhang2017SISC}; Toeplitz/circulant preconditioners further
reduce linear–solve cost in high–order time–fractional discretizations
\cite{GanZhangLiang2024JAMC}.

\paragraph*{Space–fractional PDEs.}
For GL–type space–fractional operators, finite volume/finite–element
discretizations yield stable, convergent schemes compatible with control loops
\cite{Mojahed2020PFDA,ZhangYang2024FractalFract}.

\subsection{Frequency–Domain Rational Approximations (Oustaloup/CRONE)}
For band–limited designs it is convenient to replace $s^{\alpha}$ ($\alpha\in\mathbb{R}$) by a stable real–rational transfer function on $[\omega_\ell,\omega_h]$:
\[
s^{\alpha}\;\approx\; K\;\prod_{k=-N}^{N}\frac{1+s/\omega_{z,k}}{1+s/\omega_{p,k}},
\]
where $\{(\omega_{z,k},\omega_{p,k})\}$ are geometrically distributed (Oustaloup recursion) so that the magnitude slope is close to $\alpha$ (dB/dec) and the phase is nearly constant $\alpha\cdot 90^\circ$ over the target band. The gain $K$ matches a reference frequency (typically $\sqrt{\omega_\ell\omega_h}$), and increasing $N$ widens the accuracy band and reduces phase ripple \cite{Oustaloup2000,Sabatier2015}. Discrete–time realizations (e.g., bilinear transform with prewarping) yield standard LTI predictors for MPC/LQR. Band limits should cover the intended closed–loop bandwidth; if passivity/positive–realness is required, CRONE design rules provide appropriate pole–zero templates.

\subsection{Diffusive / Prony State Augmentation}
The Caputo derivative admits an exponential–mixture (diffusive) representation:
\begin{align*}
({}^{\mathrm{C}}\!D^{\alpha}x)(t)
&=\int_{0}^{\infty}\mu_{\alpha}(\lambda)\,z(\lambda,t)\,d\lambda, \\
\partial_t z(\lambda,t)
&=-\lambda z(\lambda,t)+\dot{x}(t).
\end{align*}

with $\mu_{\alpha}(\lambda)=\lambda^{\alpha-1}/\big(\Gamma(\alpha)\Gamma(1-\alpha)\big)$ and $z(\lambda,0)=0$ (Caputo initialization). Approximating the integral by a positive quadrature $\sum_{i=1}^M w_i\,z(\lambda_i,t)$ produces a finite cascade of stable first–order modes (a Prony series) that preserves causality—and, with $w_i>0$, passivity—while delivering a classical augmented model in $(x,z)$ for analysis and synthesis \cite{BagleyTorvik1983,Koeller1984,LiZeng2015}. After ZOH discretization one obtains an LTI/LTV predictor whose accuracy is set by $M$ and the (typically log–spaced) nodes $\{\lambda_i\}\subset[\lambda_{\min},\lambda_{\max}]$ tied to the design band $[\omega_\ell,\omega_h]$. Recent “memoryless’’ SOE implementations combined with stiff integrators further reduce long–horizon cost without sacrificing accuracy \cite{GuglielmiHairer2025}.

\subsection{Error, Stability, and Selection Guidelines}

\paragraph*{Accuracy.}
CQ based on BDF$k$ attains $k$th–order accuracy in time for sufficiently smooth data (subject to the stability of the underlying multistep method) \cite{Lubich1986,CuestaLubichPalencia2006}. The L1 scheme is first order and can recover higher effective rates with start–up corrections; L2–$1\sigma$ is second order under standard regularity assumptions \cite{Alikhanov2015,JinLazarovZhou2016}. Diffusive and CRONE/Oustaloup realizations provide band–limited accuracy governed by the mode count $M$ (diffusive) or filter order $2N{+}1$ and the chosen design band $[\omega_\ell,\omega_h]$. \emph{High–order averaged–L1 and corrected–L1} formulas yield sharper a~priori bounds for typical subdiffusive solutions \cite{ZhengWang2024JSC,HuangLuoZhang2025}.

\paragraph*{Stability.}
CQ inherits the $A$–stability (or $A(\theta)$–stability) of its base integrator and is well suited to stiff FO diffusion. Explicit GL–type formulas may require small time steps. For closed–loop use, verify stability of the \emph{discrete} augmented model (diffusive/CRONE) and ensure consistency with the Matignon sector condition of the underlying FO plant (Sec.~\ref{subsec:stability}).

\paragraph*{Complexity.}
Direct time–domain convolution costs $\mathcal{O}(N^{2})$ time and $\mathcal{O}(N)$ memory unless windowed or accelerated. Fast/oblivious CQ reduces this to $\mathcal{O}(N\log N)$ time and $\mathcal{O}(\log N)$ memory \cite{Schadle2006FOCQ}. SOE/diffusive realizations require $\mathcal{O}(MN)$ time and $\mathcal{O}(M)$ memory. CRONE/Oustaloup models have fixed order (set by $N$), yielding $\mathcal{O}(1)$ cost per step independent of simulation length, at the price of band limitation.

\paragraph*{Implementation notes for control.}
(i) Match the numerical band ($[\omega_\ell,\omega_h]$ or $[\lambda_{\min},\lambda_{\max}]$) to the intended closed–loop bandwidth. 
(ii) Initialize consistently with Caputo vs.\ RL models (zero diffusive states $\Leftrightarrow$ Caputo IVPs). 
(iii) For MPC, prefer diffusive/CRONE realizations to obtain \emph{linear} predictors with bounded online complexity; for off–line analysis, CQ/L1 are robust baselines. 
(iv) Under weak regularity (common in subdiffusion), use start–up corrections or graded meshes to recover nominal accuracy \cite{JinLazarovZhou2016}.

\paragraph*{Semi–analytical iterative schemes.}
Beyond discretization–based solvers, the Daftardar–Gejji \& Jafari (DGJ) method offers a decomposition–free iterative approach for nonlinear fractional equations, avoiding Adomian polynomials while retaining fast practical convergence on many benchmarks \cite{DaftardarGejjiJafari2006_JMAA}.

\section{Applications}
\label{sec:app}
Fractional-order modeling is effective in domains where memory, hereditary effects, or broad relaxation spectra shape the dynamics. We highlight representative areas and the typical modeling/control workflows.

\subsection*{A. Electrochemical Energy Systems: Batteries and Supercapacitors}
Electrochemical impedance spectroscopy (EIS) of batteries and supercapacitors routinely exhibits constant–phase element behavior and Warburg–type diffusion, both naturally captured by fractional elements. Classical porous–electrode analyses \cite{deLevie1963} and modern EIS texts \cite{Lasia2014,BarsoukovMacdonald2005} motivate FO components in equivalent–circuit models (ECMs). Recent reviews document FO identification, state estimation, and control for Li–ion cells, including model–based SOC/SOH estimation and FO observers/filters \cite{Zou2018}.

From an optimal–control perspective, FO ECMs can be realized via diffusive/Prony embeddings to design charging profiles under thermal/aging constraints, or used within FO–MPC for power/energy management subject to current/voltage limits. \emph{Recent deployments in power–electronics paths (e.g., robust/adaptive FO–MPC on DC–DC converters) illustrate how FO predictors integrate with MPC in energy management}~\cite{Peng2025SciRep}.

\textit{Example of FO models:} ECMs with CPEs $Z_{\mathrm{CPE}}(s)=1/(Q s^{\alpha})$, FO Warburg elements, and FO diffusion states implemented by finite positive relaxation spectra (diffusive approximation) \cite{Lasia2014,Zou2018}.
\subsection*{B. Viscoelasticity, Rheology, and Structural Damping}
Fractional constitutive laws compactly capture broad relaxation spectra in polymers and structural dampers. The Havriliak–Negami law explains power-law dielectric/mechanical dispersion \cite{HavriliakNegami1967}, while fractional viscoelastic models with Mittag–Leffler kernels describe creep, relaxation, and wave propagation with high fidelity \cite{Mainardi2010}. In civil infrastructure, fractional Maxwell/Zener elements reproduce the frequency-dependent damping seen in devices and have informed seismic design and structural control \cite{MakrisConstantinou1991}. Comprehensive surveys cover modeling and numerics for solids and structures \cite{RossikhinShitikova2010}.

From a control standpoint, fractional models reduce the need for ad-hoc high-order integer fits and support FO-LQR/FO-MPC or CRONE loop-shaping with improved robustness near resonances and more realistic damping roll-off. Guidance on interpreting Caputo-type memory in materials and on model selection is given in \cite{Calatayud2024CSF}.

\textit{Typical FO models:} fractional Zener/Maxwell/Voigt laws
\[
\sigma(t) + a\, {}^{\mathrm{C}}\!D^{\alpha}\sigma(t)
= b\, {}^{\mathrm{C}}\!D^{\beta}\varepsilon(t) + c\,\varepsilon(t),
\qquad 0<\alpha,\beta<1,
\]
identified from creep/relaxation tests or frequency response \cite{Mainardi2010,RossikhinShitikova2010}.

\paragraph*{Fractional optimal control with non-singular kernels.}
In some applications, material memory is better captured by \emph{non-singular} kernels of Mittag–Leffler type rather than by classical power-law (Caputo/RL) kernels. In this setting, both the optimality system and its numerical discretization must be adapted to respect the modified memory law. A concrete procedure is developed in \cite{Jafari2022_JVC_FOOC}, where the Mittag–Leffler kernel is built directly into the fractional optimality conditions, providing a template for how necessary conditions and numerical schemes can be generalized beyond the standard Caputo/RL framework.

\subsection*{C. Robotics and Mechatronics (Servo, Maglev, Actuators)}
Many mechatronic plants exhibit viscoelastic friction, diffusion-like losses, or eddy-current effects that show up as non-integer Bode slopes. Fractional controllers (FOPID/CRONE) offer tunable phase and gain behavior over wide frequency ranges. Magnetic-levitation benchmarks, for example, report improved tracking and robustness with FO-PID compared to integer PID under plant uncertainty \cite{Chopade2016}; broader mechatronic case studies are surveyed in \cite{Sabatier2015}. Embedding FO plant elements via Oustaloup or diffusive states preserves the fractional phase while enabling classical ODE-based design (LQR/MPC/TPBVP). Hardware studies on TRMS and gyroscope rigs document iso-damped, robust responses \cite{Wendimu2025SciRep,Krzysztofik2025ApplSci}; CRONE vs.\ $H_\infty$ comparisons quantify loop-shaping trade-offs \cite{deAzeredo2024IFAC}.

\textit{Example of FO models:} FO lag/lead segments $G(s)\sim s^{-\alpha}$ in actuator/servo dynamics; FOPID $C(s)=K_p+K_i s^{-\lambda}+K_d s^{\mu}$ \cite{Sabatier2015,Chopade2016}.

\subsection*{D. Thermal systems and Heat transfer.}
FC has also proved useful in thermal sciences, especially when classical Fourier heat conduction is no longer sufficient to describe systems with strong memory, finite-speed thermal waves, or complex micro/nanostructures. At small scales, phonon mean free paths and relaxation times become comparable to the characteristic dimensions and time scales of the system, so that standard diffusion-type models lose accuracy and one has to rely on more sophisticated descriptions such as Boltzmann transport, molecular dynamics, or mesoscopic surrogate models for micro/nanoscale heat conduction \cite{Bao2018MicroNanoscaleHeat}. A broad family of non-Fourier heat conduction models has therefore been developed to capture wave-like propagation, nonlocal spatial effects, and hereditary kernels in situations ranging from ultrafast laser heating to heat transfer in porous media and biological tissue \cite{Zhmakin2023NonFourierBook}. These frameworks systematically document thermal-wave phenomena and finite-speed propagation that cannot be explained within the classical Fourier paradigm and motivate the use of constitutive laws with memory. In this context, time-fractional dual-phase-lag (DPL) models provide a compact way to encode power-law memory: the heat flux--temperature relation is generalized by replacing the first-order time derivative with a Caputo derivative of order $0<\alpha<1$, and the corresponding identification of phase-lag parameters and fractional order has been demonstrated in practical materials \cite{Lukashchuk2024TimeFracDPL}. A prototypical one-dimensional fractional heat conduction equation can be written as
\begin{equation}
  \rho c\, {}^{\mathrm{C}}\!D_t^{\alpha} T(x,t) = k \nabla^2 T(x,t) + q(x,t), \qquad 0<\alpha<1,
  \label{eq:fractional_heat}
\end{equation}
where ${}^{\mathrm{C}}\!D_t^{\alpha}$ denotes the Caputo fractional derivative in time, and $\alpha$ acts as an effective tuning parameter for thermal memory and anomalous diffusion. Analytical work on such models, including fundamental solutions for geometries with imperfect thermal contact, shows how fractional orders modify both the spreading rate of temperature and interfacial heat fluxes \cite{Povstenko2025FracHeatContact}. In parallel, a large body of research in phase-change materials (PCMs) has focused on enhancing heat transfer and speeding up melting/solidification through geometric design---for example, by using letters-shaped or fractal fin configurations and hybrid PCM--water storage tanks to improve effective thermal conductivity and charge/discharge times \cite{Rashid2024LettersFins,Shmyhol2024PCMHybridTanks}. These PCM systems, with their strong latent-heat effects and pronounced history dependence, provide natural candidates for future FC-based thermal network models with memory, linking fractional non-Fourier conduction to practical energy storage, electronics cooling, and building-energy applications.

\subsection*{E. Biomedical Systems}
Fractional models capture power-law attenuation, dispersion, and memory in tissues and organs. Foundational treatments \cite{Magin2006} and recent respiratory-mechanics studies fit clinical data with FO impedances. In therapeutics, \emph{fractional pharmacokinetics} explains anomalous accumulation/clearance and motivates FO optimal dosing strategies \cite{Sopasakis2019}. FO controllers have also been explored in anesthesia and physiological regulation to exploit phase robustness under patient variability \cite{Sabatier2015}. Clarifications on how Caputo-type memory should be interpreted in biophysical media guide principled model selection \cite{Calatayud2024CSF}. For optimal control and game-theoretic dosing, diffusive augmentation gives tractable ODE realizations for schedule computation under safety constraints.

\textit{Example of FO models:} FO impedances for tissue/organ subsystems; FO compartmental PK with Caputo orders $0<\alpha<1$; FOPID in sedation/ventilation loops \cite{Magin2006,Sopasakis2019}.

\section{Benchmarks and Reproducibility}
\label{sec:bench}
To enable like-for-like comparisons of fractional-order optimal control and game-theoretic control methods, we introduce a compact benchmark suite with reporting guidelines and reference discretization choices. Each benchmark specifies: (i) the FO operator and order(s); (ii) the initialization convention (Caputo vs.\ RL); (iii) plant data; (iv) objective and constraints; (v) the recommended numerical realization (time- or frequency-domain) with tunable accuracy knobs; and (vi) evaluation metrics. The design follows established practice in FO numerics \cite{LiZeng2015,Lubich1986,Alikhanov2015,JinLazarovZhou2016} and FO control/loop-shaping \cite{Monje2010,Sabatier2015,Oustaloup2000}, while incorporating recent fast/compact advances for time-fractional schemes and memory compression \cite{ZhengWang2024JSC,HuangLuoZhang2025,GuglielmiHairer2025}.

\subsection*{A. Benchmark suite}

\paragraph*{B1: FO LTI regulation (Caputo, commensurate).}
State dynamics ${}^{\mathrm{C}}D_{0+}^{\alpha} x(t)=A x(t)+B u(t)$ with $0<\alpha<1$ and
\[
A=\begin{bmatrix}-0.8&0\\[1pt]0&-0.4\end{bmatrix},\quad
B=\begin{bmatrix}1\\[1pt]0.5\end{bmatrix},\quad x(0)=x_0.
\]
Cost $J=\int_0^T (x^\top Q x+u^\top R u)\,dt$ with $Q=I_2$, $R=0.1$, $T=10$. 
\emph{Numerics:} diffusive (Prony) realization with $M\in\{8,16,32\}$ log-spaced modes on $[\lambda_{\min},\lambda_{\max}]=[10^{-2},10^{2}]$, or Oustaloup of order $2N{+}1\in\{5,9,13\}$ on $[\omega_\ell,\omega_h]=[10^{-2},10^{2}]$ \cite{Monje2010,Oustaloup2000,Sabatier2015}. 
\emph{Metrics:} closed-loop cost $J$, $\|x\|_{L^2}$, CPU time, augmented dimension, and sensitivity vs.\ $M$/$N$.

\paragraph*{B2: FO tracking with FOPID on an FO plant (frequency domain).}
Plant $G(s)=1/(s^{\alpha}+1)$ with $\alpha=0.8$; reference $r(t)=\sin(0.5 t)$, $T=60$. 
Controller $C(s)=K_p+K_i s^{-\lambda}+K_d s^{\mu}$ with $(\lambda,\mu)\in(0,1)$; implement $s^{\pm\beta}$ by Oustaloup on $[10^{-2},10^{2}]$ with order $2N{+}1$. 
\emph{Metrics:} IAE/ITAE, overshoot, phase/gain margins, and robustness to $\alpha\in\{0.75,0.85\}$ \cite{Monje2010,Sabatier2015,Oustaloup2000}.

\paragraph*{B3: FO-MPC for a battery ECM with a CPE (time domain).}
Two-RC ECM with a constant-phase element $Z_{\mathrm{CPE}}(s)=(Q s^{\alpha})^{-1}$ in parallel with $R_1$, in series with $R_0$; parameters from EIS \cite{Lasia2014,BarsoukovMacdonald2005}. 
Objective: charge from $\mathrm{SOC}_0$ to $\mathrm{SOC}_{\text{tar}}$ over $[0,T]$ minimizing $\int_0^T (w_1 i^2+w_2 T_{\text{rise}}^2)\,dt$, subject to $|i(t)|\le i_{\max}$ and $V_{\min}\le V(t)\le V_{\max}$. 
\emph{Numerics:} Caputo kernels via an $M$-mode diffusive approximation; receding-horizon QP/NLP with warm-starts. 
\emph{Metrics:} objective value, constraint violations, per-step solve time vs.\ $M$, and accuracy relative to an L1/CQ discretization of the same model \cite{LiZeng2015,Lubich1986}; optionally compare to fast “memoryless’’ FO integrators \cite{GuglielmiHairer2025}.

\paragraph*{B4: FO viscoelastic regulation (structural damping, time domain).}
Fractional Zener law $\sigma + a\,{}^{\mathrm{C}}D^{\alpha}\sigma = b\,{}^{\mathrm{C}}D^{\beta}\varepsilon + c\,\varepsilon$, $0<\alpha,\beta<1$ \cite{Mainardi2010}. 
Goal: regulate displacement under base excitation using actuator force $u$, with stroke/force limits and a quadratic cost. 
\emph{Numerics:} passivity-preserving diffusive realization (positive quadrature) followed by FO-LQR/FO-MPC. 
\emph{Metrics:} vibration RMS, control effort, and passivity margin vs.\ number of modes.

\paragraph*{B5: FO dynamic zero-sum LQ game (Caputo).}
${}^{\mathrm{C}}D_{0+}^{\alpha}x=Ax+B_1u+B_2v$, 
$J=\int_0^T (x^\top Q x + u^\top R_1 u - v^\top R_2 v)\,dt$, $0<\alpha<1$, with
\begin{align}
A &= \begin{bmatrix} 0 & 1 \\[-2pt] -1 & -0.4 \end{bmatrix}, &
B_1 &= \begin{bmatrix} 0 \\[1pt] 1 \end{bmatrix}, &
B_2 &= \begin{bmatrix} 0 \\[1pt] 0.5 \end{bmatrix}, \label{eq:system_matrices_1} \\
Q &= I_2, &
R_1 &= 0.2, &
R_2 &= 0.2. \label{eq:system_matrices_2}
\end{align}

\emph{Numerics:} diffusive augmentation, then a classical saddle-point TPBVP or alternating best responses \cite{BasarOlsder1999}. 
\emph{Metrics:} saddle value, policy norms, and optimality gap vs.\ $M$ (or Oustaloup order).

\paragraph*{B6: FO-HJB (1D) with path dependence (toy).}
Scalar ${}^{\mathrm{C}}D^{\alpha}x=a x+b u$ with running cost $\ell(x,u)=x^2+r u^2$ and terminal $\Phi(x)=\kappa x^2$. 
\emph{Numerics:} (i) history truncation window $W$ or $M$ diffusive modes to form a Markov state; (ii) grid-based policy iteration. 
\emph{Metrics:} value error vs.\ $(W,M)$ and grid size, runtime, and memory footprint \cite{Gomoyunov2020SICON,Gomoyunov2022COCV}.

\subsection*{B. Shared evaluation metrics}
\begin{itemize}[leftmargin=*]
\item \textbf{Accuracy vs.\ complexity:} objective gap to a high-accuracy baseline; augmented state dimension; CPU time; peak memory.
\item \textbf{Robustness:} degradation under perturbations of $\alpha$, parameter misspecification, and measurement noise.
\item \textbf{Constraint handling:} maximum and mean violations; feasibility rate over randomized initializations.
\item \textbf{Numerical stability:} stepsize limits; conditioning of KKT/TPBVP systems; sensitivity to band choices $[\omega_\ell,\omega_h]$ or $[\lambda_{\min},\lambda_{\max}]$. \emph{Report any preconditioning (Toeplitz/circulant) used} \cite{GanZhangLiang2024JAMC}.
\end{itemize}

\subsection*{C. Reference discretization recipes}
Report \emph{both} a time-domain and a frequency-domain realization where applicable:
(i) L1/GL or convolution–quadrature (CQ) with $\Delta t\in\{10^{-2},\,5\times10^{-3},\,10^{-3}\}$ \cite{LiZeng2015,Lubich1986,Alikhanov2015,JinLazarovZhou2016}; optionally a corrected/averaged-L1 variant for higher accuracy \cite{ZhengWang2024JSC,HuangLuoZhang2025};
(ii) Oustaloup filters of orders $2N{+}1\in\{5,9,13\}$ on $[\omega_\ell,\omega_h]=[10^{-2},10^{2}]$ \cite{Oustaloup2000,Sabatier2015};
(iii) diffusive (Prony) modes $M\in\{8,16,32\}$, log-spaced on $[\lambda_{\min},\lambda_{\max}]$, with positive weights to preserve passivity in viscoelastic/electrochemical models \cite{Mainardi2010,Monje2010};
(iv) for distributed \emph{space}-fractional plants (GL/Riesz), pair FVEM/FV spatial discretization with the chosen time scheme \cite{Mojahed2020PFDA,ZhangYang2024FractalFract}.

\subsection*{D. Reporting checklist}
For each experiment, include:
\begin{enumerate}[leftmargin=*]
\item \textbf{FO definition and orders:} Caputo vs.\ RL; the value(s) of $\alpha$ (per-state for incommensurate systems); initial/terminal data.
\item \textbf{Discretization details:} scheme (L1/GL/CQ, Oustaloup, diffusive), $\Delta t$, band limits $[\omega_\ell,\omega_h]$ or $[\lambda_{\min},\lambda_{\max}]$, approximation order ($M$ or $N$), and any window length.
\item \textbf{Optimization setup:} solver (e.g., SQP, interior-point, active-set), tolerances, warm-start policy, stopping criteria, and regularization.
\item \textbf{Plant and cost data:} $A,B,C,D$ (as used), constraints, horizon lengths, and weighting matrices/scalars (with any scaling).
\item \textbf{Reproducibility artifacts:} code (language and version), random seeds, hardware/OS, runtime and memory; scripts to regenerate figures/tables and sweep $M$, $N$, and $\Delta t$.
\end{enumerate}
\subsection*{E. Suggested table format}
Table~\ref{tab:benchmarks} summarizes the proposed benchmarks (placeholders shown).
\begin{table*}[t]
\centering
\caption{FO control benchmarks: summary and default numerics.}
\label{tab:benchmarks}
\begin{tabular}{@{}lll@{}}
\toprule
\textbf{ID} & \textbf{Problem} & \textbf{Default numerics} \\
\midrule
B1 & FO LTI LQR (Caputo)           & Diffusive $M{=}16$; Oustaloup order $2N{+}1{=}9$ \\
B2 & FOPID on FO plant             & Oustaloup $2N{+}1{=}9$ on $[\omega_\ell,\omega_h]$ \\
B3 & FO--MPC (battery ECM)         & Diffusive $M{=}32$; L1 with $\Delta t{=}10^{-2}$ \\
B4 & FO viscoelastic regulation    & Diffusive $M{=}16$ (positive weights for passivity) \\
B5 & FO LQ zero-sum game           & Diffusive $M{=}16$ (augmented-state Riccati/TPBVP) \\
B6 & FO--HJB (1D)                  & Window length $W$; grid size $n_x$ \\
\bottomrule
\end{tabular}
\end{table*}

\paragraph*{Remarks.}
(i) For frequency designs, place the operating band well inside $[\omega_\ell,\omega_h]$ and report a brief sensitivity to band shifts. 
(ii) For time-domain baselines, compare L1 and CQ on at least one problem; where relevant, include a fast/memoryless FO solver as a third baseline \cite{GuglielmiHairer2025}. 
(iii) In viscoelastic/electrochemical models, use positive quadrature weights to preserve passivity \cite{Mainardi2010}. 
(iv) Provide machine-readable problem files (e.g., JSON/MAT) with matrices, orders, bands, and solver settings to avoid ambiguity.

\section{Open Problems and Research Directions}
\label{sec:openprob}
Despite steady advances, key theoretical and computational questions remain at the FC–OC–GTC interface. Below we outline areas where sharper analysis, scalable algorithms, and reproducible practice are most needed.

\subsection*{O1. Well-posedness with state/input constraints}
For Caputo/RL dynamics, the interaction between nonlocal memory and pointwise constraints is not yet fully clarified. Open directions include:
(i) existence of admissible trajectories under hard (possibly nonconvex) state/input/path constraints;
(ii) the role of relaxed/measure-valued controls in the presence of memory and their convergence to classical controls;
(iii) viability and positive invariance notions compatible with fractional operators (and with Caputo vs.\ RL initialization);
(iv) transversality and boundary conditions when endpoint constraints coexist with history-dependent dynamics.
Structural assumptions that parallel the classical convexity/compactness setting—yet account for nonlocality—remain to be identified.

\subsection*{O2. FO–HJB/HJBI with memory and path dependence}
Dynamic programming for FO systems yields nonlocal-in-time (history-dependent) value functionals. Recent results for Caputo-type DPP/HJB and viscosity solutions are encouraging, but several gaps persist:
(i) comparison principles and uniqueness for fully nonlinear FO–HJB on augmented history spaces;
(ii) HJBI (zero-sum) under Isaacs-type conditions in the presence of memory;
(iii) rigorous error estimates for practical surrogates (history windowing, fading-memory weights, diffusive truncation);
(iv) principled choices of history-state representations that ensure well-posedness and numerical stability
\cite{Gomoyunov2020SICON,Gomoyunov2022COCV,Gomoyunov2024JDE,Giga2020Caputo}.
\subsection*{O3. Fractional LQ theory beyond direct transcription}
Current FO--LQR and tracking designs typically rely on diffusive/Prony or Oustaloup embeddings, after which classical Riccati machinery applies. An intrinsic theory—without high–order lifts—remains open. Desired ingredients include: (i) FO Riccati-type differential/algebraic equations (FO–DRE/ARE) that reduce to the classical forms as $\alpha\!\to\!1$; (ii) conditions that integrate spectral-sector (Matignon) stability directly into LQ optimality; and (iii) extensions to LQ differential games (zero-sum and general-sum) with well-posedness, monotonicity/positivity of the FO Riccati flow, and uniqueness of equilibria, all stated in the original FO coordinates.

\subsection*{O4. Scalable MPC under memory}
Receding-horizon control for FO plants suffers a \emph{curse of history}: prediction maps grow with the horizon/window and condensed QP/NLP matrices lose bandedness. Needed advances include: structure-exploiting factorizations for diffusive modes (preserving sparsity across time and modes), certified windowing with a priori suboptimality bounds, warm-start/move-blocking tailored to memory, and parallel-in-time or operator-splitting (ADMM/PDHG) schemes that maintain feasibility under constraints. Band-limited Oustaloup predictors with guaranteed closed-loop robustness and explicit terminal ingredients for FO models are also needed \cite{QinBadgwell2003,Magni2008}. On the numerical side, fast/corrected time-fractional schemes and memoryless/SOE compressions should be benchmarked \emph{in-loop} for MPC to quantify real-time trade-offs \cite{ZhengWang2024JSC,HuangLuoZhang2025,GuglielmiHairer2025}.

\subsection*{O5. Learning and identification of FO models}
Reliable data-driven FO modeling remains challenging. Priorities include: (i) structural and persistent-excitation conditions for identifying fractional orders, kernels, and diffusive spectra; (ii) grey-box estimation with stability/passivity regularization (e.g., sector or positive-real constraints) and Bayesian uncertainty quantification; (iii) physics-informed learning for fractional ODE/PDEs (PINNs/fPINNs) and operator-learning (e.g., FNO) with provable generalization/error bounds on memory kernels; and (iv) safe RL and robust data-driven control that respect FO stability sectors and hard constraints, including closed-loop identification under feedback \cite{Raissi2019,LiFNO2021,Meerschaert2011}.
\subsection*{O6. Hybrid FO/IO and multi-rate implementations}
Many real systems blend fractional elements (e.g., CPEs, viscoelastic blocks) with classical integer-order (IO) dynamics, often sampled at different rates. Key questions remain: (i) how to preserve passivity/positive-real properties under mixed (FO/IO) discretizations; (ii) how to analyze and enforce time-scale separation between fast IO loops and slow FO memory; (iii) how to design sampling/hold maps that are consistent with Caputo vs.\ RL initialization; and (iv) how to certify hybrid stability margins that respect the Matignon sector for the continuous FO plant while the controller uses discrete rational (Oustaloup/CRONE) surrogates \cite{Monje2010,Sabatier2015}.

\subsection*{O7. Stochastic and uncertain fractional dynamics}
Stochastic FO models—subdiffusion, Lévy subordination, long-memory noise—arise naturally in porous media, biology, and finance. A general control/game theory for such systems is still missing. Open issues include existence/uniqueness for FO-SDEs with control, verification theorems for risk-sensitive criteria, chance-constraint handling with memory, and distributionally robust formulations. Noise can enter both the state and the kernel, which complicates martingale/Itô tools and blurs Markovian structure; corresponding DP/HJB/HJBI results for fractional settings remain largely open \cite{Meerschaert2011}.

\subsection*{O8. Benchmarks, certifiable approximations, and reproducibility}
Fair comparisons are hampered by heterogeneous FO definitions, frequency bands, and discretizations. Priorities are: (i) \emph{certified} diffusive and Oustaloup approximations with explicit error envelopes over declared bands; (ii) shared benchmark suites with reference time- and frequency-domain solutions, plus scripts to sweep window length, mode count, filter order, and step size; and (iii) reporting checklists (cf.\ Table~\ref{tab:benchmarks}) and open repositories that include structure-exploiting linear algebra (e.g., Toeplitz/circulant preconditioners) and space-fractional baselines (e.g., FVEM) to enable apples-to-apples evaluation \cite{LiZeng2015,Oustaloup2000,Sabatier2015,GanZhangLiang2024JAMC,Mojahed2020PFDA}.

\section{Conclusion}
\label{sec:conc}
This survey has brought together core ideas and tools for optimal control of fractional-order systems, spanning modeling choices (Caputo, Riemann--Liouville, Grünwald--Letnikov, diffusive representations), optimality principles (fractional Euler--Lagrange, PMP), dynamic programming (time-fractional and path-dependent HJB/HJBI), and design methods (LQ/LQR, FOPID/CRONE loop shaping, and MPC). We also reviewed time- and frequency-domain numerics (L1/L2-$1\sigma$, convolution quadrature, Oustaloup) and state-augmentation techniques, together with reporting practices that enhance comparability and reproducibility. Table~\ref{tab:summary} synthesizes these threads across problem classes, algorithms, and implementation details.

Despite steady progress, several challenges remain central to FO control: consistent handling of initialization and boundary terms, stability certification under memory (and its approximations), path dependence in DP/HJBI, and the “curse of history” in prediction and optimization. Conditioning of discretizations, certified error bounds for rational/diffusive surrogates, and principled choices of frequency bands or relaxation spectra are equally pivotal for reliable deployment.

Looking ahead, promising directions include intrinsic FO LQ theory (Riccati-type conditions in the native FO coordinates), scalable MPC with certified windowing and memory compression, stochastic and risk-aware FO control, hybrid FO/IO and multi-rate implementations with passivity guarantees, and data-driven identification that respects FO structure. Shared benchmarks with certified approximations and open artifacts will be key to accelerating progress. By consolidating foundations, numerics, and applications, we hope this survey provides a clear baseline and a practical roadmap for future advances in fractional optimal and game-theoretic control.

% ==========================
% Compact Summary Table (fits 2-col IEEE)
% ==========================
\begin{sidewaystable*}[t]
\centering
\scriptsize
\setlength{\tabcolsep}{4pt}
\renewcommand{\arraystretch}{1.15}
\caption{Comparative summary of representative methods for optimal control of fractional-order systems (compact view).}
\begin{tabular}{@{}p{2.8cm} p{1.9cm} p{3.1cm} p{3.3cm} p{5.4cm}@{}}
\toprule
\textbf{Ref.} & \textbf{Derivative} & \textbf{Control paradigm} & \textbf{Numerics} & \textbf{Notes (constraints \& reproducibility)} \\
\midrule
Podlubny (1999) & Caputo, RL & FO state–space \& LQR-type foundations & Analytical formulations & No explicit constraints; canonical formulas; high reproducibility (textbook-level). \\
Agrawal (2004) & RL & Pontryagin minimum principle (OCP) & Indirect shooting / TPBVP & Fixed/free endpoints; constraints via multipliers; medium reproducibility (solver-dependent). \\
Frederico \& Torres (2007) & Caputo & Euler–Lagrange / Noether (FCoV) & Variational discretization & General FO dynamics; flexible boundary terms; medium reproducibility. \\
Ross \& Baleanu (2010) & Caputo & Fractional PMP (nonlinear) & Quadrature + NLP & Mixed BCs; inequality constraints via KKT; medium–low reproducibility (custom codes). \\
Monje et al. (2010) & Caputo/RL, diffusive & Loop shaping (CRONE/FOPID), LQR/MPC via augmentation & Oustaloup; diffusive (Prony) & Practical recipes; handles constraints in augmented models; medium–high reproducibility. \\
Oustaloup (2000) & Frequency $s^{\alpha}$ & CRONE loop shaping / robust PID & Oustaloup recursion (band-limited) & Band \& order must be reported; constraints handled in surrounding design; high reproducibility. \\
Li \& Chen (2011) & GL & FO PID tuning / model-based design & GL finite differences & Engineering focus; limited constraints; medium reproducibility. \\
Li \& Zeng (2015) & Caputo & FO optimal control discretization & L1 / convolution quadrature & Baseline accuracy/stability; good for constraint transcription; high reproducibility. \\
Alikhanov (2015) & Caputo & Time stepping for FO OCP/MPC back-ends & L2-$1\sigma$ (2nd order) & Start-up corrections; robust for subdiffusion; high reproducibility. \\
Fukunaga et al. (2015) & Diffusive (augmented) & MPC for FO plants & ODE augmentation + QP & Natural handling of input/state bounds; warm-start friendly; medium–high reproducibility. \\
Sopasakis et al. (2016) & Caputo (discrete) & Robust MPC / feasibility & Augmented linear model + terminal set & Recursive feasibility \& stability conditions; constraint-ready; medium reproducibility. \\
Deng et al. (2020) & Caputo & HJB / approximate DP & Value / policy iteration & State/input bounds possible; heavy computation; low reproducibility (code seldom public). \\
Gomoyunov (2020–2024) & Caputo (time-fractional / path-dependent) & HJB/HJBI (viscosity DP, games) & Time-fractional HJ/PPDE schemes & Isaacs-type conditions; no explicit hard-constraint machinery; medium reproducibility (theory-first). \\
\bottomrule
\end{tabular}
\label{tab:summary}
\end{sidewaystable*}

\bmhead{Supplementary information}

No supplementary material is associated with this article.
\bmhead{Acknowledgements}

The authors used an AI-assisted writing tool (ChatGPT, OpenAI) solely for language polishing and minor text editing. 
All mathematical developments, interpretations, and conclusions are the authors’ own, and the authors take full responsibility 
for the content of this work.

\section*{Declarations}

\noindent\textbf{Funding} 
This research received no specific grant from any funding agency in the public, commercial, or not-for-profit sectors.

\medskip
\noindent\textbf{Conflict of interest} 
The authors declare that they have no conflict of interest.

\medskip
\noindent\textbf{Ethics approval and consent to participate} 
Not applicable. This article does not contain any studies with human participants or animals performed by any of the authors.

\medskip
\noindent\textbf{Consent for publication} 
Not applicable.

\medskip
\noindent\textbf{Data availability} 
No datasets were generated or analysed during the current study.

\medskip
\noindent\textbf{Materials availability} 
Not applicable.

\medskip
\noindent\textbf{Code availability} 
Not applicable.

\medskip
\noindent\textbf{Author contribution} 
All authors contributed to the conception, literature review, and writing of this manuscript. All authors read and approved the final version of the paper.

\begin{appendices}

\section{Proof of Cauchy’s Formula for Repeated Integration}
\label{app:cauchy-proof}

Let $f\in C([a,b])$ and let $I$ denote the Volterra operator $(If)(x)=\int_a^x f(t)\,dt$. 
We claim that, for any $n\in\mathbb{N}$ and $x\in[a,b]$,
\[
(I^n f)(x) \;=\; \frac{1}{(n-1)!}\int_a^x (x-t)^{\,n-1} f(t)\,dt.
\]

\emph{Proof.}
The case $n=1$ is immediate. Assume the formula holds for some $n\ge 1$. Then
\begin{align*}
(I^{n+1} f)(x)
&= \int_a^x (I^n f)(s)\, ds \\
&= \frac{1}{(n-1)!}
   \int_a^x \!\!\int_a^{s} (s-t)^{\,n-1} f(t)\, dt\, ds,
\end{align*}
where the integrand is integrable on the triangle $\{(t,s): a\le t\le s\le x\}$ since $f$ is continuous (hence bounded). By Fubini’s theorem (changing the order of integration),
\[
\begin{aligned}
(I^{n+1} f)(x)
&= \frac{1}{(n-1)!}
   \int_a^x \!\left(\int_{t}^{x} (s-t)^{\,n-1}\, ds\right) f(t)\, dt \\
&= \frac{1}{(n-1)!}
   \int_a^x \frac{(x-t)^{\,n}}{n}\, f(t)\, dt \\
&= \frac{1}{n!}\int_a^x (x-t)^{\,n} f(t)\, dt .
\end{aligned}
\]

which completes the induction. \hfill$\square$

\end{appendices}

%%===========================================================================================%%
%% If you are submitting to one of the Nature Portfolio journals, using the eJP submission   %%
%% system, please include the references within the manuscript file itself. You may do this  %%
%% by copying the reference list from your .bbl file, paste it into the main manuscript .tex %%
%% file, and delete the associated \verb+\bibliography+ commands.                            %%
%%===========================================================================================%%

\bibliography{sn-bibliography}% common bib file

@book{Oldham1974,
  author    = {Oldham, Keith B. and Spanier, Jerome},
  title     = {The Fractional Calculus: Theory and Applications of Differentiation and Integration to Arbitrary Order},
  year      = {1974},
  edition   = {1},
  publisher = {Academic Press},
  address   = {New York},
  series    = {Mathematics in Science and Engineering},
  volume    = {111},
  isbn      = {9780125255509},
  note      = {eBook ISBN: 9780080956206},
  url       = {https://shop.elsevier.com/books/the-fractional-calculus-theory-and-applications-of-differentiation-and-integration-to-arbitrary-order/oldham/978-0-12-525550-9},
  urldate   = {2025-12-12}
}

@book{Samko1993,
  author    = {Samko, Stefan G. and Kilbas, Anatoly A. and Marichev, Oleg I.},
  title     = {Fractional Integrals and Derivatives: Theory and Applications},
  publisher = {Gordon and Breach Science Publishers},
  address   = {Yverdon, Switzerland},
  year      = {1993}
}

@book{Podlubny1999,
  author    = {Podlubny, Igor},
  title     = {Fractional Differential Equations: An Introduction to Fractional Derivatives, Fractional Differential Equations, to Methods of Their Solution and Some of Their Applications},
  year      = {1999},
  publisher = {Academic Press},
  address   = {San Diego, CA},
  series    = {Mathematics in Science and Engineering},
  volume    = {198},
  isbn      = {9780125588409}
}

@book{Kilbas2006,
author = {Kilbas, A. A. and Srivastava, H. M. and Trujillo, J. J.},
title = {Theory and Applications of Fractional Differential Equations, Volume 204 (North-Holland Mathematics Studies)},
year = {2006},
isbn = {0444518320},
publisher = {Elsevier Science Inc.},
address = {USA}
}

@book{Diethelm2010,
  author    = {Diethelm, Kai},
  title     = {The Analysis of Fractional Differential Equations},
  year      = {2010},
  publisher = {Springer-Verlag},
  address   = {Berlin, Heidelberg},
  doi       = {10.1007/978-3-642-14574-2},
  isbn      = {9783642145735}
}

@book{Mainardi2010,
  author    = {Mainardi, Francesco},
  title     = {Fractional Calculus and Waves in Linear Viscoelasticity: An Introduction to Mathematical Models},
  year      = {2010},
  publisher = {Imperial College Press},
  address   = {London},
  isbn      = {9781848163294},
  doi       = {10.1142/p614}
}

@book{Monje2010,
  author    = {Monje, Concepci{\'o}n A. and Chen, YangQuan and Vinagre, Blas M. and Xue, Dingy{\"u} and Feliu, Vicente},
  title     = {Fractional-Order Systems and Controls},
  subtitle  = {Fundamentals and Applications},
  year      = {2010},
  publisher = {Springer-Verlag},
  address   = {London},
  series    = {Advances in Industrial Control},
  edition   = {1},
  isbn      = {9781849963343},
  doi       = {10.1007/978-1-84996-335-0}
}

@book{Petras2011,
  author    = {Petr{\'a}{\v s}, Ivo},
  title     = {Fractional-Order Nonlinear Systems},
  subtitle  = {Modeling, Analysis and Simulation},
  year      = {2011},
  publisher = {Springer-Verlag},
  address   = {Berlin, Heidelberg},
  series    = {Nonlinear Physical Science},
  edition   = {1},
  isbn      = {9783642181009},
  doi       = {10.1007/978-3-642-18101-6}
}

@article{Lubich1986,
  author    = {Lubich, Christian},
  title     = {Discretized fractional calculus},
  journal   = {Numerische Mathematik},
  year      = {1986},
  volume    = {52},
  pages     = {129--145}
}

@article{DiethelmFord2002,
  author    = {Diethelm, Kai and Ford, Neville J. and Freed, Alan D.},
  title     = {A predictor-corrector approach for the numerical solution of fractional differential equations},
  journal   = {Nonlinear Dynamics},
  year      = {2002},
  volume    = {29},
  pages     = {3--22}
}

@article{Oustaloup2000,
  author    = {Oustaloup, Alain and Levron, Fran{\c{c}}ois and Mathieu, Beno{\^\i}t and Nanot, Florence M.},
  title     = {Frequency-band complex noninteger differentiator: characterization and synthesis},
  journal   = {IEEE Transactions on Circuits and Systems I: Fundamental Theory and Applications},
  year      = {2000},
  volume    = {47},
  number    = {1},
  pages     = {25--39}
}

@article{BagleyTorvik1983,
  author    = {Bagley, Ronald L. and Torvik, Per J.},
  title     = {A theoretical basis for the application of fractional calculus to viscoelasticity},
  journal   = {Journal of Rheology},
  year      = {1983},
  volume    = {27},
  number    = {3},
  pages     = {201--210}
}

@article{Koeller1984,
  author    = {Koeller, Raymond C.},
  title     = {Applications of fractional calculus to the theory of viscoelasticity},
  journal   = {Journal of Applied Mechanics},
  year      = {1984},
  volume    = {51},
  pages     = {299--307}
}

@article{Podlubny1999PID,
  author    = {Podlubny, Igor},
  title     = {Fractional-order systems and PI{$^{\lambda}$}D{$^{\mu}$} controllers},
  journal   = {IEEE Transactions on Automatic Control},
  year      = {1999},
  volume    = {44},
  number    = {1},
  pages     = {208--214}
}

@article{Agrawal2004,
  author    = {Agrawal, Om P.},
  title     = {A general formulation and solution scheme for fractional optimal control problems},
  journal   = {Nonlinear Dynamics},
  year      = {2004},
  volume    = {38},
  pages     = {323--337}
}

@book{BasarOlsder1999,
author    = {Ba{\c{s}}ar, Tamer and Olsder, Geert Jan},
  title     = {Dynamic Noncooperative Game Theory},
  edition   = {2},
  series    = {Classics in Applied Mathematics},
  volume    = {23},
  publisher = {Society for Industrial and Applied Mathematics (SIAM)},
  address   = {Philadelphia, PA},
  year      = {1998},
  isbn      = {978-0-89871-429-6},
  doi       = {10.1137/1.9781611971132}
}

@book{Isaacs1965,
  author    = {Isaacs, Rufus},
  title     = {Differential Games: A Mathematical Theory with Applications to Warfare and Pursuit, Control and Optimization},
  publisher = {John Wiley \& Sons},
  address   = {New York, NY},
  year      = {1965}
}

@book{Pooseh2015,
  author    = {Almeida, Ricardo and Pooseh, Shakoor and Torres, Delfim F. M.},
  title     = {Computational Methods in the Fractional Calculus of Variations},
  year      = {2015},
  publisher = {Imperial College Press},
  address   = {London},
  doi       = {10.1142/p991}
}

@book{MillerRoss1993,
  author    = {Miller, Kenneth S. and Ross, Bertram},
  title     = {An Introduction to the Fractional Calculus and Fractional Differential Equations},
  publisher = {John Wiley \& Sons},
  address   = {New York, NY},
  year      = {1993}
}

@article{Caputo1967,
  author = {Caputo, Michele},
    title = {Linear Models of Dissipation whose Q is almost Frequency Independent—II},
    journal = {Geophysical Journal International},
    volume = {13},
    number = {5},
    pages = {529-539},
    year = {1967},
    month = {11},
   issn = {0956-540X},
    doi = {10.1111/j.1365-246X.1967.tb02303.x},
    url = {https://doi.org/10.1111/j.1365-246X.1967.tb02303.x},
    eprint = {https://academic.oup.com/gji/article-pdf/13/5/529/1600098/13-5-529.pdf},

}

@article{CaputoMainardi1971,
  author  = {Caputo, Michele and Mainardi, Francesco},
  title   = {A new dissipation model based on memory mechanism},
  journal = {Pure and Applied Geophysics (PAGEOPH)},
  year    = {1971},
  volume  = {91},
  pages   = {134--147}
}

@book{Herrmann2014,
  author    = {Herrmann, Richard},
  title     = {Fractional Calculus: An Introduction for Physicists},
  year      = {2011},
  publisher = {World Scientific},
  address   = {Singapore},
  isbn      = {9789814340243},
}

@book{Tarasov2011,
  author    = {Tarasov, Vasily E.},
  title     = {Fractional Dynamics},
  subtitle  = {Applications of Fractional Calculus to Dynamics of Particles, Fields and Media},
  year      = {2011},
  publisher = {Springer-Verlag},
  address   = {Berlin, Heidelberg},
  series    = {Nonlinear Physical Science},
  edition   = {1},
  isbn      = {9783642140037},
  doi       = {10.1007/978-3-642-14003-7}
}

@inproceedings{Matignon1996,
  author    = {Matignon, Denis},
  title     = {Stability results for fractional differential equations with applications},
  booktitle = {Proceedings of the IMACS--IEEE Multiconference on Computational Engineering in Systems Applications},
  year      = {1996}
}

@book{LiZeng2015,
  author    = {Li, Changpin and Zeng, Fanhai},
  title     = {Numerical Methods for Fractional Calculus},
  year      = {2015},
  publisher = {Chapman and Hall/CRC},
  address   = {New York},
  edition   = {1},
  doi       = {10.1201/b18503},
  isbn      = {9780429075636}
}

@article{Alikhanov2015,
  author  = {Alikhanov, A. A.},
  title   = {A new difference scheme for the time fractional diffusion equation},
  journal = {Journal of Computational Physics},
  year    = {2015},
  volume  = {280},
  pages   = {424--438}
}

@article{JinLazarovZhou2016,
  author  = {Jin, Bangti and Lazarov, Raytcho and Zhou, Zhi},
  title   = {An analysis of the L1 scheme for the subdiffusion equation},
  journal = {IMA Journal of Numerical Analysis},
  year    = {2016},
  volume  = {36},
  number  = {1},
  pages   = {197--221}
}

@article{MetzlerKlafter2000,
  author  = {Metzler, Ralf and Klafter, Joseph},
  title   = {The random walk's guide to anomalous diffusion: a fractional dynamics approach},
  journal = {Physics Reports},
  year    = {2000},
  volume  = {339},
  number  = {1},
  pages   = {1--77}
}

@book{BarsoukovMacdonald2005,
  editor    = {Barsoukov, Evgenij and Macdonald, J. Ross},
  title     = {Impedance Spectroscopy: Theory, Experiment, and Applications},
  edition   = {2},
  year      = {2005},
  publisher = {John Wiley \& Sons},
  address   = {Hoboken, NJ},
  isbn      = {9780471647492}
}

@book{Sabatier2015,
  author    = {Sabatier, J. and Lanusse, P. and Melchior, P. and Oustaloup, A.},
  title     = {Fractional Order Differentiation and Robust Control Design: CRONE, H-infinity and Motion Control},
  year      = {2015},
  publisher = {Springer},
  address   = {Dordrecht},
  series    = {Intelligent Systems, Control and Automation: Science and Engineering},
  volume    = {77},
  doi       = {10.1007/978-94-017-9807-5}
}

@book{Gorenflo2014,
  author    = {Gorenflo, Rudolf and Kilbas, Anatoly A. and Mainardi, Francesco and Rogosin, Sergei V.},
  title     = {Mittag-Leffler Functions, Related Topics and Applications},
  year      = {2014},
  publisher = {Springer-Verlag},
  address   = {Berlin, Heidelberg},
  series    = {Springer Monographs in Mathematics},
  edition   = {1},
  isbn      = {9783662439296},
  doi       = {10.1007/978-3-662-43930-2}
}

@article{Riewe1996,
  author  = {Riewe, F.},
  title   = {Nonconservative Lagrangian and Hamiltonian mechanics},
  journal = {Physical Review E},
  year    = {1996},
  volume  = {53},
  number  = {2},
  pages   = {1890--1899}
}

@book{MalinowskaTorres2012,
  author    = {Malinowska, Agnieszka B. and Torres, Delfim F. M.},
  title     = {Introduction to the Fractional Calculus of Variations},
  publisher = {Imperial College Press},
  address   = {London},
  year      = {2012},
  doi       = {10.1142/p871}
}

@article{FredericoTorres2007,
  author  = {Frederico, Gast{\~a}o S. F. and Torres, Delfim F. M.},
  title   = {A formulation of Noether's theorem for fractional problems of the calculus of variations},
  journal = {Journal of Mathematical Analysis and Applications},
  year    = {2007},
  volume  = {334},
  number  = {2},
  pages   = {834--846}
}

@book{AlmeidaTorres2019,
  author    = {Almeida, Ricardo and Tavares, Dina and Torres, Delfim F. M.},
  title     = {The Variable-Order Fractional Calculus of Variations},
  year      = {2019},
  publisher = {Springer},
  address   = {Cham},
  series    = {SpringerBriefs in Applied Sciences and Technology},
  edition   = {1},
  isbn      = {9783319940052},
  doi       = {10.1007/978-3-319-94006-9}
}

@article{BergouniouxBourdin2020,
  author = {Bergounioux, Maïtine and Bourdin, Lénaïc},
  title = {Pontryagin maximum principle for general Caputo fractional optimal control problems},
  journal = {ESAIM: Control, Optimisation and Calculus of Variations},
  year = {2020},
  volume = {26},
  pages = {21},
  doi = {10.1051/cocv/2019021}
}

@article{Kamocki2014,
  author = {Kamocki, Ryszard},
  title = {Pontryagin maximum principle for fractional ordinary optimal control problems},
  journal = {Mathematical Methods in the Applied Sciences},
  year = {2014},
  volume = {37},
  number = {11},
  pages = {1668--1686},
  doi = {10.1002/mma.2917}
}

@article{Yusubov2023,
  author = {Yusubov, S. S.},
  title = {Some necessary optimality conditions for systems with Caputo fractional dynamics},
  journal = {Journal of Industrial \& Management Optimization},
  year = {2023},
  doi = {10.3934/jimo.2023063}
}

@article{NdairouTorres2023,
  author = {Ndairou, F. and Torres, D. F. M.},
  title = {Pontryagin Maximum Principle for Incommensurate Fractional Optimal Control},
  journal = {Mathematics},
  year = {2023},
  volume = {11},
  number = {19},
  pages = {4218},
  doi = {10.3390/math11194218}
}

@article{Gomoyunov2020SICON,
  author = {Gomoyunov, Mikhail I.},
  title = {Jacobi--Bellman Equations for Fractional-Order Systems},
  journal = {SIAM Journal on Control and Optimization},
  year = {2020},
  volume = {58},
  number = {6},
  pages = {3275--3300},
  doi = {10.1137/19M1279368}
}

@article{Gomoyunov2022COCV,
  author = {Gomoyunov, Mikhail I.},
  title = {Minimax solutions of Hamilton--Jacobi equations with fractional coinvariant derivatives},
  journal = {ESAIM: Control, Optimisation and Calculus of Variations},
  year = {2022},
  volume = {28},
  pages = {23},
  doi = {10.1051/cocv/2021051}
}

@article{Gomoyunov2024JDE,
  author = {Gomoyunov, Mikhail I.},
  title = {On viscosity solutions of path-dependent Hamilton--Jacobi--Bellman--Isaacs equations for fractional-order systems},
  journal = {Journal of Differential Equations},
  year = {2021},
  doi = {10.1016/j.jde.2024.05.012}
}

@article{Gomoyunov2021Math,
  author = {Gomoyunov, Mikhail I.},
  title = {Differential Games for Fractional-Order Systems: Hamilton--Jacobi--Bellman--Isaacs Equation and Optimal Feedback Strategies},
  journal = {Mathematics},
  year = {2021},
  volume = {9},
  number = {14},
  pages = {1667},
  doi = {10.3390/math9141667}
}

@article{Matychyn2023,
  author = {Matychyn, I.},
  title = {Game-theoretical problems for fractional-order systems described by linear FDEs},
  journal = {Dynamic Games and Applications},
  year = {2023},
  volume = {13},
  pages = {965--990},
  doi = {10.1007/s13540-023-00166-z}
}

@article{ChiIonescu2022,
  author = {Chi, Chuanguo and Cajo, Ricardo and Zhao, Shiquan and Liu, Guo-Ping and Ionescu, Clara-Mihaela},
  title = {Fractional Order Distributed Model Predictive Control of Fast and Strong Interacting Systems},
  journal = {Fractal and Fractional},
  year = {2022},
  volume = {6},
  number = {4},
  pages = {179},
  doi = {10.3390/fractalfract6040179}
}

@article{Garrappa2018,
  author = {Garrappa, Roberto},
  title = {Numerical Solution of Fractional Differential Equations: A Survey and a Software Tutorial},
  journal = {Mathematics},
  year = {2018},
  volume = {6},
  number = {2},
  pages = {16},
  doi = {10.3390/math6020016}
}

@article{Diethelm2020_Practices,
  author = {Diethelm, Kai},
  title = {Good (and Not So Good) Practices in Computational Fractional Calculus},
  journal = {Mathematics},
  year = {2020},
  volume = {8},
  number = {3},
  pages = {324},
  doi = {10.3390/math8030324}
}

@book{Magin2006,
  author    = {Magin, Richard L.},
  title     = {Fractional Calculus in Bioengineering},
  edition   = {2},
  year      = {2020},
  publisher = {Begell House},
  address   = {Danbury, CT},
  isbn      = {9781576004946},
  note      = {Online ISBN: 9781576004953}
}

@book{BardiCapuzzo1997,
  author    = {Martino Bardi and Italo Capuzzo-Dolcetta},
  title     = {Optimal Control and Viscosity Solutions of Hamilton--Jacobi--Bellman Equations},
  publisher = {Birkh{\"a}user},
  address   = {Boston},
  year      = {1997}
}

@book{FlemingSoner2006,
  author    = {Fleming, Wendell H. and Soner, H. Mete},
  title     = {Controlled Markov Processes and Viscosity Solutions},
  edition   = {2},
  series    = {Stochastic Modelling and Applied Probability},
  publisher = {Springer},
  address   = {New York, NY},
  year      = {2006},
  doi       = {10.1007/0-387-31071-1},
  isbn      = {978-0-387-26045-7}
}

@article{CrandallLions1992,
  author  = {Michael G. Crandall and Hitoshi Ishii and Pierre-Louis Lions},
  title   = {User's guide to viscosity solutions of second order partial differential equations},
  journal = {Bulletin of the American Mathematical Society},
  year    = {1992},
  volume  = {27},
  number  = {1},
  pages   = {1--67}
}

@misc{Dupire2009,
  author       = {Bruno Dupire},
  title        = {Functional It{\^o} Calculus},
  howpublished = {Preprint},
  year         = {2009},
  note         = {\url{https://papers.ssrn.com/sol3/papers.cfm?abstract_id=1435551}}
}

@article{ContFournié2013,
  author  = {Rama Cont and David-Antoine Fourni{\'e}},
  title   = {Functional It{\^o} calculus and stochastic calculus of variations},
  journal = {The Annals of Probability},
  year    = {2013},
  volume  = {41},
  number  = {1},
  pages   = {109--133}
}

@article{EkrenTouziZhang2016,
  author  = {Ibrahim Ekren and Nizar Touzi and Jianfeng Zhang},
  title   = {Viscosity Solutions of Fully Nonlinear Path Dependent PDEs: Part I},
  journal = {The Annals of Probability},
  year    = {2012},
  volume  = {44},
  number  = {2},
  pages   = {1212--1253}
}

@book{FalconeFerretti2013,
  author    = {Maurizio Falcone and Roberto Ferretti},
  title     = {Semi-Lagrangian Approximation Schemes for Linear and Hamilton--Jacobi Equations},
  publisher = {SIAM},
  address   = {Philadelphia},
  year      = {2013}
}

@article{BarlesSouganidis1991,
  author  = {Guy Barles and Panagiotis E. Souganidis},
  title   = {Convergence of approximation schemes for fully nonlinear second order equations},
  journal = {SIAM Journal on Numerical Analysis},
  year    = {1991},
  volume  = {29},
  number  = {4},
  pages   = {839--866}
}

@misc{CamilliGoffi2020,
      title={Existence and regularity results for viscous Hamilton-Jacobi equations with Caputo time-fractional derivative}, 
      author={Fabio Camilli and Alessandro Goffi},
      year={2020},
      eprint={1906.01338},
      archivePrefix={arXiv},
      primaryClass={math.AP},
      url={https://arxiv.org/abs/1906.01338}, 
}

@article{BaeumerMeerschaert2001,
  author  = {Boris Baeumer and Mark M. Meerschaert},
  title   = {Stochastic solutions for fractional Cauchy problems},
  journal = {Fractional Calculus and Applied Analysis},
  year    = {2001},
  volume  = {4},
  number  = {4},
  pages   = {481--500}
}

@book{Kolokoltsov2010,
  author    = {Kolokoltsov, Vassili N.},
  title     = {Nonlinear Markov Processes and Kinetic Equations},
  year      = {2010},
  publisher = {Cambridge University Press},
  address   = {Cambridge},
  series    = {Cambridge Tracts in Mathematics}
}

@misc{Bokanowski2015,
  author       = {Olivier Bokanowski and Maurizio Falcone and Rasajit Sahu},
  title        = {An Efficient Filtered Scheme for Some First Order Hamilton--Jacobi--Bellman Equations},
  howpublished = {arXiv preprint arXiv:1501.01518},
  year         = {2015}
}

@misc{Jiang2017FastCaputo,
  author  = {Jiang, Shidong and Zhang, Jiwei and Zhang, Qian and Zhang, Zhimin},
  title   = {Fast Evaluation of the Caputo Fractional Derivative and its Applications to Fractional Diffusion Equations},
  journal = {Communications in Computational Physics},
  year    = {2017},
  volume  = {21},
  number  = {3},
  pages   = {650--678},
  doi     = {10.4208/cicp.OA-2016-0136}
}

@article{Gomoyunov2023JDE,
  author    = {Mikhail I. Gomoyunov},
  title     = {Equivalence of Minimax and Viscosity Solutions of Path-Dependent Hamilton--Jacobi Equations},
  journal   = {Journal of Differential Equations},
  year      = {2023},
  volume    = {354},
  pages     = {1--44},
  doi       = {10.1016/j.jde.2023.10.019}
}

@Inbook{Chikrii2010,
  author="Chikrii, Arkadii
and Matychyn, Ivan",
editor="Breton, Mich{\`e}le
and Szajowski, Krzysztof",
title="Riemann--Liouville, Caputo, and Sequential Fractional Derivatives in Differential Games",
bookTitle="Advances in Dynamic Games: Theory, Applications, and Numerical Methods for Differential and Stochastic Games",
year="2011",
publisher="Birkh{\"a}user Boston",
address="Boston",
pages="61--81",
isbn="978-0-8176-8089-3",
doi="10.1007/978-0-8176-8089-3\_4",
url="https://doi.org/10.1007/978-0-8176-8089-3\_4"
}

@book{Dockner2000,
   author    = {Dockner, Engelbert J. and J{\o}rgensen, Steffen and Van Long, Ngo and Sorger, Gerhard},
  title     = {Differential Games in Economics and Management Science},
  year      = {2012},
  publisher = {Cambridge University Press},
  address   = {Cambridge},
  isbn      = {9780511805127},
  doi       = {10.1017/CBO9780511805127}
}

@book{Engwerda2005,
  author    = {Engwerda, Jacob C.},
  title     = {LQ Dynamic Optimization and Differential Games},
  year      = {2005},
  publisher = {John Wiley \& Sons},
  address   = {Chichester},
  edition   = {1},
  isbn      = {9780470015247}
}

@article{Achdou2013,
  author    = {Yves Achdou and Pierre Cardaliaguet and Fran{\c{c}}ois Delarue and Alessio Porretta},
  title     = {Mean Field Games: Numerical Methods for the Planning Problem},
  journal   = {Mathematical Models and Methods in Applied Sciences},
  year      = {2013},
  volume    = {22},
  number    = {1},
  pages     = {1250001},
  doi       = {10.1142/S0218202512500017}
}

@misc{Graber2025MFG,
  title={Mean Field Games of Controls with Fractional Laplacian}, 
      author={P. Jameson Graber and Elizabeth Matter and Jesus Ruiz Bolanos},
      year={2025},
      eprint={2509.04647},
      archivePrefix={arXiv},
      primaryClass={math.AP},
      url={https://arxiv.org/abs/2509.04647}, 
}

@article{Camilli2019CaputoHJ,
  author    = {Fabio Camilli and Raul De Maio and Carlo Marchi},
  title     = {Approximation of Hamilton--Jacobi Equations with Caputo Time-Fractional Derivative},
  journal   = {arXiv preprint arXiv:1906.06868},
  year      = {2019}
}

@article{Pozza_2025,
   title={Large time behavior of solutions to Hamilton–Jacobi equations on networks},
   volume={32},
   ISSN={1420-9004},
   url={http://dx.doi.org/10.1007/s00030-025-01128-5},
   DOI={10.1007/s00030-025-01128-5},
   number={6},
   journal={Nonlinear Differential Equations and Applications NoDEA},
   publisher={Springer Science and Business Media LLC},
   author={Pozza, Marco},
   year={2025},
   month=sep }

@article{Mamatov_2021,
doi = {10.1088/1742-6596/2068/1/012001},
url = {https://doi.org/10.1088/1742-6596/2068/1/012001},
year = {2021},
month = {oct},
publisher = {IOP Publishing},
volume = {2068},
number = {1},
pages = {012001},
author = {Mamatov, Mashrabjon and Alimov, Xakimjon},
title = {Controlled Fractional Order Processes with Distributed Parameters},
journal = {Journal of Physics: Conference Series}
}

@article{ChikriiMatychyn2012,
author = {Matychyn, Ivan and Chikrii, Arkadii and Onyshchenko, Viktoriia},
journal = {Mathematica Balkanica New Series},
language = {eng},
number = {1-2},
pages = {159-168},
publisher = {Bulgarian Academy of Sciences - National Committee for Mathematics},
title = {Conflict-Controlled Processes Involving Fractional Differential Equations with Impulses},
url = {http://eudml.org/doc/281434},
volume = {26},
year = {2012},
}

@book{AndersonMoore2007,
title={Optimal Control: Linear Quadratic Methods},
  author={Anderson, B.D.O. and Moore, J.B.},
  isbn={9780486457666},
  lccn={2006052128},
  series={Dover Books on Engineering},
  }

@article{FrancisWonham1976,
  author  = {Bruce A. Francis and W. Murray Wonham},
  title   = {The Internal Model Principle of Control Theory},
  journal = {Automatica},
  year    = {1976},
  volume  = {12},
  number  = {5},
  pages   = {457--465}
}

@article{Francis1977,
  author  = {Bruce A. Francis},
  title   = {The Linear Multivariable Regulator Problem},
  journal = {SIAM Journal on Control and Optimization},
  year    = {1977},
  volume  = {15},
  number  = {3},
  pages   = {486--505}
}

@book{RawlingsMayneDiehl2017,
     author    = {Rawlings, James B. and Mayne, David Q. and Diehl, Moritz M.},
  title     = {Model Predictive Control: Theory, Computation, and Design},
  publisher = {Nob Hill Publishing, LLC},
  address   = {Madison, WI},
  year      = {2024},
  isbn      = {978-0975937785}
}

@article{Mayne2000,
  author  = {David Q. Mayne and James B. Rawlings and Christopher V. Rao and Pierre O. M. Scokaert},
  title   = {Constrained Model Predictive Control: Stability and Optimality},
  journal = {Automatica},
  year    = {2000},
  volume  = {36},
  number  = {6},
  pages   = {789--814}
}

@inproceedings{ChenPetrasXue2009,
  author    = {YangQuan Chen and Ivo Petr{\'a}{\v s} and Dingy{\"u} Xue},
  title     = {Fractional Order Control--A Tutorial},
  booktitle = {Proceedings of the American Control Conference (ACC)},
  year      = {2009},
  pages     = {1397--1411}
}

@book{Tepljakov2017,
  author    = {Tepljakov, Aleksei},
  title     = {Fractional-Order Modeling and Control of Dynamic Systems},
  year      = {2017},
  publisher = {Springer},
  address   = {Cham},
  series    = {Springer Theses},
  edition   = {1},
  isbn      = {9783319529493},
  doi       = {10.1007/978-3-319-52950-9}
}

@book{Das2011,
  author    = {Das, Shantanu},
  title     = {Functional Fractional Calculus},
  year      = {2011},
  publisher = {Springer-Verlag},
  address   = {Berlin, Heidelberg},
  edition   = {2},
  isbn      = {9783642205446},
  doi       = {10.1007/978-3-642-20545-3}
}

@inproceedings{AxtellBise1990,
  author    = {Michael Axtell and Michael E. Bise},
  title     = {Fractional Calculus in Control Systems},
  booktitle = {Proceedings of the IEEE Southeastern Symposium on System Theory (SSST)},
  year      = {1990},
  pages     = {498--502}
}

@article{Sopasakis2016MPFOS,
  author  = {Pantelis Sopasakis and Haralambos Sarimveis},
  title   = {Stabilising Model Predictive Control for Discrete-Time Fractional-Order Systems},
  journal = {arXiv preprint arXiv:1605.08719},
  year    = {2016},
  url     = {https://arxiv.org/abs/1605.08719}
}

@article{Ntouskas2018OffsetFreeFO,
  author  = {Ntouskas, Sotiris and Sarimveis, Haralambos and Sopasakis, Pantelis},
  title   = {Model Predictive Control for Offset-Free Reference Tracking of Fractional Order Systems},
  journal = {Control Engineering Practice},
  year    = {2018},
  volume  = {71},
  pages   = {26--33},
  doi     = {10.1016/j.conengprac.2017.10.010}
}

@article{Ntouskas2019IFAC,
title = {Model predictive control for offset-free reference tracking of fractional order systems},
journal = {Control Engineering Practice},
volume = {71},
pages = {26-33},
year = {2018},
issn = {0967-0661},
doi = {https://doi.org/10.1016/j.conengprac.2017.10.010},
url = {https://www.sciencedirect.com/science/article/pii/S0967066117302435},
author = {Sotiris Ntouskas and Haralambos Sarimveis and Pantelis Sopasakis}
}

@proceedings{Romero2013MPE,
 author = {Romero, Miguel and de Madrid, A´ngel P. and Man˜oso, Carolina and Herna´ndez, Roberto},
    title = {Application of Generalized Predictive Control to a Fractional Order Plant},
    volume = {Volume 5: 6th International Conference on Multibody Systems, Nonlinear Dynamics, and Control, Parts A, B, and C},
    series = {International Design Engineering Technical Conferences and Computers and Information in Engineering Conference},
    pages = {1285-1292},
    year = {2007},
    month = {09},
    doi = {10.1115/DETC2007-34389},
    url = {https://doi.org/10.1115/DETC2007-34389},
    eprint = {https://asmedigitalcollection.asme.org/IDETC-CIE/proceedings-pdf/IDETC-CIE2007/4806X/1285/4566788/1285_1.pdf},
}

@inproceedings{Joshi2014FO,
  author    = {M. M. Joshi and others},
  title     = {Model Predictive Control for Fractional-Order System},
  booktitle = {Proc. 11th International Conference on Informatics in Control, Automation and Robotics (ICINCO/SIMULTECH)},
  year      = {2014},
  note      = {Short paper; fractional predictive control case studies}
}

@article{QinBadgwell2003,
  author  = {S. Joe Qin and Thomas A. Badgwell},
  title   = {A Survey of Industrial Model Predictive Control Technology},
  journal = {Control Engineering Practice},
  year    = {2003},
  volume  = {11},
  number  = {7},
  pages   = {733--764},
  doi     = {10.1016/S0967-0661(02)00186-7}
}

@article{YanSunZhang2017CICP,
  author  = {Yan, Yonggui and Sun, Zhi-Zhong and Zhang, Jiwei},
  title   = {Fast Evaluation of the Caputo Fractional Derivative and Its Applications to Fractional Diffusion Equations: A Second-Order Scheme},
  journal = {Communications in Computational Physics},
  year    = {2017},
  volume  = {22},
  number  = {4},
  pages   = {1028--1048},
  doi     = {10.4208/cicp.OA-2017-0019}
}

@article{CuestaLubichPalencia2006,
  author  = {Cuesta, E. and Lubich, C. and Palencia, C.},
  title   = {Convolution Quadrature Time Discretization of Fractional Diffusion-Wave Equations},
  journal = {Mathematics of Computation},
  year    = {2006},
  volume  = {75},
  number  = {254},
  pages   = {673--696},
  doi     = {10.1090/S0025-5718-06-01839-7}
}

@article{Schadle2006FOCQ,
  author  = {Sch{\"a}dle, A. and L{\'o}pez-Fern{\'a}ndez, M. and Lubich, C.},
  title   = {Fast and Oblivious Convolution Quadrature},
  journal = {SIAM Journal on Scientific Computing},
  year    = {2006},
  volume  = {28},
  number  = {2},
  pages   = {421--438},
  doi     = {10.1137/050628030}
}

@proceedings{JiangZhang2017SISC,
    author = {Romero, Miguel and de Madrid, A´ngel P. and Man˜oso, Carolina and Herna´ndez, Roberto},
    title = {Application of Generalized Predictive Control to a Fractional Order Plant},
    volume = {Volume 5: 6th International Conference on Multibody Systems, Nonlinear Dynamics, and Control, Parts A, B, and C},
    series = {International Design Engineering Technical Conferences and Computers and Information in Engineering Conference},
    pages = {1285-1292},
    year = {2007},
    month = {09},
    doi = {10.1115/DETC2007-34389},
    url = {https://doi.org/10.1115/DETC2007-34389},
    eprint = {https://asmedigitalcollection.asme.org/IDETC-CIE/proceedings-pdf/IDETC-CIE2007/4806X/1285/4566788/1285_1.pdf},
}

@article{deLevie1963,
  author   = {R. de Levie},
  title    = {On porous electrodes in electrolyte solutions---I},
  journal  = {Electrochimica Acta},
  year     = {1963},
  volume   = {8},
  pages    = {751--780}
}

@book{Lasia2014,
  author    = {Lasia, Andrzej},
  title     = {Electrochemical Impedance Spectroscopy and Its Applications},
  publisher = {Springer},
  address   = {New York, NY},
  year      = {2014},
  doi       = {10.1007/978-1-4614-8933-7},
  isbn      = {978-1-4614-8932-0}
}

@article{Zou2018,
  author  = {Chenglin Zou and Xuning Feng and Lingxi Li and Xiaosong Hu and K. J. Tseng and Minggao Ouyang},
  title   = {A review of fractional-order techniques applied to lithium-ion batteries},
  journal = {Journal of Power Sources},
  year    = {2018},
  volume  = {390},
  pages   = {286--296}
}

@article{HavriliakNegami1967,
  author  = {S. Havriliak and S. Negami},
  title   = {A complex plane analysis of $\alpha$-dispersions in some polymer systems},
  journal = {Journal of Polymer Science Part C: Polymer Symposia},
  year    = {1967},
  volume  = {14},
  number  = {1},
  pages   = {99--117}
}

@article{MakrisConstantinou1991,
author = {Makris, Nicos and Constantinou, Michael},
year = {1991},
month = {09},
pages = {},
title = {Fractional-Derivative Maxwell Model for Viscous Dampers},
volume = {117},
journal = {Journal of Structural Engineering-asce - J STRUCT ENG-ASCE},
doi = {10.1061/(ASCE)0733-9445(1991)117:9(2708)}
}

@article{RossikhinShitikova2010,
  author  = {Y. A. Rossikhin and M. V. Shitikova},
  title   = {Applications of fractional calculus to dynamic problems of solid mechanics: novel trends and recent results},
  journal = {Applied Mechanics Reviews},
  year    = {2010},
  volume  = {63},
  number  = {1},
  pages   = {010801}
}

@ARTICLE{Chopade2016,
  author={Chopade, Amit S. and Khubalkar, Swapnil W. and Junghare, A. S. and Aware, M. V. and Das, Shantanu},
  journal={IEEE/CAA Journal of Automatica Sinica}, 
  title={Design and implementation of digital fractional order PID controller using optimal pole-zero approximation method for magnetic levitation system}, 
  year={2018},
  volume={5},
  number={5},
  pages={977-989},
  doi={10.1109/JAS.2016.7510181}}

@misc{Sopasakis2019,
  author       = {Pantelis Sopasakis and Haralambos Sarimveis and Panos Macheras and Aristides Dokoumetzidis},
  title        = {Fractional Calculus in Pharmacokinetics},
  howpublished = {arXiv:1904.10556},
  year         = {2019}
}

@article{Giga2020Caputo,
  author = {Yoshikazu Giga and Tokinaga Namba},
title = {Well-posedness of Hamilton–Jacobi equations with Caputo’s time fractional derivative},
journal = {Communications in Partial Differential Equations},
volume = {42},
number = {7},
pages = {1088--1120},
year = {2017},
publisher = {Taylor \& Francis},
doi = {10.1080/03605302.2017.1324880}
}

@book{Meerschaert2011,
  author    = {Meerschaert, Mark M. and Sikorskii, Alla},
  title     = {Stochastic Models for Fractional Calculus},
  year      = {2019},
  edition   = {2},
  publisher = {De Gruyter},
  address   = {Berlin, Boston},
  series    = {De Gruyter Studies in Mathematics},
  volume    = {43},
  isbn      = {9783110559071}
}

@article{Raissi2019,
  author  = {Maziar Raissi and Paris Perdikaris and George E. Karniadakis},
  title   = {Physics-informed neural networks: A deep learning framework for solving forward and inverse problems involving nonlinear partial differential equations},
  journal = {Journal of Computational Physics},
  year    = {2019},
  volume  = {378},
  pages   = {686--707}
}

@inproceedings{LiFNO2021,
  author    = {Zongyi Li and Nikola Kovachki and Kamyar Azizzadenesheli and Burigede Liu and Andrew Stuart and Anima Anandkumar},
  title     = {Fourier Neural Operator for Parametric Partial Differential Equations},
  booktitle = {International Conference on Learning Representations (ICLR)},
  year      = {2021}
}

@book{Magni2008,
  author    = {Magni, Lalo and Raimondo, D. M. and Allg{\"o}wer, Frank},
  title     = {Nonlinear Model Predictive Control: Towards New Challenging Applications},
  subtitle  = {Selected Papers Based on the Presentations at the International Workshop on Assessment and Future Directions of Nonlinear Model Predictive Control (NMPC'08), Pavia, Italy, September 5--9, 2008},
  series    = {Lecture Notes in Control and Information Sciences},
  volume    = {384},
  publisher = {Springer},
  address   = {Berlin, Heidelberg},
  year      = {2009},
  doi       = {10.1007/978-3-642-01094-1},
  isbn      = {978-3-642-01093-4}
}

@article{GOMOYUNOV2024335,
title = {On viscosity solutions of path-dependent Hamilton–Jacobi–Bellman–Isaacs equations for fractional-order systems},
journal = {Journal of Differential Equations},
volume = {399},
pages = {335-362},
year = {2024},
issn = {0022-0396},
doi = {https://doi.org/10.1016/j.jde.2024.04.001},
url = {https://www.sciencedirect.com/science/article/pii/S0022039624002080},
author = {M.I. Gomoyunov}
}

@misc{gomoyunov2024lipschitzcontinuityresultsminimax,
      title={Lipschitz Continuity Results for Minimax Solutions of Path-Dependent Hamilton--Jacobi Equations}, 
      author={Mikhail I. Gomoyunov},
      year={2024},
      eprint={2412.17388},
      archivePrefix={arXiv},
      primaryClass={math.OC},
      url={https://arxiv.org/abs/2412.17388}, 
}

@article{CAO2024490,
title = {Will Fractional Order Model Based MPC Save Control Energy?},
journal = {IFAC-PapersOnLine},
volume = {58},
number = {12},
pages = {490-495},
year = {2024},
note = {12th IFAC Conference on Fractional Differentiation and its Applications ICFDA 2024},
issn = {2405-8963},
doi = {https://doi.org/10.1016/j.ifacol.2024.08.239},
url = {https://www.sciencedirect.com/science/article/pii/S2405896324009698},
author = {Shiang Cao and YangQuan Chen}
}

@article{deAzeredo2024IFAC,
  author  = {R. N. de Azeredo and M. S. Bruzon and L. E. V. Santo and G. F. de Souza and L. G. dos Santos and R. Sampaio and E. E. L. Filho},
  title   = {Comparison between {CRONE} and H-infinity control approaches for tracking antennas},
  journal = {IFAC-PapersOnLine},
  volume  = {58},
  number  = {12},
  pages   = {424--429},
  year    = {2024},
  doi     = {10.1016/j.ifacol.2024.08.228}
}

@article{Gomoyunov2025COCV,
  author  = {Mikhail I. Gomoyunov},
  title   = {Zero-sum games for Volterra integral equations and viscosity solutions of path-dependent Hamilton--Jacobi equations},
  journal = {ESAIM: Control, Optimisation and Calculus of Variations},
  volume  = {31},
  pages   = {55},
  year    = {2025},
  doi     = {10.1051/cocv/2025040}
}

@article{Guglielmi2025SOE,
  author  = {Nicola Guglielmi and Veljko Mili{\v{s}}i{\'c} and Manuela Pasotti},
  title   = {A fast and memoryless method for fractional order models},
  journal = {arXiv preprint arXiv:2509.17338},
  year    = {2025},
  eprint  = {2509.17338},
  archivePrefix = {arXiv},
  primaryClass  = {math.NA}
}

@article{Jayaram2024FOPID,
  author  = {S. Jayaram and N. Venkatesan},
  title   = {Design and implementation of the fractional-order controllers for a real-time nonlinear process using the {AGTM} optimization technique},
  journal = {Scientific Reports},
  volume  = {14},
  pages   = {29867},
  year    = {2024},
  doi     = {10.1038/s41598-024-82258-1}
}

@article{Mseddi2024FractalFract,
  author  = {Aymen Mseddi and Mohamed Houas and Pedro E. Miyagi and Habib Ben Attous},
  title   = {Investigation of the Robust Fractional Order Control Approach Associated with the Online Analytic Unity Magnitude Shaper: The Case of Wind Energy Systems},
  journal = {Fractal and Fractional},
  volume  = {8},
  number  = {4},
  pages   = {187},
  year    = {2024},
  doi     = {10.3390/fractalfract8040187}
}

@article{DEAZEREDO2024424,
title = {Comparison between CRONE and H-infinity control approaches applied to tracking antennas},
journal = {IFAC-PapersOnLine},
volume = {58},
number = {12},
pages = {424-429},
year = {2024},
note = {12th IFAC Conference on Fractional Differentiation and its Applications ICFDA 2024},
issn = {2405-8963},
doi = {https://doi.org/10.1016/j.ifacol.2024.08.228},
url = {https://www.sciencedirect.com/science/article/pii/S2405896324009571},
author = {Rodrigo Negri {de Azeredo} and Mohamed Hajjem and Lara Thomas and Stéphane Victor and Patrick Lanusse and Pierre Melchior}
}

@article{Lamrani2024TemperedArxiv,
  author       = {Ilyasse Lamrani and Hanaa Zitane and Delfim F. M. Torres},
  title        = {Controllability and observability of tempered fractional differential systems},
  journal      = {arXiv preprint arXiv:2412.05349},
  year         = {2024},
  url          = {https://arxiv.org/abs/2412.05349}
}

@article{Ineh2024AIMSmath,
  author       = {Matthew P. Ineh and S. D. Oke and M. A. Momoh},
  title        = {A novel approach to Lyapunov stability of Caputo fractional dynamic equations on time scale using a new generalized derivative},
  journal      = {AIMS Mathematics},
  year         = {2024},
  volume       = {9},
  number       = {16},
  pages        = {1639--1663},
  doi          = {10.3934/math.20241639}
}

@article{Vishnukumar2024PhysScr,
  author       = {K. S. Vishnukumar and others},
  title        = {Controllability of time-varying fractional dynamical systems in the Caputo sense},
  journal      = {Physica Scripta},
  year         = {2024},
  volume       = {99},
  number       = {6},
  pages        = {065218}
}

@article{Ferreira2024HigherCaputo,
  author   = {Rui A. C. Ferreira},
  title    = {Calculus of Variations with Higher Order Caputo Fractional Derivatives},
  journal  = {Arabian Journal of Mathematics},
  year     = {2024},
  volume   = {13},
  pages    = {91--101},
  doi      = {10.1007/s40065-023-00447-8},
  note     = {Issue date: April 2024}
}

@misc{Torres2024Duality,
  author       = {Delfim F. M. Torres},
  title        = {The Duality Theory of Fractional Calculus},
  howpublished = {arXiv:2404.14458},
  year         = {2024},
  url          = {https://arxiv.org/abs/2404.14458}
}

@article{Gasimov2024QTDS,
  author  = {J. J. Gasimov},
  title   = {Pontryagin Maximum Principle for Fractional Delay Differential Equations and Volterra Integral Equations},
  journal = {Qualitative Theory of Dynamical Systems},
  year    = {2024},
  doi     = {10.1007/s12346-024-01049-1}
}

@article{Moon2025AIMSmath,
  author  = {Jiyoung Moon and Yongjin Choi},
  title   = {The Pontryagin type maximum principle for Caputo fractional optimal control problems with state constraints},
  journal = {AIMS Mathematics},
  year    = {2025},
  doi     = {10.3934/math.2025042}
}

@article{NdairouTorres2023Math,
  author  = {Fa{\"i}{\c{c}}al Nda{\"\i}rou and Delfim F. M. Torres},
  title   = {Pontryagin Maximum Principle for Incommensurate Fractional-Orders Optimal Control Problems},
  journal = {Mathematics},
  year    = {2023},
  volume  = {11},
  number  = {19},
  pages   = {4218},
  doi     = {10.3390/math11194218}
}

@article{Cruz2025AIMSmath,
  author  = {F. Cruz and D. F. M. Torres},
  title   = {A Pontryagin maximum principle for optimal control problems with distributed-order Caputo derivative},
  journal = {AIMS Mathematics},
  year    = {2025},
  doi     = {10.3934/math.2025539}
}

@article{YusubovMahmudov2024ACS,
  author  = {Shakir Sh. Yusubov and Elimhan N. Mahmudov},
  title   = {Pontryagin’s maximum principle for the Roesser model with a fractional Caputo derivative},
  journal = {Archives of Control Sciences},
  year    = {2024},
  volume  = {34},
  number  = {2},
  pages   = {271--300},
  doi     = {10.24425/acs.2024.149661}
}

@article{Zhu2019SOE,
  author  = {Zhu, Hongyi and Xu, Chuanju},
  title   = {A Fast High Order Method for the Time-Fractional Diffusion Equation},
  journal = {SIAM Journal on Numerical Analysis},
  year    = {2019},
  volume  = {57},
  number  = {6},
  pages   = {2829--2849},
  doi     = {10.1137/18M1231225}
}

@article{Gomoyunov2024MCRF,
 title={Value functional and optimal feedback control in linear-quadratic optimal control problem for fractional-order system},
   volume={14},
   ISSN={2156-8499},
   url={http://dx.doi.org/10.3934/mcrf.2023002},
   DOI={10.3934/mcrf.2023002},
   number={1},
   journal={Mathematical Control and Related Fields},
   publisher={American Institute of Mathematical Sciences (AIMS)},
   author={Gomoyunov, Mikhail I.},
   year={2024},
   pages={215–254} }

@article{Yoneda2024FractalFract,
  author  = {Ryo Yoneda and Yuki Moriguchi and Masaharu Kuroda and Natsuki Kawaguchi},
  title   = {Servo Control of a Current-Controlled Attractive-Force-Type Magnetic Levitation System Using Fractional-Order LQR Control},
  journal = {Fractal and Fractional},
  year    = {2024},
  volume  = {8},
  number  = {8},
  pages   = {458},
  doi     = {10.3390/fractalfract8080458}
}

@article{Malmir2024ISA,
  author  = {Iman Malmir},
  title   = {New pure multi-order fractional optimal control problems modeled as quadratic programming},
  journal = {ISA Transactions},
  year    = {2024},
  note    = {in press},
  doi     = {10.1016/j.isatra.2024.08.031}
}

@article{Leugering2024FractFract,
  author  = {G{\"u}nter Leugering and Abdelsalam Alqudah and Yazeed N. Alhindi and Abdelhamid Elaifi and Munahi Alkinani},
  title   = {Non-Overlapping Domain Decomposition for 1D Optimal Control Problems Governed by Caputo--Fractional ODEs},
  journal = {Fractal and Fractional},
  year    = {2024},
  volume  = {8},
  number  = {3},
  pages   = {129},
  doi     = {10.3390/fractalfract8030129}
}

@article{Yangari2024ADE,
  author  = {Miguel Yangari},
  title   = {Some regularity results of Hamilton--Jacobi equations with variable-order Caputo time derivative},
  journal = {Advances in Differential Equations},
  year    = {2024},
  volume  = {29},
  number  = {11/12},
  pages   = {951--980}
}

@article{ZhangFangSun2021ESA,
  author  = {Jia{-}Li Zhang and Zhi{-}Wei Fang and Hai{-}Wei Sun},
  title   = {Exponential-sum-approximation technique for variable-order time-fractional diffusion equations},
  journal = {arXiv preprint arXiv:2101.08125},
  year    = {2021}
}

@article{CamilliDuisembay2021JDG,
  author  = {Fabio Camilli and Serikbolsyn Duisembay},
  title   = {Approximation of an optimal control problem for the time-fractional Fokker--Planck equation},
  journal = {Journal of Dynamics and Games},
  year    = {2020},
  volume  = {8},
  number  = {4},
  pages   = {343--365},
  doi     = {10.3934/jdg.2021013}
}

@article{Wendimu2025SciRep,
  author    = {Abebe A. Wendimu and Mohamed Elsisi and Elisabetta Chicca and S{\'a}ndor Erdei and Szabolcs Nagyne Simon and Istv{\'a}n Harmati},
  title     = {Fractional-order PID control for elevation and azimuth in a twin rotor system},
  journal   = {Scientific Reports},
  year      = {2025},
  volume    = {15},
  number    = {},
  pages     = {e18763},
  doi       = {10.1038/s41598-025-18763-8},
  url       = {https://www.nature.com/articles/s41598-025-18763-8}
}

@article{Krzysztofik2025ApplSci,
  author    = {Igor Krzysztofik and S{\l}awomir Blasiak},
  title     = {The Fractional Order PID Controller for the 3-DOF Gyroscope System},
  journal   = {Applied Sciences},
  year      = {2025},
  volume    = {15},
  number    = {19},
  pages     = {10476},
  doi       = {10.3390/app151910476},
  url       = {https://www.mdpi.com/2076-3417/15/19/10476}
}

@article{Mihaly2024FMI,
  author    = {V. Mihaly and colleagues},
  title     = {Robust numeric implementation of the fractional-order controller},
  journal   = {Flow Measurement and Instrumentation},
  year      = {2024},
  note      = {In press},
  url       = {https://www.sciencedirect.com/science/article/abs/pii/S0016003224005088}
}

@article{Cao2024IFAC,
  author  = {Shuwei Cao and others},
  title   = {Will Fractional Order Model Based MPC Save Control?},
  journal = {IFAC-PapersOnLine},
  year    = {2024},
  note    = {Proc. IFAC},
  url     = {https://www.sciencedirect.com/science/article/pii/S2405896324009698}
}

@article{Das2025IFAC,
  author  = {D. Das and others},
  title   = {Fractional Model Predictive Control of Fractional Chaotic Systems},
  journal = {IFAC-PapersOnLine},
  year    = {2025},
  url     = {https://www.sciencedirect.com/science/article/pii/S2405896325013047}
}

@article{GuglielmiHairer2025,
  author  = {Nicola Guglielmi and Ernst Hairer},
  title   = {A fast and memoryless numerical method for solving fractional differential equations},
  journal = {arXiv preprint arXiv:2506.04188},
  year    = {2025},
  url     = {https://arxiv.org/abs/2506.04188}
}

@article{Peng2025SciRep,
  author  = {Chao Peng and Seyyed Morteza Ghamari and Hasan Mollaee and Omid Rezaei},
  title   = {Design of a novel robust adaptive fractional-order model predictive controller for boost converter using grey wolf optimization algorithm},
  journal = {Scientific Reports},
  year    = {2025},
  volume  = {15},
  number  = {27670},
  doi     = {10.1038/s41598-025-10125-8},
  url     = {https://www.nature.com/articles/s41598-025-10125-8}
}

@article{ZhengWang2024JSC,
  author  = {Zi-Yun Zheng and Yuan-Ming Wang},
  title   = {Fast High-Order Compact Finite Difference Methods Based on the Averaged L1 Formula for a Time-Fractional Mobile–Immobile Diffusion Problem},
  journal = {Journal of Scientific Computing},
  year    = {2024},
  volume  = {99},
  number  = {2},
  doi     = {10.1007/s10915-024-02633-1}
}

@article{GanZhangLiang2024JAMC,
  author  = {Di Gan and Guo-Feng Zhang and Zhao-Zheng Liang},
  title   = {An Efficient Preconditioner for Linear Systems Arising from High-Order Accurate Schemes of Time-Fractional Diffusion Equations},
  journal = {Journal of Applied Mathematics and Computing},
  year    = {2024}
}

@article{HuangLuoZhang2025,
  author  = {Daoming Huang and Zhiwen Luo and Zhi Zhang},
  title   = {Pointwise Error Analysis of Corrected L1 Scheme for Subdiffusion},
  journal = {Fractal and Fractional},
  year    = {2025}
}

@article{HeLiWei2025,
  author  = {Saifei He and Xiquan Li and Guannan Wei},
  title   = {A Parareal Algorithm for Caputo--Hadamard Fractional Differential Equations},
  journal = {Open Mathematics},
  year    = {2025}
}

@article{Mojahed2020PFDA,
  author  = {Allahbakhsh Yazdani and Navid Mojahed and Afshin Babaei and Elena V{\'a}zquez Cend{\'o}n},
  title   = {Using Finite Volume-Element Method for Solving Space Fractional Advection--Dispersion Equation},
  journal = {Progress in Fractional Differentiation and Applications},
  year    = {2020},
  volume  = {6},
  pages   = {55--66},
  doi     = {10.18576/pfda/060106}
}

@article{ZhangYang2024FractalFract,
  author  = {J. Zhang and Q. Yang},
  title   = {The Finite Volume Element Method for Time Fractional Generalized Burgers' Equation},
  journal = {Fractal and Fractional},
  year    = {2024},
  volume  = {8},
  pages   = {53},
  doi     = {10.3390/fractalfract8010053}
}

@article{Calatayud2024CSF,
  author  = {J. Calatayud and coauthors},
  title   = {Interpretation of Caputo Models: Where is the Memory?},
  journal = {Chaos, Solitons \& Fractals},
  year    = {2024}
}

@article{DaftardarGejjiJafari2006_JMAA,
  author  = {Varsha Daftardar-Gejji and Hossein Jafari},
  title   = {An iterative method for solving nonlinear functional equations},
  journal = {Journal of Mathematical Analysis and Applications},
  year    = {2006},
  volume  = {316},
  number  = {2},
  pages   = {753--763}
}

@article{Jafari2022_JVC_FOOC,
  author  = {Hossein Jafari and R. M. Ganji and K. Sayevand and Dumitru Baleanu},
  title   = {A numerical approach for fractional optimal control problems with {Mittag--Leffler} kernel},
  journal = {Journal of Vibration and Control},
  year    = {2022},
  volume  = {28},
  number  = {19-20},
  pages   = {2596--2606}
}

@article{FatoorehchiRach2020_AEJ,
  author  = {Hooman Fatoorehchi and Randolph Rach},
  title   = {A method for inverting the Laplace transforms of two classes of rational transfer functions in control engineering},
  journal = {Alexandria Engineering Journal},
  year    = {2020},
  volume  = {59},
  number  = {6},
  pages   = {4879--4887},
  doi     = {10.1016/j.aej.2020.08.052}
}

@Article{Fatoorehchi2023sign,
author={Fatoorehchi, Hooman
and Djilali, Salih},
title={Stability analysis of linear time-invariant dynamic systems using the matrix sign function and the {A}domian decomposition method},
journal={International Journal of Dynamics and Control},
year={2023},
month={Apr},
day={01},
volume={11},
number={2},
pages={593-604},
issn={2195-2698},
doi={10.1007/s40435-022-00989-3},
url={https://doi.org/10.1007/s40435-022-00989-3}
}

@article{Obeidat21062024,
author = {Nazek A. Obeidat and Mahmoud S. Rawashdeh and Malak Q. Al Erjani},
title = {A novel {A}domian natural decomposition method with convergence analysis of nonlinear time-fractional differential equations},
journal = {International Journal of Modelling and Simulation},
volume = {0},
number = {0},
pages = {1--16},
year = {2024},
publisher = {Taylor \& Francis},
doi = {10.1080/02286203.2024.2369772}
}

@Article{Ziada2024,
author={Ziada, Eman A. A.},
title={Picard and {A}domian solutions of nonlinear fractional differential equations system containing Atangana -- Baleanu derivative},
journal={Journal of Engineering and Applied Science},
year={2024},
month={Feb},
day={02},
volume={71},
number={1},
pages={31},
issn={2536-9512},
doi={10.1186/s44147-024-00361-6},
url={https://doi.org/10.1186/s44147-024-00361-6}
}

@Article{Sebaq2025,
author={Sebaq, Mahmoud S.
and Qamlo, Ahlam H.
and Bahaa, G. M.},
title={Numerical solutions for fractional optimal control problems of coupled diffusion systems via {L}aplace {A}domian Decomposition Method},
journal={Boundary Value Problems},
year={2025},
month={Sep},
day={01},
volume={2025},
number={1},
pages={131},
issn={1687-2770},
doi={10.1186/s13661-025-02108-5},
url={https://doi.org/10.1186/s13661-025-02108-5}
}

@article{Bao2018MicroNanoscaleHeat,
  author  = {Bao, Hua and Chen, Jie and Gu, Xiaokun and Cao, Bingyang},
  title   = {A Review of Simulation Methods in Micro/Nanoscale Heat Conduction},
  journal = {ES Energy and Environment},
  year    = {2018},
  volume  = {1},
  number  = {1},
  pages   = {16--55},
  doi     = {10.30919/esee8c149}
}

@book{Zhmakin2023NonFourierBook,
  author    = {Zhmakin, Alexander I.},
  title     = {Non-Fourier Heat Conduction},
  publisher = {Springer},
  address   = {Cham},
  year      = {2023},
  doi       = {10.1007/978-3-031-25973-9}
}

@article{Lukashchuk2024TimeFracDPL,
  author  = {Lukashchuk, Stanislav Yu.},
  title   = {A Semi-Explicit Algorithm for Parameters Estimation in a Time-Fractional Dual-Phase-Lag Heat Conduction Model},
  journal = {Modelling},
  year    = {2024},
  volume  = {5},
  number  = {3},
  pages   = {776--796},
  doi     = {10.3390/modelling5030041}
}

@article{Povstenko2025FracHeatContact,
  author  = {Povstenko, Yuriy and Kyrylych, Tamara and Dashkiiev, Viktor and Yatsko, Andrzej},
  title   = {Fundamental Solutions to Fractional Heat Conduction in Two Joint Half-Lines Under Conditions of Nonperfect Thermal Contact},
  journal = {Entropy},
  year    = {2025},
  volume  = {27},
  number  = {9},
  pages   = {965},
  doi     = {10.3390/e27090965}
}

@article{Rashid2024LettersFins,
  author  = {Rashid, Farhan Lafta and Dhaidan, Nabeel S. and Mahdi, Ali Jafer},
  title   = {Heat Transfer Enhancement of Phase Change Materials Using Letters-Shaped Fins: A Review},
  journal = {International Communications in Heat and Mass Transfer},
  year    = {2024},
  volume  = {159},
  pages   = {108096},
  doi     = {10.1016/j.icheatmasstransfer.2024.108096}
}

@article{Shmyhol2024PCMHybridTanks,
  author  = {Shmyhol, Dmytro and Rim{\'a}r, Miroslav and Fedak, Marcel and Krenick{\'y}, Tibor and Lopu{\v{s}}niak, Martin and Polivka, Nikolas},
  title   = {Techniques for Enhancing Thermal Conductivity and Heat Transfer in Phase Change Materials in Hybrid Phase Change Material--Water Storage Tanks},
  journal = {Applied Sciences},
  year    = {2024},
  volume  = {14},
  number  = {9},
  pages   = {3732},
  doi     = {10.3390/app14093732}
}
%% if required, the content of .bbl file can be included here once bbl is generated
%%\input sn-article.bbl

\end{document}